\documentclass[12pt, preprint]{elsarticle}
\usepackage{amssymb}
\usepackage{amsfonts}
\usepackage{amsmath}
\usepackage[dvips,bottom=1.4in,right=1in,top=1in, left=1in,nohead]{geometry}

\newtheorem{theorem}{Theorem}[section]
\newtheorem{acknowledgement}[theorem]{Acknowledgement}

\newtheorem{corollary}[theorem]{Corollary}

\newtheorem{definition}[theorem]{Definition}

\newtheorem{lemma}[theorem]{Lemma}

\newtheorem{proposition}[theorem]{Proposition}
\newtheorem{remark}[theorem]{Remark}

\newenvironment{proof}[1][Proof]{\noindent\textbf{#1.} }{\ \rule{0.5em}{0.5em}}
\numberwithin{equation}{section}

\begin{document}

\title{\textbf{On a class of degenerate parabolic equations with dynamic
boundary conditions}}
\author{Ciprian G. Gal \\
Department of Mathematics\\
Florida International University,\\
Miami, FL 33199, USA\\
cgal@fiu.edu}

\begin{abstract}
\noindent {\footnotesize We consider a quasi-linear parabolic (possibly,
degenerate) equation with nonlinear dynamic boundary conditions. The
corresponding class of initial and boundary value problems has already been
studied previously, proving well-posedness of weak solutions and the
existence of the global attractor, assuming that the nonlinearities are
\emph{subcritical }to a given exponent. The goal of this article is to show
that the previous analysis can be redone for \emph{supercritical}
nonlinearities by proving an additional }$L^{\infty }${\footnotesize %
-estimate on the solutions. In particular, we derive new conditions which
reflect an exact balance between the internal and the boundary mechanisms
involved, even when both the nonlinear sources contribute in opposite
directions. Then, we show how to construct a trajectory attractor for the
weak solutions of the associated parabolic system, and prove that any
solution belonging to the attractor is bounded, which implies uniqueness.
Finally, we also prove for the (semilinear) reaction-diffusion equation with
nonlinear dynamic boundary conditions, that the fractal dimension of the
global attractor is of the order }$\nu ^{-\left( N-1\right) },$%
{\footnotesize \ as diffusion }$\nu \rightarrow 0^{+},${\footnotesize \ in
any space dimension }$N\geq 2${\footnotesize , improving some recent results
in} {\footnotesize \cite{GalNS}}.{\footnotesize \ }
\end{abstract}

\maketitle

\section{Introduction}

Let us consider the following partial differential equation%
\begin{equation}
\partial _{t}u-\text{div}\left( a\left( \left\vert \nabla u\right\vert
^{2}\right) \nabla u\right) +f\left( u\right) =h_{1}(x),  \label{1.11}
\end{equation}%
in $\Omega \times (0,+\infty ),$ where $\Omega $ is a bounded domain in $%
\mathbb{R}^{N},$ $N\geq 1$, with smooth boundary $\Gamma :=\partial \Omega $%
, $a$ is a given nonnegative function, and $f$ and $h_{1}$ are suitable
functions. The mathematical literature regarding equation (\ref{1.11})
subject to all kinds of homogeneous boundary conditions is fairly vast. We
recall that global well-posedness results for (\ref{1.11}) with Dirichlet or
Neumann type of boundary conditions can be found in \cite{AEY, BV, CD, CF}
(see also \cite{Be, Br, CG, MO, TY, YSZ}). In addition, the analysis of
dissipative dynamical systems generated by equations like (\ref{1.11}) was
carried out in a number of papers mainly devoted to the asymptotic behavior
of strong solutions \cite{Ao, CF, Il, Op}, and to establish the existence of
global and/or exponential attractors (see, for instance, \cite{AEY, CG, EO,
MO, TY, YSZ}). For other classical results concerning the long term dynamics
of (\ref{1.11}) we also refer the reader to \cite{BV, EO}.

All the mentioned results are mainly concerned with standard boundary
conditions (that is, Dirichlet's and/or Neumann's). Let us now consider
\emph{dynamic} boundary conditions of the form
\begin{equation}
\partial _{t}u+b\left( x\right) a\left( \left\vert \nabla u\right\vert
^{2}\right) \partial _{\mathbf{n}}u+g\left( u\right) =h_{2}(x),  \label{1.12}
\end{equation}%
on $\Gamma \times (0,+\infty ),$ where $g$ and\ $h_{2}$ are suitable
functions defined on $\Gamma $, and $b\in L^{\infty }\left( \Gamma \right) ,$
$b\geq b_{0}>0$. This type of boundary conditions{\normalsize \ }arises for
many known equations of mathematical physics. They are motivated by problems
in diffusion phenomena \cite{Crank, CS, FGGR1, FGGR2, FZ, GR, Wa, Wu},
reaction-diffusion systems in phase-transition phenomena \cite{CGGM, GG1,
GG2, S}, special flows in hydrodynamics \cite{Langer, MW}, models in
climatology \cite{MW2}, and many others. For possible physical
interpretations of (\ref{1.12}) for problem (\ref{1.11}), we refer the
reader to \cite{Gal2} (cf. \cite{GR} also).

Problems such as \eqref{1.11}-\eqref{1.12} have already been investigated in
a number of papers \cite{CJ, Escher, Escher2, GW, Ma, Wa}. Constantin and
Escher deal with non-degenerate boundary value problems with smooth
nonlinearities (in particular, $\left\langle a\left( \left\vert \xi
\right\vert ^{2}\right) \xi ,\xi \right\rangle \geq c\left\vert \xi
\right\vert ^{2},$ with $c>0$) and show that unique (classical) maximal
solutions exist in some Bessel potential spaces \cite{Escher, Escher2}. Such
results enable the authors to investigate other qualitative properties
concerning global existence and blow-up phenomena (see, also \cite{CJ}).
These results are also improved by Meyries \cite{Ma}, still in the
non-degenerate case, by assuming more general boundary conditions and by
requiring that $f\left( s\right) /s$ and $g\left( s\right) /s$ are
dissipative as $\left\vert s\right\vert \rightarrow \infty $. A first
analysis, which aimed at deducing only a minimal number of assumptions on
the data and nonlinearities, was done in \cite{GW} by assuming that $f$ and $%
g$ are subcritical polynomial nonlinearities and by allowing $a\left(
s\right) $ to have a polynomial degeneracy at zero. For instance, one can
take%
\begin{equation}
a\left( s\right) =\left\vert s\right\vert ^{\left( p-2\right) /2},\text{ for
}p\neq 2.  \label{a}
\end{equation}%
In particular, we proved that problem (\ref{1.11})-(\ref{1.12}) with $%
a\left( s\right) $ as in (\ref{a}), subject to square-integrable initial
data $u_{\mid t=0}=u_{0}$ is well-posed,\ and then we established the
existence of a global attractor bounded in $W^{1,p}\left( \Omega \right) $.
Well-posedness for problem \eqref{1.11}, \eqref{1.12} for $a\left( s\right)
=\left\vert s\right\vert ^{\left( p-2\right) /2}$, assuming monotone
functions $f,$ $g$ was considered in \cite{Wa}. The non-degenerate case $%
a\left( s\right) \equiv \nu >0$ when $g=0,$ is discussed in detail in \cite%
{GalNS}. The stationary case associated with (\ref{1.11})-(\ref{1.12}) is
treated in \cite{GWs}.

It is well-known that when at least one of the source terms, the bulk
nonlinear term $f$ or the boundary term $g$ is present in (\ref{1.11})-(\ref%
{1.12}), conditions can be derived on their growth rates which imply either
the global existence of solutions or blow-up in finite time \cite{FQ}.
Namely in the non-degenerate case, for $\lambda ,\mu \in \left\{ 0,\pm
1\right\} $ with $\max \left\{ \lambda ,\mu \right\} =1$, $f\left( s\right)
:=-\lambda \left\vert s\right\vert ^{r_{1}-1}s$ and $g\left( s\right) :=-\mu
\left\vert s\right\vert ^{r_{2}-1}s$, solutions of%
\begin{equation}
\partial _{t}u-\nu \Delta u+f\left( u\right) =h_{1}\left( x\right) ,\text{
in }\Omega \times (0,+\infty ),  \label{1.11bb}
\end{equation}%
subject to the dynamic condition%
\begin{equation}
\partial _{t}u+\nu b\partial _{\mathbf{n}}u+g\left( u\right) =h_{2}\left(
x\right) ,\text{ on }\Gamma \times \left( 0,\infty \right) ,  \label{1.12bb}
\end{equation}%
are globally well-defined, for every given (sufficiently smooth) initial
data $u_{\mid t=0}=u_{0},$\ if $r_{1}r_{2}>1$ and $\lambda r_{1}+\mu r_{2}>0$%
. Furthermore, \cite{FQ} shows that if we further restrict the growths of $%
r_{1}$,$r_{2}$ so that $r_{1}<\left( N+2\right) /\left( N-2\right) $ and $%
r_{2}\leq N/\left( N-2\right) $, then the global solutions are also bounded.
On the other hand, if $\lambda =0$, $\mu =1,$ then some solutions blowup in
finite time with blowup occurring in the $L^{\infty }$-norm at a rate $%
\left( t-T_{\ast }\right) ^{-\left( r_{2}-1\right) },$ for some additional
conditions on $u_{0}$ and $r_{2}$. In the same way, when $\mu =0$ and $%
\lambda =1$, then some solutions blowup in finite time with a blowup rate
which depends on $r_{1}$ and $u_{0}$ (see \cite{Be2}). The occurrence of
blow up phenomena is closely related to the blowup problem for the ordinary
differential equation
\begin{equation}
u_{t}+h\left( u\right) =0,  \label{ODE}
\end{equation}%
where either $h=f$ or $h=g$. More precisely, it is easy to see that
solutions of the ODE (\ref{ODE}) are spatially homogeneous solutions of
either equation (\ref{1.11bb}) or equation (\ref{1.12bb}), and so if these
solutions blowup in finite time so do the solutions of (\ref{1.11bb}), (\ref%
{1.12bb}) (see \cite{RB} for further details, and additional references).
Similar results showing the same behavior are also derived for the parabolic
system (\ref{1.11})-(\ref{1.12}) in \cite{V}, stating sufficient conditions
for the functions $a$, $f$ and $g$ so that blowup in finite time occurs\ in
the $L^{\infty }$-norm. In particular, it was shown, for odd functions $f,$ $%
g$ and initial data $u_{0}\in L^{\infty }\left( \Omega \right) \cap
W^{1,p}\left( \Omega \right) $ and some additional conditions on $u_{0}$,
that there are solutions that blow up in finite time with an upper bound on
the blowup time which can be determined precisely.

The main goal of this paper is to deduce more general conditions (when
compared to conditions deduced in \cite{CJ, Escher, Escher2, FZ, GW, Ma, Wa,
Y0, YY}) on the reactive and radiation terms $f$ and $g$, respectively,
which imply that problem (\ref{1.11})-(\ref{1.12}) is dissipative in a
suitable sense, and that it possesses a (possibly, finite dimensional)
global attractor which characterizes the long-term behavior of the parabolic
system under consideration. Recently in \cite{RB} (see also \cite{YY2} for
some extensions), the authors have considered the semilinear parabolic
equation (\ref{1.11bb}) subject to nonlinear Robin boundary conditions%
\begin{equation}
\nu \partial _{\mathbf{n}}u+g\left( u\right) =0\text{ on }\Omega \times
(0,+\infty ),  \label{p2}
\end{equation}%
and they derived sufficient conditions on $f$ and $g,$ which imply
dissipativity for such problems. In particular, they have obtained a general
balance between $f$ and $g$, allowing for a real competition between both
the two nonlinear mechanisms which may work in opposite directions, one
fighting for blow-up in finite time, the other for dissipativity. Then, they
also proved the existence of a compact attractor in $H^{1}\left( \Omega
\right) ,$ assuming that the growth of $f$ and $g$ is subcritical. Their
method relies essentially on the fact that problem (\ref{1.11bb}), (\ref{p2}%
) possesses a Lyapunov functional, which can then be used to show either
dissipativity of (\ref{1.11bb}), (\ref{p2}), by exploiting some Poincare
type inequality (see (\ref{PI}) below), or blow up of some solutions.

Our goal is to extend these results in several directions, by working
instead with a class of degenerate parabolic equations, such as (\ref{1.11}%
), and then by subjecting (\ref{1.11}) to dynamic boundary conditions of the
form (\ref{1.12}). Moreover, we also wish to consider nonlinearities with
\emph{arbitrary} polynomial growth at infinity. We aim to construct weak
(energy)\ solutions with the help from a different (than in \cite{GW})
approximation scheme, which is based on the existence of classical (smooth)
solutions for a (strictly) non-degenerate system associated with (\ref{1.11}%
)-(\ref{1.12}). Let $F$ and $G$ be the primitives of $f$ and $g$,
respectively, such that $F\left( 0\right) =0,$ $G\left( 0\right) =0$. Even
though a natural energy functional exists for suitable approximates of the
problem (\ref{1.11})-(\ref{1.12}), at the moment it is not clear how to
prove that this energy (see \cite[(1.5)]{GW}), namely,

\begin{equation*}
\mathcal{E}_{\Omega ,\Gamma }\left( u\right) :=\int\limits_{\Omega }\left[
a\left( \left\vert \nabla u\right\vert ^{2}\right) \left\vert \nabla
u\right\vert ^{2}+F\left( u\right) -h_{1}\left( x\right) u\right]
dx+\int\limits_{\Gamma }\left[ G\left( u\right) -h_{2}\left( x\right) u%
\right] \frac{dS}{b},
\end{equation*}%
is in fact a Lyapunov function for (\ref{1.11})-(\ref{1.12}) when $a\left(
s\right) =\left\vert s\right\vert ^{\left( p-2\right) /2}$, $p\neq 2$, due
to a lack of regularity of the weak solutions (see, however \cite{GalNS},
and its references, when $a\left( s\right) =\nu >0$). Therefore, the method
in \cite{RB} which relies on the use of a Lyapunov function does not seem
applicable to our situation here. Indeed, when one is dealing with gradient
systems with a set of equilibria which is bounded in the phase space where $%
\mathcal{E}_{\Omega ,\Gamma }\left( u_{0}\right) <\infty $, one could avoid
to prove the existence of a bounded absorbing set and directly show the
existence of the global attractor for \emph{subcritical} nonlinearities.
However, since we wish to construct global attractors for (\ref{1.11})-(\ref%
{1.12}), under no essential growth assumptions on the nonlinearities, we
prefer to prove the existence of a uniform dissipative estimate which can be
also easily adapted to nonautonomous generalizations. Another difficult step
that we need to overcome is the uniqueness problem for the weak (energy)
solutions considered here. Indeed, we wish to deduce sufficiently general
conditions on $f$ and $g$ without excluding the scenario based on which
these functions are competing at infinity. When $f$ and $g$ are both
dissipative, i.e., when $f$ and $g$ are both polynomial potentials of odd
degree with a positive leading coefficient, uniqueness holds, and thus, the
dynamical system associated with the parabolic system (\ref{1.11})-(\ref%
{1.12}) can be defined in the classical sense. However, if at least one of
the polynomial nonlinearities possesses a negative leading coefficient (for
instance, suppose that%
\begin{equation}
g\left( s\right) \sim c_{g}\left\vert s\right\vert ^{r_{2}-1}s,\text{ as }%
\left\vert s\right\vert \rightarrow \infty ,  \label{ggg}
\end{equation}%
for $c_{g}<0$), uniqueness is not know, and in this case the classical
semigroup can be defined as a semigroup of multi-valued maps only. In order
to investigate the long-term behavior of the degenerate parabolic system (%
\ref{1.11})-(\ref{1.12}), we will employ the trajectory dynamical approach,
which allows us to avoid the use of unfriendly multivalued maps, and to
apply the usual theory of global attractors (see, e.g., \cite{CV} for the
general theory). We strongly emphasize that non-uniqueness of the weak
solutions constructed here is only a feature of the nonlinear interplay
between the two nonlinear mechanisms, and is not related to the smoothness
of the functions involved, as this is usually the case for other PDE's. See,
e.g., \cite{CV, MZ, T} and references therein, for applications\ for which
the uniqueness is not yet solved, such as, hyperbolic equations with
supercritical nonlinearities, reaction-diffusion systems, and so on. In
fact, even when our functions $f\left( s\right) $ and $g\left( s\right) $
are assumed to be (locally)\ Lipschitz, uniqueness of the weak energy
solutions is not known unless $f\left( s\right) $ and $g\left( s\right) $
are monotone increasing for all $\left\vert s\right\vert \geq s_{0}$ (for
some $s_{0}>0$), which is quite restrictive. Finally, to make matters worse,
for boundary nonlinearities that satisfy (\ref{ggg}), there may be solutions
of (\ref{1.11})-(\ref{1.12}) which blowup in finite time at some points in $%
\overline{\Omega },$ unless the internal mechanism governed by nonlinear
flux and reaction is sufficiently strong to overcome the boundary reaction.
Therefore, it is also essential to deduce some kind of \emph{optimal}
general conditions between the bulk and boundary nonlinearities that will
only produce bounded (non-singular)\ solutions for all positive times and
arbitrary initial data. The main difficulty here is, of course, to establish
the asymptotic compactness for the system (\ref{1.11})-(\ref{1.12}) under
some general conditions involving $f,g,$ even when at least one has a bad
sign at infinity, and to verify that any solution belonging to the
attractor, is bounded. Consequently, we obtain the uniqueness on the
trajectory attractor even for competing nonlinear mechanisms.

To better understand the larger scope of our results, we shall illustrate
their application to the reaction-diffusion equation (\ref{1.11bb}), (\ref%
{1.12bb}) for \emph{supercritical} nonlinearities, that is, for functions $f$
and $g$ satisfying the growth assumptions%
\begin{equation}
\lim_{\left\vert y\right\vert \rightarrow \infty }\frac{f^{^{\prime }}\left(
y\right) }{\left\vert y\right\vert ^{r_{1}-2}}=\left( r_{1}-1\right) c_{f}%
\text{, }\lim_{\left\vert y\right\vert \rightarrow \infty }\frac{g^{^{\prime
}}\left( y\right) }{\left\vert y\right\vert ^{r_{2}-2}}=\left(
r_{2}-1\right) c_{g},  \label{asss}
\end{equation}%
for some arbitrary $r_{2},r_{1}\geq 1,$ with $\max \left( r_{1},r_{2}\right)
\geq 2$, and some $c_{f},c_{g}\in \mathbb{R}\backslash \left\{ 0\right\} $.
Of course, our results below hold under more general assumptions on $f,g,$
see Section 3. In (\ref{asss}), we say that $f$ is dissipative if $c_{f}>0$
and non-dissipative if $c_{f}<0$ (the same applies to $g$). Let us assume
bounded $h_{1}\left( x\right) $ and $h_{2}\left( x\right) $. When both
nonlinear terms cooperate, i.e., both $f$ and $g$ are dissipative and%
\begin{equation}
f^{\prime }\left( y\right) \geq -\widetilde{c}_{f},\text{ }g^{\prime }\left(
y\right) \geq -\widetilde{c}_{g},\text{ for all }y\in \mathbb{R}
\label{limi}
\end{equation}%
(for some $\widetilde{c}_{f},\widetilde{c}_{g}>0$), then problem (\ref%
{1.11bb}), (\ref{1.12bb}) is well-posed and possesses a global attractor $%
\mathcal{A}_{gl}$ in the classical sense, bounded in $\mathcal{Z}%
:=W^{2,2}\left( \Omega \right) \cap L^{\infty }\left( \overline{\Omega }%
\right) $, regardless of the size of $r_{1}$ and $r_{2}$ (see Section 3.1;
cf. also \cite{GalNS}, when $g=0$). When the conditions (\ref{limi}) do
\emph{not} hold, we recall that uniqueness of weak solutions is not known in
general. However, if $f$ and $g$ are still dissipative, we can prove that
the reaction-diffusion system (\ref{1.11bb}), (\ref{1.12bb}) possesses a
(strong) trajectory attractor $\mathcal{A}_{tr}$, which is bounded in $%
L^{\infty }\left( \mathbb{R}_{+};\mathcal{Z}\right) .$ Moreover, uniqueness
holds on the attractor $\mathcal{A}_{tr}$ so that the long-term behavior of (%
\ref{1.11bb}), (\ref{1.12bb}) can be also characterized by a regular global
attractor%
\begin{equation}
\mathcal{A}_{gl}:=\mathcal{A}_{tr}\left( 0\right) ,  \label{gl}
\end{equation}%
which can be defined in the usual sense of dynamical systems (cf. Section
3.2).

For the case of competing nonlinearities, the following scenarios are
possible:

\noindent \textbf{Case (i)}: For the case of bulk dissipation (i.e., $%
c_{f}>0 $) and anti-dissipative behavior at the boundary $\Gamma $ (i.e., $%
c_{g}<0$), problem (\ref{1.11bb}), (\ref{1.12bb}) has at least one
globally-defined weak solution, which is bounded, if%
\begin{equation*}
\max \left( r_{2},2\left( r_{2}-1\right) \right) <r_{1}.
\end{equation*}%
Equality can be also allowed if the boundary condition is homogeneous, i.e.,
if $h_{2}=0$. Moreover, there exists a (strong)\ trajectory attractor $%
\mathcal{A}_{tr},$ bounded in $L^{\infty }\left( \mathbb{R}_{+};\mathcal{Z}%
\right) $, such that solutions are unique on the attractor. Thus, (\ref%
{1.11bb})-(\ref{1.12bb}) also possesses the (smooth) global attractor $%
\mathcal{A}_{gl},$ defined as in (\ref{gl}).

\noindent \textbf{Case (ii)}: On the other hand, in the case of boundary
dissipation ($c_{g}>0$) and internal non-dissipation ($c_{f}<0$), for every $%
L^{2}$-data we obtain that, if $r_{2}=r_{1}=2$ (which imply that $f$ and $g$
are sublinear)\ and $h_{2}=0,$ and%
\begin{equation}
(c_{f}+\left( \left\vert \Omega \right\vert \alpha \right) ^{-1}c_{g})\nu >2(%
\widetilde{C}_{\Omega }c_{g}r_{2})^{2},  \label{dom}
\end{equation}%
where%
\begin{equation*}
\alpha ^{-1}:=\int\limits_{\Gamma }\left( b\left( x\right) \right) ^{-1}dS,
\end{equation*}%
and $\widetilde{C}_{\Omega }$ is a proper Sobolev-Poincar\'{e} constant (see
Section 2, assumption (H4)), then (\ref{1.11bb})-(\ref{1.12bb}) is
well-posed in the classical sense and has a global attractor $\mathcal{A}%
_{gl}$, bounded in $\mathcal{Z}$ (in that case, $\left\vert f^{^{\prime
}}\left( y\right) \right\vert $ and $\left\vert g^{^{\prime }}\left(
y\right) \right\vert $ are bounded for all $y\in \mathbb{R}$ by (\ref{asss}%
), see Proposition \ref{uniq}). We note that the nonlinear balance condition
established in Section 3 can only be used to show the existence of a
(strong)\ trajectory attractor, whenever the interior non-dissipative term $%
f $ is sublinear as long as the boundary mechanism\ stays\ dissipative
(i.e., $c_{g}>0$), such that $g$\ suitably dominates $f$ (see, (\ref{dom})).
It would be interesting to see whether one can still construct attractors
for the case of a superlinear non-dissipative function $f,$ and some
dissipative boundary function $g$ of arbitrary growth. However, as we show
at the end of Section 3, we will see that, in this case, the superlinear
growth of the nondissipative function $f$ produces blowup\ in $L^{\infty }$%
-norm\ of some solutions even for arbitrary nonlinearities $g$ (see Section
3.4). Hence, the nonlinear boundary conditions (\ref{1.12bb})\ for equation (%
\ref{1.11bb}) \emph{cannot} prevent blowup of some solutions of (\ref{1.11bb}%
)-(\ref{1.12bb}) as long as the internal nonlinearity is strongly
non-dissipative (for instance, when $f$ satisfies (\ref{asss}) with $r_{1}>2$%
, such that $c_{f}<0$). Our main result in Section 3.4 (see, Theorem \ref%
{blowup}) extends some results in \cite{Be2} for the parabolic equation (\ref%
{1.11bb}) subject to linear dynamic boundary conditions (i.e., when $g\equiv
0$ in (\ref{1.12bb})), and the special cases treated in \cite{FQ, KZ}. In
this sense, the nonlinear balance conditions derived in Section 3, which
imply that the system (\ref{1.11bb})-(\ref{1.12bb}) is dissipative, are
optimal.

Finally, exploiting known parabolic regularity theory for PDE's of the form (%
\ref{1.11bb})-(\ref{1.12bb}), the regularity of the solution for (\ref%
{1.11bb})-(\ref{1.12bb}) increases as the functions $f$, $g$ and the domain $%
\Omega $ become more regular (see Remark \ref{reg}; cf. also \cite{GalNS}
and references therein). In particular, the global attractor $\mathcal{A}%
_{gl}$ consists of (smooth) classical solutions which are defined for all
times. Thus, using this additional regularity that solutions of (\ref{1.11bb}%
)-(\ref{1.12bb}) enjoy on $\mathcal{A}_{gl}$ (for all the above cases), we
obtain an explicit upper bound on the fractal dimension of $\mathcal{A}_{gl}$
for this reaction-diffusion system by imposing weaker assumptions on the
nonlinearities than in \cite{GalNS}. In particular, for any $N\geq 2$ and
for as long as (\ref{asss}) with $c_{f}>0$ and $c_{g}>0$ holds, we have%
\begin{equation}
\mathcal{C}_{0}\nu ^{-\left( N-1\right) }\leq \dim _{F}\mathcal{A}_{gl}\leq
\mathcal{C}_{1}\left( 1+\nu ^{-\left( N-1\right) }\right) ,  \label{dim}
\end{equation}%
for some positive constants $\mathcal{C}_{0},\mathcal{C}_{1}$ which can be
computed explicitly (see Section 3.3). The lower bound in (\ref{dim}) was
established in \cite{GalNS}, assuming dissipative $f$ and (homogeneous)\
linear boundary equations (i.e., $g=0$). We note that, for each fixed $\nu
>0 $, there is a discrepancy between the upper and lower bounds in (\ref{dim}%
) as $\mathcal{C}_{0}$ depends only on $\Omega $, $\Gamma $, $f,$ $g$,
whereas $\mathcal{C}_{1}$ is also a function of the $L^{\infty }$-norms of
the sources $h_{1},h_{2}.$ However, we observe that both the upper and lower
bounds are of the order $\nu ^{-\left( N-1\right) }$ as $\nu \rightarrow
0^{+}$, cf. (\ref{dim}). When $N=1$, the dimension of $\mathcal{A}_{gl}$ is
of the order $\nu ^{-1/2}$, as $\nu \rightarrow 0^{+}$. We recall that, for
the reaction-diffusion equation (\ref{1.11bb}) with the usual Dirichlet or
Neumann-Robin boundary condition (\ref{p2}), we have upper and lower bounds
of the order $\nu ^{-N/2},$ for any $N\geq 1$ (see, e.g., \cite{R, T}).
Thus, we have a much larger estimate (as a function of diffusion, as $\nu
\rightarrow 0^{+}$) for the global attractor $\mathcal{A}_{gl}$ in dimension
$N\geq 3$ (see also \cite{GalNS}).

We outline the plan of the paper, as follows. In Section \ref{well-posed},
we introduce some notations and preliminary facts, then we recall how to
prove the existence and stability of $L^{2}$-energy solutions to our system (%
\ref{1.11})-(\ref{1.12}). Section \ref{attractor} is devoted to the
existence of a bounded absorbing set and, then, of the global attractor $%
\mathcal{A}_{gl}$ for unique $L^{2}$-energy solutions. Then, we show that
weak solutions possess the $L^{2}$-$L^{\infty }$ smoothing property
exploiting some kind of iterative argument, and consequently, deduce the
existence of an absorbing set in $L^{\infty }$. In the final part of Section %
\ref{attractor}, we establish the existence of the trajectory attractor $%
\mathcal{A}_{tr}$\ for our problem and deduce some additional properties for
$\mathcal{A}_{tr}$, especially in the non-degenerate case when $a\left(
s\right) \equiv \nu >0$. A blowup result for (\ref{1.11bb})-(\ref{1.12bb})\
is also established in the case when $g\neq 0$. Finally, in the Appendix we
give some auxiliary results which are essential in the proofs.

\section{Well-posedness in $L^{2}$-space}

\label{well-posed}

We use the standard notation and facts from the dynamic theory of parabolic
equations (see, for instance, \cite{GW}). The natural space for our problem
is%
\begin{equation*}
\mathbb{X}^{s_{1},s_{2}}:=L^{s_{1}}(\Omega )\oplus L^{s_{2}}(\Gamma )=\{U=%
\binom{u_{1}}{u_{2}}:\;u_{1}\in L^{s_{1}}(\Omega ),\;u_{2}\in
L^{s_{2}}(\Gamma )\},
\end{equation*}%
$s_{1},s_{2}\in \left[ 1,+\infty \right] ,$ endowed with norm%
\begin{equation}
\left\Vert U\right\Vert _{\mathbb{X}^{s_{1},s_{2}}}=\left( \int_{\Omega
}\left\vert u_{1}\left( x\right) \right\vert ^{s_{1}}dx\right)
^{1/s_{1}}+\left( \int_{\Gamma }\left\vert u_{2}(x)\right\vert ^{s_{2}}\frac{%
dS_{x}}{b\left( x\right) }\right) ^{1/s_{2}},  \label{2.1}
\end{equation}%
if $s_{1},s_{2}\in \lbrack 1,\infty ),$ and
\begin{align*}
\Vert U\Vert _{\mathbb{X}^{\infty }}& :=\max \{\Vert u_{1}\Vert _{L^{\infty
}(\Omega )},\Vert u_{2}\Vert _{L^{\infty }(\Gamma )}\} \\
& \simeq \Vert u_{1}\Vert _{L^{\infty }(\Omega )}+\Vert u_{2}\Vert
_{L^{\infty }(\Gamma )}.
\end{align*}%
We agree to denote by $\mathbb{X}^{s}$ the space $\mathbb{X}^{s,s}.$
Identifying each function $u\in W^{1,p}(\Omega )$, with the vector $U:=%
\binom{u}{u_{\mid \Gamma }}$, it is easy to see that $W^{1,p}(\Omega )$ is a
dense subspace of $\mathbb{X}^{s}$ for $s\in \lbrack 1,\infty )$. Moreover,
we have%
\begin{equation*}
\mathbb{X}^{s}=L^{s}\left( \overline{\Omega },d\mu \right) ,\text{ }s\in %
\left[ 1,+\infty \right] ,
\end{equation*}%
where the measure $d\mu =dx_{\mid \Omega }\oplus dS_{x}/b\left( x\right)
_{\mid \Gamma }$ on $\overline{\Omega }$ is defined for any measurable set $%
A\subset \overline{\Omega }$ by $\mu (A)=|A\cap \Omega |+S(A\cap \Gamma )$.
Identifying each function $\theta \in C\left( \overline{\Omega }\right) $
with the vector $\Theta =\binom{\theta _{\mid \Omega }}{\theta _{\mid \Gamma
}}$, we have that $C(\overline{\Omega })$ is a dense subspace of $\mathbb{X}%
^{s}$ for every $s\in \lbrack 1,\infty )$ and a closed subspace of $\mathbb{X%
}^{\infty }$. In general, any vector $\theta \in \mathbb{X}^{s}$ will be of
the form $\binom{\theta _{1}}{\theta _{2}}$ with $\theta _{1}\in L^{s}\left(
\Omega ,dx\right) $ and $\theta _{2}\in L^{s}\left( \Gamma ,dS/b\left(
x\right) \right) ,$ and there need not be any connection between $\theta
_{1} $ and $\theta _{2}$. For domains $\Omega $ with Lipschitz boundary $%
\Gamma $, recall that we have $W^{k,p}\left( \Omega \right) \subset
L^{p_{s}}\left( \Omega \right) $, with $p_{s}=\frac{pN}{N-pk}$ if $pk<N,$
and $1\leq p_{s}<\infty ,$ if $N=pk$. Moreover the trace operator $%
Tr_{D}(u):=u_{|_{\Gamma }},$ initially defined for $u\in C^{k}(\overline{%
\Omega }),$ has an extension to a bounded linear operator from $%
W^{k,p}(\Omega )$ into $L^{q_{s}}(\partial \Omega ),$ where $q_{s}:=\frac{%
p(N-1)}{N-pk}$ if $pk<N$, and $1\leq q_{s}<\infty $ if $N=pk$. For $p>Nk$,
we have $W^{k,p}\left( \Omega \right) \subset C^{\zeta ,\widetilde{l}}\left(
\overline{\Omega }\right) ,$ for some $\zeta ,\widetilde{l}$. We also
recall, on account of well-known generalized Poincar\'{e}-type inequalities
(see, e.g., \cite{Le}), that%
\begin{equation}
\left\Vert u\right\Vert _{W^{1,p}\left( \Omega \right) }:=\left\Vert \nabla
u\right\Vert _{\left( L^{p}\left( \Omega \right) \right) ^{N}}+l\left(
u\right)  \label{equiv}
\end{equation}%
is a norm on $W^{1,p}\left( \Omega \right) $, which is equivalent to the
usual one, for any of the following expressions for $l$:%
\begin{equation*}
l\left( u\right) :=\left( \int_{\Gamma }\left\vert u\right\vert ^{s}\frac{%
dS_{x}}{b}\right) ^{1/s},\text{ }l\left( u\right) :=\left( \int_{\Omega
}\left\vert u\right\vert ^{s}dx\right) ^{1/s},
\end{equation*}%
for any $1\leq s\leq p$. Next, for each $p>1$, we let%
\begin{equation*}
\mathbb{V}^{k,p}=\{U:=\binom{u}{u_{\mid \Gamma }}:\;u\in W^{k,p}(\Omega )\}
\end{equation*}%
and endow it with the norm $\left\Vert \cdot \right\Vert _{\mathbb{V}^{k,p}}$
given by%
\begin{equation*}
\left\Vert U\right\Vert _{\mathbb{V}^{k,p}}=\left\Vert u\right\Vert
_{W^{k,p}\left( \Omega \right) }+\left\Vert u_{\mid \Gamma }\right\Vert
_{W^{k-1/p,p}\left( \Gamma \right) }.
\end{equation*}%
It easy to see that we can identify $\mathbb{V}^{k,p}$ with $W^{k,p}\left(
\Omega \right) \oplus W^{k-1/p,p}\left( \Gamma \right) $ under this norm.
Moreover, we emphasize that $\mathbb{V}^{k,p}$ is not a product space and
since $W^{k,p}(\Omega )\hookrightarrow W^{k-1/p,p}\left( \Gamma \right) $ by
trace theory, $\mathbb{V}^{k,p}$ is topologically isomorphic to $%
W^{k,p}\left( \Omega \right) $ in the obvious way. It is also immediate that
$\mathbb{V}^{k,p}$ is compactly embedded into $\mathbb{X}^{2},$ for any $%
p>p_{0}:=2N/(N+2)$ and $k\geq 1$. From now on, we denote by $\left\Vert
\cdot \right\Vert _{W^{k,p}\left( \Omega \right) }$ and $\left\Vert \cdot
\right\Vert _{W^{k,q}\left( \Gamma \right) }$ the norms on $W^{k,p}\left(
\Omega \right) $ and $W^{k,q}\left( \Gamma \right) ,$ respectively. Also, $%
\langle \cdot ,\cdot \rangle _{s}$ and $\langle \cdot ,\cdot \rangle
_{s,\Gamma }$ stand for the usual scalar product in $L^{s}\left( \Omega
\right) $ and $L^{s}\left( \Gamma \right) $, respectively. We also agree to
denote by $\mathbb{V}^{p}$ the space $\mathbb{V}^{1,p},$ and $\langle \cdot
,\cdot \rangle $ the duality between $X$ and $X^{\ast },$ for some generic
Banach space $X$.

Our first goal in this paper is to give a nonlinear balance between $f$ and $%
g$ which implies dissipativity of (\ref{1.11})-(\ref{1.12}), even when both
the nonlinear terms contribute in opposite directions. More precisely, we
wish to prove the existence of (globally well-defined) weak solutions,
provided that the nonlinearities satisfy (possibly part of) the assumptions
listed below:

(H1) Let $b\left( y\right) :=a(\left\vert y\right\vert ^{2})y,$ $y\in
\mathbb{R}^{N}$ and assume that $a\in {\normalsize C\left( \mathbb{R},%
\mathbb{R}\right) }$, $b\in C^{1}${\normalsize $\left( \mathbb{R}^{N},%
\mathbb{R}^{N}\right) $ }satisfy the following conditions:{\normalsize
\begin{equation}
\left\{
\begin{array}{ll}
\left\vert a(\left\vert y\right\vert ^{2})\right\vert \leq
c_{1}(1+\left\vert y\right\vert ^{p-2}), & \forall y\in \mathbb{R}^{N}, \\
\left\langle b\left( y\right) ,y\right\rangle _{\mathbb{R}^{N}}=a(\left\vert
y\right\vert ^{2})\left\vert y\right\vert ^{2}\geq \nu \left\vert
y\right\vert ^{p}, & \forall y\in \mathbb{R}^{N},%
\end{array}%
\right.  \label{1.1.2}
\end{equation}%
}for some constants{\normalsize \ $c_{1},\nu >0$. }Moreover, assume that $b$
is monotone nondecreasing, i.e.,%
\begin{equation}
\left\langle a(\left\vert y_{1}\right\vert ^{2})y_{1}-a(\left\vert
y_{2}\right\vert ^{2})y_{2},y_{1}-y_{2}\right\rangle _{\mathbb{R}^{N}}\geq 0,%
\text{ for all }y_{1},y_{2}\in \mathbb{R}^{N}.  \label{1.1.7}
\end{equation}

\medskip (H2) $f,g\in C^{1}\left( \mathbb{R},\mathbb{R}\right) $ satisfy%
\begin{equation}
\lim_{\left\vert y\right\vert \rightarrow +\infty }\inf f^{\prime }\left(
y\right) >0,\quad \lim_{\left\vert y\right\vert \rightarrow +\infty }\inf
g^{\prime }\left( y\right) >0.  \label{2.33}
\end{equation}

\medskip (H3) (a) $f,g\in C\left( \mathbb{R},\mathbb{R}\right) $ satisfy the
growth assumptions%
\begin{equation}
\left\vert f\left( y\right) \right\vert \leq c_{f}\left( 1+\left\vert
y\right\vert ^{r_{1}-1}\right) ,\quad \left\vert g\left( y\right)
\right\vert \leq c_{g}\left( 1+\left\vert y\right\vert ^{r_{2}-1}\right)
,\quad \forall y\in {\mathbb{R}}\text{,}  \label{2.5}
\end{equation}%
for some positive constants $c_{f},$ $c_{g}$, and some $r_{1},r_{2}\geq 1$.

\ \ \ \ \ \ \ (b) $f,g\in C\left( \mathbb{R},\mathbb{R}\right) $ satisfy%
\begin{equation*}
\begin{array}{ll}
\left\vert f\left( y\right) \right\vert \leq c_{f}\left( 1+\left\vert
y\right\vert ^{r_{1}-1}\right) , &  \\
c_{g}\left\vert y\right\vert ^{r_{2}}-c\leq g\left( y\right) y\leq
\widetilde{c}_{g}\left\vert y\right\vert ^{r_{2}}+c, & \forall y\in {\mathbb{%
R}}\text{,}%
\end{array}%
\end{equation*}%
for some appropriate positive constants and some $r_{1},r_{2}\geq 1$.

(H4) Let $\lambda =\left( \int_{\Gamma }b^{-1}dS\right) ^{-1}$, and suppose
that $g\in C^{1}\left( \mathbb{R},\mathbb{R}\right) $. There exists $%
\varepsilon \in \left( 0,\nu /q\right) ,$ with $\nu $ as in (\ref{1.1.2}),
such that,%
\begin{equation}
\lim_{\left\vert y\right\vert \rightarrow +\infty }\inf \frac{f\left(
y\right) y+\left( \left\vert \Omega \right\vert \lambda \right) ^{-1}g\left(
y\right) y-\frac{\widetilde{C}_{\Omega }^{q}}{\left( \varepsilon p\right)
^{q/p}q}\left\vert g^{^{\prime }}\left( y\right) y+g\left( y\right)
\right\vert ^{q}}{\left\vert y\right\vert ^{r_{1}}}>0,  \label{NB}
\end{equation}%
for some $r_{1}\geq p.$ Here%
\begin{equation}
\widetilde{C}_{\Omega }=\left\{
\begin{array}{ll}
C_{\Omega ,b}\left( \lambda \left\vert \Omega \right\vert \right) ^{-1}, &
\text{if (H3a) holds} \\
\frac{C_{\Omega ,b}}{2}\left( \lambda \left\vert \Omega \right\vert \right)
^{-1}, & \text{if (H3b) holds,}%
\end{array}%
\right.  \label{NBbis}
\end{equation}%
and $C_{\Omega ,b}$ is the best Sobolev constant in the following Poincar%
\'{e}'s inequality:%
\begin{equation}
\left\Vert \phi -\lambda \langle \phi /b,1\rangle _{1,\Gamma }\right\Vert
_{L^{s}\left( \Omega \right) }\leq C_{\Omega ,b}\left\Vert \nabla \phi
\right\Vert _{L^{s}\left( \Omega \right) },\text{ }s\geq 1  \label{PI}
\end{equation}%
(see, e.g., \cite[Lemma 3.1]{RB}).

We observe that condition (H4) provides an exact balance between the two
nonlinear mechanisms. As we shall see, this balance will depend both upon
the sign and growth rate of $f$ and $g$ at infinity (cf. also, \cite{RB}).

We have the following rigorous notion of weak solution to (\ref{1.11})-(\ref%
{1.12}), with initial condition $u\left( 0\right) =u_{0},$ similar to \cite%
{GW}.

\begin{definition}
\label{weak}Let $p\in (\frac{2N}{N+2},+\infty )\cap \left( 1,+\infty \right)
$, and let $h_{1}\left( x\right) \in {\normalsize L}^{r_{1}^{^{\prime
}}}\left( \Omega \right) $, $h_{2}\left( x\right) \in {\normalsize L}%
^{r_{2}^{^{\prime }}}\left( \Gamma \right) ,$ where $r_{i}^{^{\prime }}$ is
the dual conjugate of $r_{i}$. The pair $U\left( t\right) =\binom{u\left(
t\right) }{v\left( t\right) }$ is said to be a weak solution if $v\left(
t\right) =u\left( t\right) _{\mid \Gamma }$, in the trace sense, for a.e. $%
t\in \left( 0,T\right) ,$ for any $T>0$, and $U$ fulfills
\end{definition}

\begin{equation*}
\left\{
\begin{array}{l}
U\left( t\right) \in L^{\infty }\left( \left[ 0,+\infty \right) ;\mathbb{X}%
^{2}\right) \cap W_{loc}^{1,s}\left( \left[ 0,+\infty \right) ;\left(
\mathbb{V}^{k,p}\right) ^{\ast }\right) \\
u\left( t\right) \in L_{\text{loc}}^{p}\left( \left[ 0,+\infty \right)
;W^{1,p}\left( \Omega \right) \right) , \\
v\left( t\right) \in L_{\text{loc}}^{p}\left( \left[ 0,+\infty \right)
;W^{1-1/p,p}\left( \Gamma \right) \right) ,%
\end{array}%
\right.
\end{equation*}%
for $s=\min (q,r_{1}^{^{\prime }},r_{2}^{^{\prime }}),$ $q:=p/\left(
p-1\right) $, and%
\begin{equation}
k=\max (1,\frac{N}{p}-\frac{N}{r_{1}},\frac{N}{p}-\frac{N-1}{r_{2}}).
\label{kcond}
\end{equation}%
Moreover, the following identity%
\begin{align}
& \left\langle \partial _{t}U,\Xi \right\rangle _{\mathbb{X}%
^{2}}+\left\langle a\left( \left\vert \nabla u\right\vert ^{2}\right) \nabla
u,\nabla \sigma \right\rangle _{2}+\left\langle f\left( u\right) ,\sigma
\right\rangle _{2}+\left\langle g\left( v\right) /b,\sigma _{\mid \Gamma
}\right\rangle _{2,\Gamma }  \label{VF} \\
& =\left\langle h_{1},\sigma \right\rangle _{2}+\left\langle h_{2}/b,\sigma
_{\mid \Gamma }\right\rangle _{2,\Gamma },  \notag
\end{align}%
holds for all $\Xi =\binom{\sigma }{\sigma _{\mid \Gamma }}\in \mathbb{V}%
^{k,p},$ a.e. $t\in \left( 0,T\right) $. Finally, we have, in the space $%
\mathbb{X}^{2}$,%
\begin{equation}
U\left( 0\right) =\binom{u_{0}}{v_{0}}=:U_{0},  \label{ini}
\end{equation}%
where $u\left( 0\right) =u_{0}$ almost everywhere in $\Omega $, and $v\left(
0\right) =v_{0}$ almost everywhere in $\Gamma $. Note that in this setting, $%
v_{0}$ need not be the trace of $u_{0}$ at the boundary.

We can cast the weak formulation (\ref{VF}) into a proper functional
equation by defining suitable operators. To this end, let $\langle \cdot
,\cdot \rangle $ denote the duality between $\mathbb{V}^{p}$ and $\left(
\mathbb{V}^{p}{}\right) ^{\ast }$. Define the form%
\begin{equation*}
\widetilde{\mathcal{B}}_{p}(U_{1},U_{2}):=\int_{\Omega }a\left( \left\vert
\nabla u_{1}\right\vert ^{2}\right) \nabla u_{1}\cdot \nabla
u_{2}dx+\int_{\Omega }|u_{1}|^{p-2}u_{1}u_{2}dx,
\end{equation*}%
for all $U_{i}=\binom{u_{i}}{u_{i\mid \Gamma }}\in \mathbb{V}^{p},$ $i=1,2$.
Note that%
\begin{align}
\widetilde{\mathcal{B}}_{p}(U_{1},U_{2})& =-\int_{\Omega }div\left( a\left(
\left\vert \nabla u_{1}\right\vert ^{2}\right) \nabla u_{1}\right) u_{2}dx
\label{calc} \\
& +\int_{\Gamma }b\left( x\right) a\left( \left\vert \nabla u_{1}\right\vert
^{2}\right) \partial _{\mathbf{n}}u_{1}u_{2}\frac{dS}{b\left( x\right) }%
+\int_{\Omega }|u_{1}|^{p-2}u_{1}u_{2}dx.  \notag
\end{align}%
It follows from Lemma \ref{lem-hemi-2} (see Appendix), that for each $U=%
\binom{u}{v}\in \mathbb{V}^{p}$, there exists $\mathcal{B}_{p}(U)\in \left(
\mathbb{V}^{p}{}\right) ^{\ast }$ such that
\begin{equation}
\widetilde{\mathcal{B}}_{p}(U,W)=\langle \mathcal{B}_{p}(U),W\rangle ,
\label{deff}
\end{equation}%
for every $W\in \mathbb{V}^{p}$. Hence, this relation defines an operator $%
\mathcal{B}_{p}:\mathbb{V}^{p}\rightarrow \left( \mathbb{V}^{p}\right)
^{\ast },$ which is bounded. Exploiting Lemma \ref{lem-hemi-2} once again,
it is easy to see that $\mathcal{B}_{p}$ is monotone and coercive. It also
follows that $\mathcal{B}_{p}\left( \mathbb{V}^{p}\right) =\left( \mathbb{V}%
^{p}\right) ^{\ast }$ (see, e.g., \cite{Br}). Thus, we end up with the
following functional form%
\begin{equation}
\partial _{t}U+\mathcal{B}_{p}U+\mathcal{F}\left( U\right) =\mathcal{G}%
\left( x\right) ,  \label{eq}
\end{equation}%
where $\mathcal{G}\left( x\right) \mathcal{=}\binom{h_{1}\left( x\right) }{%
h_{2}\left( x\right) },$ and the operators $\mathcal{B}_{p}:D\left( \mathcal{%
B}_{p}\right) \rightarrow \mathbb{X}^{2},$ $\mathcal{F}:D(\mathcal{F}%
)\subset \mathbb{X}^{2}\rightarrow \mathbb{X}^{2}$ are given, formally, by%
\begin{equation}
\mathcal{B}_{p}U=\left(
\begin{array}{c}
-\text{div}\left( a\left( \left\vert \nabla u\right\vert ^{2}\right) \nabla
u\right) +\left\vert u\right\vert ^{p-2}u \\
b\left( x\right) a\left( \left\vert \nabla u\right\vert ^{2}\right) \partial
_{{\mathbf{n}}}u%
\end{array}%
\right) ,  \label{2.18}
\end{equation}%
\begin{equation*}
\mathcal{F}\left( U\right) =\left(
\begin{array}{c}
f\left( u\right) -\left\vert u\right\vert ^{p-2}u \\
g\left( v\right)%
\end{array}%
\right) .
\end{equation*}

We aim to prove some regularity results for the weak solutions constructed
in Definition \ref{weak}. In \cite{GW}, solutions were constructed with aid
from a Galerkin approximation scheme by imposing additional growth
restrictions on the nonlinearities $f,g$ so that they are essentially
dominated by a monotone operator associated with the $p$-Laplacian. However,
the application of this scheme seems to be problematic in our context since
the solutions constructed with Definition \ref{weak} are much weaker than
those constructed in \cite[Definition 2.3]{GW}. Therefore, we need to rely
on another scheme which is based on the construction of classical (smooth)\
solutions to the non-degenerate analogue of (\ref{1.11})-(\ref{1.12}). One
of the advantages of this construction is that now every weak solution can
be approximated by regular ones and the justification of our estimates for
such solutions is immediate. Thus, for each $\epsilon >0$, let us consider
the following non-degenerate parabolic problem:%
\begin{align}
\partial _{t}u-\text{div}\left( a_{\epsilon }\left( \left\vert \nabla
u\right\vert ^{2}\right) \nabla u\right) +f\left( u\right) & =h_{1}(x),\text{
in }\Omega \times \left( 0,\infty \right) ,  \label{ap} \\
\partial _{t}u+b\left( x\right) a_{\epsilon }\left( \left\vert \nabla
u\right\vert ^{2}\right) \partial _{\mathbf{n}}u+g\left( u\right) &
=h_{2}(x),\text{ on }\Gamma \times \left( 0,\infty \right) ,  \notag
\end{align}%
where $a_{\epsilon }\left( s\right) =a\left( s+\epsilon \right) >0,$ for all
$s\in \mathbb{R}$, subject to the initial conditions%
\begin{equation}
u_{\epsilon }\left( 0\right) =u_{0\epsilon },\text{ }v_{\epsilon }\left(
0\right) =u_{0\epsilon \mid \Gamma }.  \label{ini2}
\end{equation}%
Let $u_{0\epsilon }\in C^{\infty }\left( \overline{\Omega }\right) $ such
that%
\begin{equation*}
U_{\epsilon }\left( 0\right) \rightarrow U\left( 0\right) =U_{0}\text{ in }%
\mathbb{X}^{2}.
\end{equation*}%
Then, the approximate problem (\ref{ap})-(\ref{ini2}) admits a unique
(smooth) classical solution with%
\begin{equation}
u_{\epsilon }\in C^{1}\left( \left[ 0,t_{\ast }\right] ;C^{\infty }\left(
\overline{\Omega }\right) \right)  \label{smooth}
\end{equation}%
for some $t_{\ast }>0$ and each $\epsilon >0$ (see \cite{Escher, Escher2, CJ}%
). Being pedants, we cannot apply the main results of \cite{Escher} (cf.
also \cite{Escher2}) directly to equations (\ref{ap})-(\ref{ini2}) since the
functions $a_{\epsilon },$ $f,$ $g$ and the external forces $h_{1},$ $h_{2}$
are not smooth enough. Moreover, the solutions constructed this way may only
exists locally in time for some interval $\left[ 0,t_{\ast }\right) $.
However, taking sequences $h_{1\epsilon }\in C^{\infty }\left( \overline{%
\Omega }\right) ,$ $h_{2\epsilon }\in C^{\infty }\left( \Gamma \right) $
such that $h_{1\epsilon }\rightarrow h_{1}$ in $L^{\infty }\left( \Omega
\right) ,$ and $h_{2\epsilon }\rightarrow h_{2}$ in $L^{\infty }\left(
\Gamma \right) ,$ respectively, and by approximating the functions $%
a_{\epsilon },$ $f$, $g$ by smooth ones, say, in $C^{\infty }\left( \mathbb{R%
},\mathbb{R}\right) $, we may apply Remark \ref{rem} below for the solutions
of the approximate equations, and deduce the existence of a globally defined
in $\mathbb{X}^{\infty }$-norm solution to (\ref{ap})-(\ref{ini2}). Indeed,
taking advantage of the fact that $u_{0\epsilon }\in C^{\infty }\left(
\overline{\Omega }\right) ,$ the global $\mathbb{X}^{\infty }$-a priori
bound for $u_{\epsilon }$ guarantees its global existence in at least $%
\mathbb{V}^{p}\cap \mathbb{X}^{\infty }$-norm, which turns out to be
sufficient for our purpose. As we shall see in the next section, this bound
can be naturally obtained under the above assumptions on the nonlinearities
by performing a modified Alikakos-Moser iteration argument (see Theorem \ref%
{linff} and Remark \ref{rem} below).

We shall now deduce the first result concerning the solvability of problem (%
\ref{1.11})-(\ref{1.12}).

\begin{theorem}
\label{T1}Let $a,$ $f$ and $g$ satisfy either the assumptions (H1), (H3a),
(H4) with%
\begin{equation}
\max \left( r_{2},q\left( r_{2}-1\right) \right) <r_{1},  \label{maxi}
\end{equation}%
or (H1), (H3b) and (H4). Then, for any initial data $U_{0}\in \mathbb{X}%
^{2}, $ there exists at least one (globally defined) weak solution $U\left(
t\right) $\ in the sense of Definition \ref{weak}.
\end{theorem}

\begin{proof}
We divide the proof into several steps according to the different hypotheses
being used.

\textbf{Step 1}. (i) We shall now derive some basic apriori estimates for $%
U=U_{\epsilon }\left( t\right) $ when $p\leq N$, assuming that (H1), (H3a)
and (H4) are satisfied. The case $p>N$ can be treated analogously. The
following estimates will be deduced by a formal argument, which can be
justified by means of the approximation procedure devised above. Also, for
practical purposes $c$ will denote a positive constant that is independent
of time, $\epsilon >0$ and initial data, but which only depends on the other
structural parameters. Such a constant may vary even from line to line. Note
that the smooth solutions constructed in (\ref{ap})-(\ref{ini2}) also
satisfy the weak formulation (\ref{VF}). Thus, the key choices $\sigma
=u_{\epsilon }\left( t\right) $, $\sigma _{\mid \Gamma }=v_{\epsilon }\left(
t\right) /b$ in (\ref{VF}) are justified. After standard transformations, in
view of assumption (H1), we obtain
\begin{align}
& \frac{1}{2}\frac{d}{dt}\left\Vert U\left( t\right) \right\Vert _{\mathbb{X}%
^{2}}^{2}+\nu \left\Vert \nabla u\left( t\right) \right\Vert _{L^{p}\left(
\Omega \right) }^{p}+\left\langle f\left( u\left( t\right) \right) ,u\left(
t\right) \right\rangle _{2}+\left\langle g\left( v\left( t\right) \right)
,v\left( t\right) /b\right\rangle _{2,\Gamma }  \label{2.7} \\
& \leq \left\langle h_{1},u\left( t\right) \right\rangle _{2}+\left\langle
h_{2}/b,v\left( t\right) \right\rangle _{2,\Gamma }.  \notag
\end{align}%
Following \cite{RB}, we can now write%
\begin{align}
& \left\langle f\left( u\right) ,u\right\rangle _{2}+\left\langle g\left(
v\right) ,v/b\right\rangle _{2,\Gamma }  \label{2.7bis} \\
& =\left\langle f\left( u\right) u+c_{B}g\left( u\right) u,1\right\rangle
_{2}-c_{B}\left\langle g\left( u\right) u-\lambda \left\langle g\left(
v\right) v/b,1\right\rangle _{2,\Gamma }\right\rangle _{2},  \notag
\end{align}%
where%
\begin{equation*}
c_{B}:=\frac{\lambda ^{-1}}{\left\vert \Omega \right\vert },\text{ }\lambda
=\left( \int_{\Gamma }b^{-1}dS\right) ^{-1}.
\end{equation*}%
Applying inequality (\ref{PI}) to the last term on the right-hand side of (%
\ref{2.7bis}) yields%
\begin{align}
& c_{B}\left\vert \left\langle g\left( u\right) u-\lambda \left\langle
g\left( v\right) v/b,1\right\rangle _{2,\Gamma }\right\rangle _{2}\right\vert
\label{2.7tris} \\
& \leq \widetilde{C}_{\Omega }\left\Vert \nabla \left( g\left( u\right)
u\right) \right\Vert _{L^{1}\left( \Omega \right) }=\widetilde{C}_{\Omega
}\left\Vert \left( g^{^{\prime }}\left( u\right) u+g\left( u\right) \right)
\nabla u\right\Vert _{L^{1}\left( \Omega \right) }  \notag \\
& \leq \varepsilon \left\Vert \nabla u\right\Vert _{L^{p}\left( \Omega
\right) }^{p}+\frac{\widetilde{C}_{\Omega }^{q}}{\left( \varepsilon p\right)
^{q/p}q}\left\Vert g^{^{\prime }}\left( u\right) u+g\left( u\right)
\right\Vert _{L^{q}\left( \Omega \right) }^{q}  \notag
\end{align}%
with $\widetilde{C}_{\Omega }=C_{\Omega ,b}\left( \lambda \left\vert \Omega
\right\vert \right) ^{-1},$ and we recall that $q$ is conjugate to $p$.
Since assumption (H4) holds for some $\varepsilon \in \left( 0,\nu /q\right)
,$ then from (\ref{2.7bis})-(\ref{2.7tris}), we obtain
\begin{align}
& \frac{1}{2}\frac{d}{dt}\left\Vert U\left( t\right) \right\Vert _{\mathbb{X}%
^{2}}^{2}+\left( \nu -\varepsilon \right) \left\Vert \nabla u\left( t\right)
\right\Vert _{L^{p}\left( \Omega \right) }^{p}+c\left\Vert u\left( t\right)
\right\Vert _{L^{r_{1}}\left( \Omega \right) }^{r_{1}}  \label{2.7q} \\
& \leq \left\langle h_{1},u\left( t\right) \right\rangle _{2}+\left\langle
h_{2}/b,v\left( t\right) \right\rangle _{2,\Gamma }+c,  \notag
\end{align}%
for some positive constant $c>0$, independent of $U,$ $t$ and $\epsilon $.
Exploiting the estimate in Lemma \ref{boundary} (see Appendix), and then
using H\"{o}lder and Young inequalities, we can bound the term on the
right-hand side of (\ref{2.7q}) by%
\begin{align*}
& \left( c\left\Vert h_{1}\right\Vert _{L^{r_{1}^{^{\prime
}}}}^{r_{1}^{^{\prime }}}+\frac{c}{2}\left\Vert u\right\Vert
_{L^{r_{1}}\left( \Omega \right) }^{r_{1}}\right) \\
& +\left[ \delta \left\Vert \nabla u\right\Vert _{L^{p}\left( \Omega \right)
}^{p}+C_{\delta }\left( \left\Vert u\right\Vert _{L^{\gamma }\left( \Omega
\right) }^{\gamma }+\left\Vert h_{2}/b\right\Vert _{L^{r_{2}^{^{\prime
}}}\left( \Gamma \right) }^{r_{2}^{^{\prime }}}+1\right) \right] ,
\end{align*}%
with $\gamma =\max \left( r_{2},p\left( r_{2}-1\right) /\left( p-1\right)
\right) $, for a sufficiently small $\delta >0$, and sufficiently large $%
C_{\delta }>0$. Since by assumption $\gamma <r_{1}$, we can control the $%
L^{\gamma }$-norm of $u$ in terms of the $L^{r_{1}}$-norm of the solution
(i.e., $\left\Vert u\right\Vert _{L^{\gamma }}^{\gamma }\leq \kappa
\left\Vert u\right\Vert _{L^{r_{1}}}^{r_{1}}+C_{\kappa },$ $\kappa \ll 1$).
Thus, we get for suitable choices of $\varepsilon \in \left( 0,\nu /q\right)
$ and $\delta \in \left( 0,\nu /q\right) ,$ the following inequality
\begin{align}
& \frac{1}{2}\frac{d}{dt}\left\Vert U\left( t\right) \right\Vert _{\mathbb{X}%
^{2}}^{2}+\frac{c}{p}\left\Vert \nabla u\left( t\right) \right\Vert
_{L^{p}\left( \Omega \right) }^{p}+\frac{c}{2}\left\Vert u\left( t\right)
\right\Vert _{L^{r_{1}}\left( \Omega \right) }^{r_{1}}  \label{dissi} \\
& \leq c\left( 1+\left\Vert h_{1}\right\Vert _{L^{r_{1}^{^{\prime
}}}}^{r_{1}^{^{\prime }}}+\left\Vert h_{2}/b\right\Vert _{L^{r_{2}^{^{\prime
}}}\left( \Gamma \right) }^{r_{2}^{^{\prime }}}\right) ,  \notag
\end{align}%
for almost all $t\geq 0$. Recalling that $L^{r_{1}}\left( \Omega \right)
\subset L^{p}\left( \Omega \right) $ (we have, $r_{1}\geq p$), we can now
integrate this inequality over $\left[ 0,T\right] $ to deduce%
\begin{align}
U_{\epsilon }& \in L^{\infty }\left( \left[ 0,T\right] ;\mathbb{X}%
^{2}\right) \cap L^{p}\left( \left[ 0,T\right] ;\mathbb{V}^{p}\right) ,
\label{est1} \\
u_{\epsilon }& \in L^{r_{1}}\left( \left[ 0,T\right] \times \Omega \right) ,
\notag
\end{align}%
uniformly with respect to $\epsilon >0$. On account of these bounds, we get%
\begin{equation*}
\mathcal{B}_{p,\epsilon }\left( U_{\epsilon }\right) \in L^{q}\left( \left[
0,T\right] ;\left( \mathbb{V}^{p}\right) ^{\ast }\right) \subseteq
L^{q}\left( \left[ 0,T\right] ;\left( \mathbb{V}^{k,p}\right) ^{\ast
}\right) ,
\end{equation*}%
uniformly in $\epsilon >0$, for any $k\geq 1$ (cf. Lemma \ref{lem-hemi-2};
see also \cite{GW}). Here $\mathcal{B}_{p,\epsilon }$ is the monotone
operator associated with the function $a_{\epsilon }$ (see (\ref{deff}), (%
\ref{2.18})). Moreover, by Lemma \ref{boundary} in the Appendix, we get at
once%
\begin{equation}
v_{\epsilon }=u_{\epsilon \mid \Gamma }\in L^{r_{2}}\left( \left[ 0,T\right]
\times \Gamma \right) \cap L^{p}\left( \left[ 0,T\right] \times \Gamma
\right) ,  \label{est2}
\end{equation}%
uniformly in $\epsilon $. Due to assumption (H3a), from (\ref{est1})-(\ref%
{est2}), we deduce%
\begin{equation}
\binom{f\left( u_{\epsilon }\right) }{g\left( v_{\epsilon }\right) }\in
L^{r_{1}^{^{\prime }}}(\left[ 0,T\right] \times \Omega )\times
L^{r_{2}^{^{\prime }}}(\left[ 0,T\right] \times \Gamma ).  \label{est3}
\end{equation}%
Thus, $\mathcal{F}\left( U_{\epsilon }\right) $ is uniformly (in $\epsilon $%
) bounded in $L^{s}(\left[ 0,T\right] ;\mathbb{X}^{s}),$ which implies%
\begin{equation*}
\mathcal{G}\left( x\right) -\mathcal{B}_{p,\epsilon }U_{\epsilon }-\mathcal{F%
}\left( U_{\epsilon }\right) \in L^{s}\left( \left[ 0,T\right] ;\left(
\mathbb{V}^{k,p}\right) ^{\ast }\right) ,
\end{equation*}%
with $s=\min (q,r_{1}^{^{\prime }},r_{2}^{^{\prime }})>1$. Therefore, $%
\partial _{t}U_{\epsilon }$ is bounded in $L^{s}\left( \left[ 0,T\right]
;\left( \mathbb{V}^{k,p}\right) ^{\ast }\right) $, uniformly with respect to
$\epsilon >0$, for some $k\geq 1$. Indeed, having chosen $k$ so that $%
W^{k,p}\subset L^{r_{1}}\left( \Omega \right) \subset L^{p}\left( \Omega
\right) $ and $W^{k-1/p,p}\subset L^{r_{2}}\left( \Gamma \right) $ (in
particular, it holds $\mathbb{V}^{k,p}\subset \mathbb{X}^{r_{1},r_{2}}$ with
continuous inclusion), so if $V\in \left( \mathbb{X}^{r_{1},r_{2}}\right)
^{\ast },$ then $V\in \left( \mathbb{V}^{k,p}\right) ^{\ast }$. Thus,
equation (\ref{eq}) holds as an equality in $L^{s}\left( \left[ 0,T\right]
;\left( \mathbb{V}^{k,p}\right) ^{\ast }\right) $ and it can be considered
distributionally in the space $\mathcal{D}^{^{\prime }}\left( \left[ 0,T%
\right] ;\left( \mathbb{V}^{k,p}\right) ^{\ast }\right) $. The existence of
a weak solution is then based on monotone operator arguments, followed by a
passage to limit as $\epsilon \rightarrow 0,$ and can be carried over
exactly as in the proof of \cite[Theorem 2.6]{GW}. We shall briefly describe
the details below in \textbf{Step 2}.

\textbf{Step 1}. (ii) We will now deduce the apriori bounds (\ref{est1})-(%
\ref{est2}), if one assumes (H3b) instead of (H3a). According to (\ref{2.7}%
), in light of inequality (H3b) for $g,$ we have the following:
\begin{align}
& \frac{1}{2}\frac{d}{dt}\left\Vert U\left( t\right) \right\Vert _{\mathbb{X}%
^{2}}^{2}+\nu \left\Vert \nabla u\left( t\right) \right\Vert _{L^{p}\left(
\Omega \right) }^{p}+\left\langle f\left( u\left( t\right) \right) ,u\left(
t\right) \right\rangle _{2}  \label{2.7s} \\
& +\frac{1}{2}\left\langle g\left( v\left( t\right) \right) ,v\left(
t\right) /b\right\rangle _{2,\Gamma }+\frac{c_{g}}{2}\left\Vert v\left(
t\right) \right\Vert _{L^{r_{2}}\left( \Gamma ,b^{-1}dS\right) }^{r_{2}}
\notag \\
& \leq \left\langle h_{1},u\left( t\right) \right\rangle _{2}+\left\langle
h_{2}/b,v\left( t\right) \right\rangle _{2,\Gamma }+c,  \notag
\end{align}%
for some positive constant $c.$ We can write, as in (\ref{2.7bis}),%
\begin{align*}
& \left\langle f\left( u\right) ,u\right\rangle _{2}+\frac{1}{2}\left\langle
g\left( v\right) ,v/b\right\rangle _{2,\Gamma } \\
& =\left\langle f\left( u\right) u+\frac{c_{B}}{2}g\left( u\right)
u,1\right\rangle _{2}-\frac{c_{B}}{2}\left\langle g\left( u\right) u-\lambda
\left\langle g\left( v\right) v/b,1\right\rangle _{2,\Gamma }\right\rangle
_{2},
\end{align*}%
and argue exactly as above to get the following estimate:%
\begin{align}
& \frac{1}{2}\frac{d}{dt}\left\Vert U\left( t\right) \right\Vert _{\mathbb{X}%
^{2}}^{2}+\left( \nu -\varepsilon \right) \left\Vert \nabla u\left( t\right)
\right\Vert _{L^{p}\left( \Omega \right) }^{p}+c\left( \left\Vert u\left(
t\right) \right\Vert _{L^{r_{1}}\left( \Omega \right) }^{r_{1}}+\left\Vert
v\left( t\right) \right\Vert _{L^{r_{2}}\left( \Gamma ,b^{-1}dS\right)
}^{r_{2}}\right)  \label{dissi2} \\
& \leq \left\langle h_{1},u\left( t\right) \right\rangle _{2}+\left\langle
h_{2}/b,v\left( t\right) \right\rangle _{2,\Gamma }+c.  \notag
\end{align}%
The desired control of $U=U_{\epsilon }\left( t\right) $ in (\ref{est1})-(%
\ref{est2}) can be obtained immediately from a simple application of H\"{o}%
lder and Young inequalities on the terms on the right-hand side of (\ref%
{dissi2}). Thus the proof is the same as in Step 1-(i).

\textbf{Step 2}. It is obvious that%
\begin{equation}
\partial _{t}U_{\epsilon }+\mathcal{B}_{p,\epsilon }U_{\epsilon }+\mathcal{F}%
\left( U_{\epsilon }\right) =\mathcal{G}  \label{fform}
\end{equation}%
holds as an equality in $L^{s}\left( \left[ 0,T\right] ;\left( \mathbb{V}%
^{k,p}\right) ^{\ast }\right) $ (this is the same as equation (\ref{eq})
with $U$ and $a\left( \cdot \right) $ replaced by $U_{\epsilon }$ and $%
a_{\epsilon }\left( \cdot \right) ,$ respectively). From the estimates that
we deduced in Step 1, (i)-(ii), we see that there exists a subsequence $%
\left\{ U_{\epsilon }\right\} =\left\{ \binom{u_{\epsilon }}{v_{\epsilon }}%
\right\} $ (still denoted by $\left\{ U_{\epsilon }\right\} $), such that as
$\epsilon \rightarrow 0,$%
\begin{equation}
\begin{array}{l}
U_{\epsilon }\rightarrow U\text{ weakly star in }L^{\infty }\left( \left[ 0,T%
\right] ;\mathbb{X}^{2}\right) , \\
U_{\epsilon }\rightarrow U\text{ weakly in }L^{p}\left( \left[ 0,T\right] ;%
\mathbb{V}^{p}\right) , \\
\partial _{t}U_{\epsilon }\rightarrow \rightarrow \partial _{t}U\text{
weakly in }L^{s}\left( \left[ 0,T\right] ;\left( \mathbb{V}^{k,p}\right)
^{\ast }\right) .%
\end{array}
\label{2.27}
\end{equation}%
On the basis on standard interpolation and compact embedding results for
vector valued functions (see \cite{GW}), we also have%
\begin{equation}
U_{\epsilon }\rightarrow U\text{ strongly in }L^{p}\left( \left[ 0,T\right] ;%
\mathbb{X}^{2}\right) .  \label{2.30}
\end{equation}%
Clearly, $U\in C\left( \left[ 0,T\right] ;\left( \mathbb{V}^{k,p}\right)
^{\ast }\right) $. By refining in \eqref{2.30}, $u_{\epsilon }$ converges to
$u$ a.e. in $\Omega \times \left( 0,T\right) $ and $v_{\epsilon }$ converges
to $v$ a.e. in $\Gamma \times \left( 0,T\right) ,$ respectively. Then, by
means of known results in measure theory (see, e.g., \cite{T}), the
continuity of $f$, $g,$ and the convergence of \eqref{2.30} imply that $%
f\left( u_{\epsilon }\right) $ converges weakly to $f\left( u\right) $ in $%
L^{r_{1}^{^{\prime }}}(\Omega \times \left( 0,T\right) ).$ Moreover, $%
g\left( v_{\epsilon }\right) $ converges weakly to $g\left( v\right) $ in $%
L^{r_{2}^{^{\prime }}}(\Gamma \times \left( 0,T\right) ),$ and thus, $%
\mathcal{F}\left( U_{\epsilon }\right) $ converges weakly star to $\mathcal{F%
}\left( U\right) $ in $L^{s}\left( \left[ 0,T\right] ;\left( \mathbb{V}%
^{k,p}\right) ^{\ast }\right) $. Since $\mathcal{B}_{p,\epsilon }U_{\epsilon
}$ is bounded in $L^{q}\left( \left[ 0,T\right] ;\left( \mathbb{V}%
^{p}\right) ^{\ast }\right) $, we further see that%
\begin{equation}
\mathcal{B}_{p,\epsilon }U_{\epsilon }\rightarrow \Xi \text{ weakly star in }%
L^{q}\left( \left[ 0,T\right] ;\left( \mathbb{V}^{p}\right) ^{\ast }\right) ,
\label{2.32}
\end{equation}%
and thus weakly star in $L^{s}\left( \left[ 0,T\right] ;\left( \mathbb{V}%
^{k,p}\right) ^{\ast }\right) ,$ since $s\leq q$ and $k\geq 1$. We are now
ready to pass to the limit as $\epsilon \rightarrow 0$ in equation (\ref%
{fform}). We have%
\begin{equation}
\partial _{t}U+\Xi +\mathcal{F}\left( U\right) =\mathcal{G}  \label{equality}
\end{equation}%
as an equality in $L^{s}\left( \left[ 0,T\right] ;\left( \mathbb{V}%
^{k,p}\right) ^{\ast }\right) $. It remains to show $\Xi =\mathcal{B}_{p}U$,
which can be proved by a standard monotonicity argument exactly as in \cite[%
Theorem 2.6]{GW}. We leave the details to the interested reader.
\end{proof}

The following proposition is also immediate.

\begin{proposition}
\label{T1bis}Let $a\left( \cdot \right) $ satisfy (H1). In addition, assume
that $f$, $g\in C\left( \mathbb{R},\mathbb{R}\right) $ satisfy%
\begin{equation}
\begin{array}{ll}
c_{1}\left\vert y\right\vert ^{r_{1}}-c\leq f\left( y\right) y\leq
c_{2}\left\vert y\right\vert ^{r_{1}}+c, & \forall y\in {\mathbb{R}}\text{,}
\\
c_{3}\left\vert y\right\vert ^{r_{2}}-c\leq g\left( y\right) y\leq
c_{4}\left\vert y\right\vert ^{r_{2}}+c, & \forall y\in {\mathbb{R}}\text{,}%
\end{array}
\label{arbp}
\end{equation}%
for some appropriate positive constants and some $r_{1},r_{2}\geq 1$ such
that $\max \left( r_{1},r_{2}\right) \geq p$. Then, for any initial data $%
U_{0}\in \mathbb{X}^{2},$ there exists at least one (globally defined) weak
solution $U\left( t\right) $\ in the sense of Definition \ref{weak}.
\end{proposition}

\begin{proof}
In this case, both $f$ and $g$ are dissipative so that we do \emph{not} need
to exploit the validity of assumption (H4). Indeed, it follows from (\ref%
{2.7}) that
\begin{align}
& \frac{1}{2}\frac{d}{dt}\left\Vert U\left( t\right) \right\Vert _{\mathbb{X}%
^{2}}^{2}+\nu \left\Vert \nabla u\left( t\right) \right\Vert _{L^{p}\left(
\Omega \right) }^{p}+c_{1}\left\Vert u\left( t\right) \right\Vert
_{L^{r_{1}}\left( \Omega \right) }^{r_{1}}+c_{3}\left\Vert v\left( t\right)
\right\Vert _{L^{r_{2}}(\Gamma ,b^{-1}dS)}^{r_{2}} \\
& \leq \left\langle h_{1},u\left( t\right) \right\rangle _{2}+\left\langle
h_{2}/b,v\left( t\right) \right\rangle _{2,\Gamma }+c,  \notag
\end{align}%
which yields the desired control of $U\left( t\right) $ in the corresponding
spaces (\ref{est1})-(\ref{est2}) with relative ease (see \cite{GW}, for
further details). Thus, the proof is the same as in Theorem \ref{T1}.
\end{proof}

\begin{remark}
\label{imp}(i) If $U\left( t\right) $ is a weak solution of problem (\ref%
{1.11})-(\ref{1.12}), in the sense of Definition \ref{weak}, then clearly $%
U\left( t\right) \in C\left( \left[ 0,T\right] ;\left( \mathbb{V}%
^{k,p}\right) ^{\ast }\right) .$ Since by duality,%
\begin{equation*}
\mathbb{V}^{k,p}\subset \mathbb{X}^{2}=\left( \mathbb{X}^{2}\right) ^{\ast
}\subset \left( \mathbb{V}^{k,p}\right) ^{\ast },
\end{equation*}%
for any $p\in (p_{0},\infty )\cap \left( 1,\infty \right) ,$ $k\geq 1$, and
recalling that $U\left( t\right) \in L^{\infty }\left( \left[ 0,T\right] ;%
\mathbb{X}^{2}\right) ,$ it follows
\begin{equation*}
U\left( t\right) \in C_{w}\left( \left[ 0,T\right] ;\mathbb{X}^{2}\right)
\end{equation*}%
(see, e.g., \cite[Theorem II.1.7]{CV}). Therefore the initial value $U_{\mid
t=0}=U_{0}$ is meaningful when $U_{0}\in \mathbb{X}^{2}$. Finally, we note
that in general, the assumptions (H3)-(H4) alone do not ensure the
uniqueness of the Cauchy problem (\ref{1.11})-(\ref{1.12}).

(ii) We are also allowed to have equality in (\ref{maxi}), i.e., $\gamma
=\max \left( r_{2},q\left( r_{2}-1\right) \right) \leq r_{1}$ if $%
h_{2}\equiv 0$. Indeed, this follows once again from (\ref{2.7q}) and Lemma %
\ref{boundary}, which allows us to control surface integrals in terms of
volume integrals.
\end{remark}

\begin{proposition}
\label{ppp}Let the assumptions of Theorem \ref{T1} be satisfied. Then any
weak solution $U\left( t\right) =\binom{u\left( t\right) }{v\left( t\right) }
$ of (\ref{1.11})-(\ref{1.12}) belongs to $C\left( \left[ 0,T\right] ;%
\mathbb{X}^{2}\right) ,$ such that $\left\Vert U\left( t\right) \right\Vert
_{\mathbb{X}^{2}}^{2}$ is absolutely continuous on $\left[ 0,T\right] ,$ and%
\begin{equation}
\frac{1}{2}\frac{d}{dt}\left\Vert U\left( t\right) \right\Vert _{\mathbb{X}%
^{2}}^{2}=-\left\langle \mathcal{B}_{p}\left( U\left( t\right) \right)
,U\left( t\right) \right\rangle -\left\langle \mathcal{F}\left( U\left(
t\right) \right) -\mathcal{G},U\left( t\right) \right\rangle ,  \label{coo}
\end{equation}%
for almost all $t\in \left[ 0,T\right] .$ We will refer to (\ref{coo}) as
the \emph{energy identity} for the parabolic system (\ref{1.11})-(\ref{1.12}%
).
\end{proposition}

\begin{proof}
This statement follows from a generalization of a known interpolation result
(see, e.g., \cite{CV, LM, T}). Indeed, identifying the space $H^{\ast
}=\left( \mathbb{X}^{2}\right) ^{\ast }$ with $H=\mathbb{X}^{2},$ we have%
\begin{equation*}
V:=\mathbb{V}^{1,p}\subset \mathbb{X}^{2}\subset \left( \mathbb{V}%
^{1,p}\right) ^{\ast }\subseteq \left( \mathbb{V}^{k,p}\right) ^{\ast }=:W,
\end{equation*}%
for any $k\geq 1.$ Moreover, the following inclusion $E:=\mathbb{X}%
^{r_{1},r_{2}}\subseteq \mathbb{X}^{2}\subseteq E^{\ast }\subseteq W$ also
holds (indeed, the dual of $E$ is the space $E^{\ast }=\mathbb{X}%
^{r_{1}^{^{\prime }},r_{2}^{^{\prime }}},$ and by (\ref{kcond}), $\mathbb{V}%
^{k,p}\subset E$). By virtue of equation (\ref{eq}), any distributional
derivative $\partial _{t}U\left( t\right) $ from $\mathcal{D}^{^{\prime
}}\left( \left[ 0,T\right] ;W\right) $ can be represented as $\partial
_{t}U\left( t\right) =\mathcal{Z}_{1}\left( t\right) +\mathcal{Z}_{2}\left(
t\right) ,$ where%
\begin{equation*}
\mathcal{Z}_{1}\left( t\right) :=-\mathcal{B}_{p}\left( U\left( t\right)
\right) ,\text{ }\mathcal{Z}_{2}\left( t\right) :=-\mathcal{F}\left( U\left(
t\right) \right) +\mathcal{G}.
\end{equation*}%
According to (\ref{est1})-(\ref{est3}), $\mathcal{Z}_{1}\left( t\right) \in
L^{q}\left( \left[ 0,T\right] ;V^{\ast }\right) =\left( L^{p}\left( \left[
0,T\right] ;V\right) \right) ^{\ast },$ $q=p/\left( p-1\right) ,$ for%
\begin{equation}
U\left( t\right) \in L^{p}\left( \left[ 0,T\right] ;V\right) +\left(
L^{r_{1}}\left( \left[ 0,T\right] \times \Omega \right) \times
L^{r_{2}}\left( \left[ 0,T\right] \times \Gamma \right) \right) ,
\label{regu}
\end{equation}%
while from (\ref{est3}), $\mathcal{Z}_{2}\left( t\right) \in
L^{r_{1}^{^{\prime }}}\left( \left[ 0,T\right] \times \Omega \right) \times
L^{r_{2}^{^{\prime }}}\left( \left[ 0,T\right] \times \Gamma \right) ,$
which is precisely the dual of the product space in (\ref{regu}). Thus, the
claim follows, for instance, from \cite[Theorem II.1.8]{CV} (see also \cite[%
Proposition 23.23]{Ze}).
\end{proof}

We will now state some results which reflect the applicability of assumption
(H4) to a wide range of situations. In particular, it applies to the case of
competing nonlinearities $f$ and $g$, that is, nonlinearities with arbitrary
\emph{polynomial} growth which satisfy (H3), but when either one exhibits a
non-dissipative behavior at infinity. Recall that $p\in (\frac{2N}{N+2}%
,\infty )\cap \left( 1,\infty \right) $. In the interesting case of an
internal dissipation mechanism, and non-dissipative boundary conditions, we
have the following.

\begin{corollary}
\label{C1}Assume that $f,g\in C^{1}\left( \mathbb{R},\mathbb{R}\right) $
satisfy%
\begin{equation}
\lim_{\left\vert y\right\vert \rightarrow \infty }\frac{f^{^{\prime }}\left(
s\right) }{\left\vert s\right\vert ^{r_{1}-2}}=\left( r_{1}-1\right) c_{f}>0%
\text{ and }\lim_{\left\vert y\right\vert \rightarrow \infty }\frac{%
g^{^{\prime }}\left( s\right) }{\left\vert s\right\vert ^{r_{2}-2}}=\left(
r_{2}-1\right) c_{g}<0  \label{idnb}
\end{equation}%
with $r_{1}\geq p,$ $r_{2}>1.$ Suppose that one of the following conditions
holds:

(i) (\ref{maxi}) holds, i.e., $\max \left( r_{2},q\left( r_{2}-1\right)
\right) <r_{1}$.

(ii) $h_{2}=0$, $p<r_{2}<q\left( r_{2}-1\right) =r_{1}$ and%
\begin{equation}
c_{f}\nu p^{-q}q>\widetilde{C}_{\Omega }^{q}c_{g}^{q}\left( r_{2}\right)
^{q}.  \label{ccoef}
\end{equation}%
Then, in each case the conclusion of Theorem \ref{T1} applies.
\end{corollary}

\begin{proof}
We begin by noting that (H3a) is immediately satisfied. For sufficiently
large $y,$ we have%
\begin{equation*}
f\left( y\right) \sim c_{f}\left\vert y\right\vert ^{r_{1}-2}y,\text{ }%
g\left( y\right) \sim c_{g}\left\vert y\right\vert ^{r_{2}-2}y,
\end{equation*}%
and $f\left( y\right) y\sim c_{f}\left\vert y\right\vert ^{r_{1}},$ $g\left(
y\right) y\sim c_{g}\left\vert y\right\vert ^{r_{2}}.$ Thus, the leading
terms in (\ref{NB}) are%
\begin{equation}
c_{f}\left\vert y\right\vert ^{r_{1}}+\left( \left\vert \Omega \right\vert
\lambda \right) ^{-1}c_{g}\left\vert y\right\vert ^{r_{2}}-\frac{\widetilde{C%
}_{\Omega }^{q}}{\left( \varepsilon p\right) ^{q/p}q}c_{g}^{q}\left(
r_{2}\right) ^{q}\left\vert y\right\vert ^{q\left( r_{2}-1\right) },
\label{list}
\end{equation}%
for some $\varepsilon \in \left( 0,\nu /q\right) $. By assumption (i), it
holds $\gamma =\max \left( r_{2},q\left( r_{2}-1\right) \right) <r_{1},$ so
the coefficient of the highest order term in (\ref{list}) is $c_{f}$, which
is positive. If (ii) holds, it is obvious that $r_{2}<\max \left(
r_{2},q\left( r_{2}-1\right) \right) =r_{1},$ so the coefficient of the
highest order term in (\ref{list}) is%
\begin{equation*}
c_{f}-\frac{\widetilde{C}_{\Omega }^{q}}{\left( \varepsilon p\right) ^{q/p}q}%
c_{g}^{q}\left( r_{2}\right) ^{q},
\end{equation*}%
which is positive for some $\varepsilon \in \left( 0,\nu /q\right) ,$ if (%
\ref{ccoef}) is satisfied. Therefore, the conditions of Theorem \ref{T1} are
met, and the proof is finished.
\end{proof}

In the case of non-dissipative \emph{polynomial} behavior for $f$, we have
the following.

\begin{corollary}
\label{C1bis}Assume $h_{2}=0$ and $f,$ $g\in C^{1}\left( \mathbb{R},\mathbb{R%
}\right) $ satisfy%
\begin{equation}
\lim_{\left\vert y\right\vert \rightarrow \infty }\frac{f^{^{\prime }}\left(
s\right) }{\left\vert s\right\vert ^{p-2}}=\left( r_{1}-1\right) c_{f}<0%
\text{ and }\lim_{\left\vert y\right\vert \rightarrow \infty }\frac{%
g^{^{\prime }}\left( s\right) }{\left\vert s\right\vert ^{p-2}}=\left(
r_{2}-1\right) c_{g}>0  \label{inbd}
\end{equation}%
for some $p\in (\frac{2N}{N+2},\infty )\cap \left( 1,\infty \right) ,$ and
let%
\begin{equation}
\left( c_{f}+\left( \left\vert \Omega \right\vert \lambda \right)
^{-1}c_{g}\right) \nu p^{-q}q>\widetilde{C}_{\Omega }^{q}c_{g}^{q}\left(
r_{2}\right) ^{q}.  \label{ccoef2}
\end{equation}

Then, the conclusion of Theorem \ref{T1} applies.
\end{corollary}

\begin{proof}
First, it is immediate that (H3b) holds with $r_{1}=r_{2}=p$. Obviously, in
this case $r_{2}=r_{1}=q\left( r_{2}-1\right) .$ The coefficient of the
highest order term in (\ref{list}) is then%
\begin{equation*}
\left( c_{f}+\left( \left\vert \Omega \right\vert \lambda \right)
^{-1}c_{g}\right) -\frac{\widetilde{C}_{\Omega }^{q}}{\left( \varepsilon
p\right) ^{q/p}q}c_{g}^{q}\left( r_{2}\right) ^{q},
\end{equation*}%
which is positive for some $\varepsilon \in \left( 0,\nu /q\right) ,$ if (%
\ref{ccoef2}) is satisfied. Therefore, condition (H4) holds true and, thus,
the assumptions of Theorem \ref{T1} are verified.
\end{proof}

In the case of simultaneous internal and boundary dissipation, we can prove
the following stability result. Note that in this proposition, the
uniqueness holds in the class of all solutions which are constructed by
means of Definition \ref{weak}, and not only for solutions which can be
obtained as the limit, as $\epsilon \rightarrow 0,$ of the (strictly)
non-degenerate parabolic system introduced earlier.

\begin{proposition}
\label{uniq}Let the assumptions of Theorem \ref{T1} be satisfied, and in
addition, assume that (H2) holds. Then, there exists a unique weak solution
to problem (\ref{1.11})-(\ref{1.12}), which depends continuously on the
initial data in a Lipschitz way.
\end{proposition}

\begin{proof}
Let us consider two weak solutions $U_{1}=\binom{u_{1}}{u_{1\mid \Gamma }}$
and $U_{2}=\binom{u_{2}}{u_{2\mid \Gamma }}$, which according to Proposition %
\ref{ppp} belong to the spaces in (\ref{regu}), and set $U\left( t\right)
=U_{1}\left( t\right) -U_{2}\left( t\right) $. Thus, $U\left( t\right) $
satisfies the equation%
\begin{equation*}
\partial _{t}U+\left( \mathcal{B}_{p}\left( U_{1}\left( t\right) \right) -%
\mathcal{B}_{p}\left( U_{2}\left( t\right) \right) \right) +\mathcal{F}%
\left( U_{1}\left( t\right) \right) -\mathcal{F}\left( U_{2}\left( t\right)
\right) =0,
\end{equation*}%
and%
\begin{equation*}
U\left( 0\right) =U_{1}\left( 0\right) -U_{2}\left( 0\right) .
\end{equation*}%
Since%
\begin{equation*}
U\left( t\right) \in L^{p}\left( \left[ 0,T\right] ;\mathbb{V}^{p}\right)
\cap \left( L^{r_{1}}\left( \left[ 0,T\right] \times \Omega \right) \times
L^{r_{2}}\left( \left[ 0,T\right] \times \Gamma \right) \right) ,
\end{equation*}%
and%
\begin{equation*}
\partial _{t}U\left( t\right) \in L^{q}\left( \left[ 0,T\right] ;\left(
\mathbb{V}^{p}\right) ^{\ast }\right) +(L^{r_{1}^{^{\prime }}}\left( \left[
0,T\right] \times \Omega \right) \times L^{r_{2}^{^{\prime }}}\left( \left[
0,T\right] \times \Gamma \right) ),
\end{equation*}%
Proposition \ref{ppp} is indeed applicable, and we have%
\begin{eqnarray*}
&&\frac{1}{2}\frac{d}{dt}\left\Vert U\left( t\right) \right\Vert _{\mathbb{X}%
^{2}}^{2}+\left\langle \mathcal{B}_{p}\left( U_{1}\left( t\right) \right) -%
\mathcal{B}_{p}\left( U_{2}\left( t\right) \right) ,U\left( t\right)
\right\rangle \\
&=&-\left\langle \mathcal{F}\left( U_{1}\left( t\right) \right) -\mathcal{F}%
\left( U_{2}\left( t\right) \right) ,U\left( t\right) \right\rangle ,
\end{eqnarray*}%
for almost all $t\in \left[ 0,T\right] .$ Recalling that $\mathcal{B}%
_{p}\left( \cdot \right) \in \mathcal{L}\left( \mathbb{V}^{p},\left( \mathbb{%
V}^{p}\right) ^{\ast }\right) $ is monotone and coercive (see (\ref{deff});
cf. also Appendix), we get%
\begin{align}
& \frac{1}{2}\frac{d}{dt}\left\Vert U\left( t\right) \right\Vert _{\mathbb{X}%
^{2}}^{2}  \label{uniqqq} \\
& \leq -\left\langle f\left( u_{1}\left( t\right) \right) -f\left(
u_{2}\left( t\right) \right) ,u\left( t\right) \right\rangle -\left\langle
g\left( u_{1}\left( t\right) \right) -g\left( u_{2}\left( t\right) \right)
,u\left( t\right) \right\rangle .  \notag
\end{align}%
Exploiting assumption (H2) (which implies, $f^{^{\prime }}\left( y\right)
\geq -c_{f}$ and $g^{^{\prime }}\left( y\right) \geq -c_{g},$ $\forall y\in
\mathbb{R}$, for some $c_{f},$ $c_{g}>0$), we obtain%
\begin{equation*}
\frac{d}{dt}\left\Vert U\left( t\right) \right\Vert _{\mathbb{X}%
^{2}}^{2}\leq 2\left( c_{f}+c_{g}\right) \left\Vert U\left( t\right)
\right\Vert _{\mathbb{X}^{2}}^{2},
\end{equation*}%
for almost all $t\in \left[ 0,T\right] .$ Integrating this inequality over $%
\left[ 0,T\right] $ and applying Gronwall's inequality, we deduce%
\begin{equation}
\left\Vert U\left( t\right) \right\Vert _{\mathbb{X}^{2}}^{2}\leq
e^{ct}\left\Vert U\left( 0\right) \right\Vert _{\mathbb{X}^{2}}^{2},
\label{unc2}
\end{equation}%
which yields the desired result.
\end{proof}

\begin{remark}
(i) Note that assumption (H2) is only required to prove uniqueness of the
weak solution, and is usually not required for the theory of attractors.
Moreover, this assumption is actually too restrictive so that nonlinearities
that satisfy it are \emph{not} allowed to carry a bad sign at infinity, and
thus this would automatically eliminate the scenario proposed by the
statements of Corollaries \ref{C1} and \ref{C1bis}. Indeed, a simple
observation that will be made in Section 3 is that actually uniqueness is
necessary on the attractor only, and this can be obtained by deducing
additional regularity estimates for the solutions. This observation is in
particular very useful if one needs to consider entropy-related
nonlinearities of the form $f\left( y\right) =y^{l}\log \left( y\right) $,
for $y>0$, and $f\left( y\right) =0$, for $y\leq 0$, for some $l\geq 1.$

(ii) Obviously, estimate (\ref{unc2}) also holds if we assume that $%
\left\vert f^{^{\prime }}\left( y\right) \right\vert $ and $\left\vert
g^{^{\prime }}\left( y\right) \right\vert $ are bounded for all $y\in
\mathbb{R}$.
\end{remark}

As an immediate consequence of the stability result just proven above,
problem (\ref{1.11})-(\ref{1.12}), (\ref{ini}) defines a dynamical system in
the classical sense.

\begin{corollary}
Let the assumptions of Proposition \ref{uniq} be satisfied. The problem (\ref%
{1.11})-(\ref{1.12}), (\ref{ini}) defines a (nonlinear) continuous semigroup
$\mathcal{S}_{2}\left( t\right) $ on the phase space $\mathbb{X}^{2}$,
\begin{equation*}
\mathcal{S}_{2}\left( t\right) :\mathbb{X}^{2}\rightarrow \mathbb{X}^{2},
\end{equation*}%
given by
\begin{equation}
\mathcal{S}_{2}\left( t\right) U_{0}=U\left( t\right) ,  \label{2.15}
\end{equation}%
where $U\left( t\right) $ is the (unique)\ weak solution which satisfies the
energy identity (\ref{coo}).
\end{corollary}

\section{Global Attractors}

\label{attractor}

\subsection{Attractors for $\left( \mathcal{S}_{2}\left( t\right) ,\mathbb{X}%
^{2}\right) $\ revisited}

In order to study the asymptotic behavior of (\ref{1.11})-(\ref{1.12}), (\ref%
{ini}), we need to derive some additional apriori estimates for the
solutions. We shall focus our study on the case $p\geq 2$ only, since for $%
p\in (\frac{2N}{N+2},2)$ we need to impose slightly different assumptions on
the nonlinearities, and so we will pursue this question elsewhere. We first
aim to improve some results from \cite{GW} for the weak solutions
constructed in Definition \ref{weak}, which are unique by Proposition \ref%
{uniq}, and to show the existence of the (classical) global attractor,
bounded in $\mathbb{X}^{\infty }\cap \mathbb{V}^{p}$. We emphasize again
that all the results below hold for any $p\geq 2$.

The next result is a direct consequence of estimate (\ref{dissi}) of Theorem %
\ref{T1} (see \cite[Section 2, Proposition 3.3]{GW}, for details).

\begin{proposition}
\label{P1}Let the assumptions of either Theorem \ref{T1} or Proposition \ref%
{T1bis} be satisfied. The solution semigroup $\left\{ \mathcal{S}_{2}\left(
t\right) \right\} _{t\geq 0}$ has a $\left( \mathbb{X}^{2},\mathbb{X}%
^{2}\right) $-bounded absorbing set. More precisely, there is a positive
constant $C_{0},$ depending only on the physical parameters of the problem,
such that for any bounded subset $B\subset \mathbb{X}^{2}$, there exists a
positive constant $t^{\#}=t^{\#}\left( \left\Vert B\right\Vert _{\mathbb{X}%
^{2}}\right) $ such that
\begin{equation}
\sup_{t\geq t^{\#}}\left[ \left\Vert U\left( t\right) \right\Vert _{\mathbb{X%
}^{2}}+\int_{t}^{t+1}\left( \int_{\Omega }a(\left\vert \nabla u\left(
s\right) \right\vert ^{2})\left\vert \nabla u\left( s\right) \right\vert
^{2}+\left\vert u\left( s\right) \right\vert ^{r_{1}}\right) dxds\right]
\leq C_{0}.  \label{3.2}
\end{equation}
\end{proposition}

Our next goal is to establish the existence of a bounded absorbing set in $%
\mathbb{X}^{\infty }$, which has an interest on its own. The following
result extends \cite[Theorem 3.7]{GW} by removing the additional growth
conditions that were imposed on $f,g$ in \cite[Theorem 3.7, (3.17)]{GW}.

\begin{theorem}
\label{linff}Let the assumptions of Proposition \ref{P1} be satisfied. Let $%
h_{1}\in L^{\infty }\left( \Omega \right) $, $h_{2}\in L^{\infty }\left(
\Gamma \right) $, and suppose%
\begin{equation}
\lim_{\left\vert y\right\vert \rightarrow \infty }\inf \frac{f\left(
y\right) }{y}>0,\text{ }\lim_{\left\vert y\right\vert \rightarrow \infty
}\inf \frac{g\left( y\right) }{y}>0.  \label{weaker2}
\end{equation}%
Then, given any initial data $U_{0}$\ in $\mathbb{X}^{2}$, the corresponding
solution $U\left( t\right) $ of (\ref{1.11})-(\ref{1.12}), (\ref{ini})
belongs to $\mathbb{X}^{\infty }$, for each $t>0$. Moreover, there exists a
positive constant $C_{1},$ independent of $t$ and the initial data, and a
positive constant $t_{+}$ depending on $t^{\#}$, such that%
\begin{equation}
\sup_{t\geq t_{+}}\left\Vert U\left( t\right) \right\Vert _{\mathbb{X}%
^{\infty }}\leq C_{1}.  \label{3.3}
\end{equation}
\end{theorem}

\begin{proof}
All the calculations below are formal. However, they can be rigorously
justified by means of the approximation procedure devised in Section 2 (see (%
\ref{ap})-(\ref{ini2})). From now on, $c$ will denote a positive constant
that is independent of $t,$ $\epsilon $, $m$ and initial data, which only
depends on the other structural parameters of the problem. Such a constant
may vary even from line to line. Moreover, we shall denote by $Q_{\tau
}\left( m\right) $ a monotone nondecreasing function in $m$ of order $\tau ,$
for some nonnegative constant $\tau ,$ independent of $m.$ More precisely, $%
Q_{\tau }\left( m\right) \sim cm^{\tau }$ as $m\rightarrow +\infty .$

We begin by showing that the $\mathbb{X}^{m}$-norm of $U$ satisfies a local
recursive relation which can be used to perform an iterative argument. We
divide the proof of (\ref{3.3}) into several steps.

\noindent \textbf{Step 1} (The basic energy estimate in $\mathbb{X}^{m+1}$).
We multiply (\ref{1.11}) by $\left\vert u\right\vert ^{m-1}u,$ $m\geq 1,$
and integrate over $\Omega $. We obtain%
\begin{eqnarray}
&&\frac{1}{\left( m+1\right) }\frac{d}{dt}\left\Vert u\right\Vert
_{m+1}^{m+1}+\left\langle f\left( u\right) ,\left\vert u\right\vert
^{m-1}u\right\rangle _{2}+m\int_{\Omega }a\left( \left\vert \nabla
u\right\vert ^{2}\right) \left\vert \nabla u\right\vert ^{2}\left\vert
u\right\vert ^{m-1}dx  \label{eqn2} \\
&=&\int_{\Gamma }a\left( \left\vert \nabla u\right\vert ^{2}\right) \partial
_{\mathbf{n}}u\left\vert v\right\vert ^{m-1}vdS+\left\langle h_{1}\left(
x\right) ,\left\vert u\right\vert ^{m-1}u\right\rangle _{2}.  \notag
\end{eqnarray}%
Similarly, we multiply (\ref{1.12}) by $\left\vert v\right\vert ^{m-1}v/b$
and integrate over $\Gamma $. We have%
\begin{align}
& \frac{1}{\left( m+1\right) }\frac{d}{dt}\left\Vert \frac{v}{b}\right\Vert
_{m+1,\Gamma }^{m+1}+\int_{\Gamma }b\left( x\right) a\left( \left\vert
\nabla u\right\vert ^{2}\right) \partial _{\mathbf{n}}u\left\vert
v\right\vert ^{m-1}v\frac{dS}{b\left( x\right) }+\left\langle g\left(
v\right) ,\frac{\left\vert v\right\vert ^{m-1}v}{b}\right\rangle _{2,\Gamma }
\label{eqn3} \\
& =\left\langle h_{2}\left( x\right) ,\frac{\left\vert v\right\vert ^{m-1}v}{%
b}\right\rangle _{2,\Gamma }.  \notag
\end{align}%
Let us first observe that, in light of assumption (\ref{1.1.2}), it is easy
to check%
\begin{equation}
m\int_{\Omega }a\left( \left\vert \nabla u\right\vert ^{2}\right) \left\vert
\nabla u\right\vert ^{2}\left\vert u\right\vert ^{m-1}dx\geq \nu m\left(
\frac{p}{p+m-1}\right) ^{p}\int_{\Omega }\left\vert \nabla \left\vert
u\right\vert ^{\frac{p+m-1}{p}}\right\vert ^{p}dx.  \label{eqn3b}
\end{equation}%
Adding relations (\ref{eqn2})-(\ref{eqn3}), we deduce on account of the
assumptions (\ref{weaker2}) (indeed, it holds%
\begin{equation}
f\left( y\right) y\geq -c_{1}s^{2}-c_{2},g\left( y\right) y\geq
-c_{3}s^{2}-c_{4},  \label{weakly}
\end{equation}%
for all $y\in \mathbb{R}$, and some $c_{i}>0$) and an application of basic H%
\"{o}lder and Young inequalities, the following inequality%
\begin{equation}
\frac{d}{dt}\left\Vert U\right\Vert _{\mathbb{X}^{m+1}}^{m+1}+\gamma \nu
\int_{\Omega }\left\vert \nabla \left\vert u\right\vert ^{\frac{p+m-1}{p}%
}\right\vert ^{p}dx\leq Q_{1}\left( m\right) \left( \left\Vert U\right\Vert
_{\mathbb{X}^{m+1}}^{m+1}+1\right) .  \label{eq4}
\end{equation}%
Here the positive constant $\gamma \sim m^{-\left( p-2\right) }$, $p\geq 2$,
and the function $Q_{1}\left( m\right) \sim m$ depends also on the $%
L^{\infty }$-norms of $h_{1},$ $h_{2},$ and of $b_{0}$.

\noindent \textbf{Step 2} (The local relation). Set $m_{k}=p^{k},$ and define%
\begin{equation}
\mathcal{Y}_{k}\left( t\right) :=\int_{\Omega }\left\vert u\left( t,\cdot
\right) \right\vert ^{1+m_{k}}dx+\int_{\Gamma }\left\vert v\left( t,\cdot
\right) \right\vert ^{1+m_{k}}\frac{dS}{b}=\left\Vert U\left( t\right)
\right\Vert _{\mathbb{X}^{m_{k}+1}}^{m_{k}+1},  \label{def}
\end{equation}%
for all $k\geq 0$. Let $t,\mu $ be two positive constants such that $t-\mu
/m_{k}>0$. Their precise values will be chosen later. We claim that%
\begin{equation}
\mathcal{Y}_{k}\left( t\right) \leq M_{k}\left( t,\mu \right) :=c\left(
m^{k}\right) ^{\sigma }(\sup_{s\geq t-\mu /m_{k}}\mathcal{Y}_{k-1}\left(
s\right) +1)^{n_{k}},\text{ }\forall k\geq 1,  \label{claim2}
\end{equation}%
where $c,$ $\sigma $ are positive constants independent of $k,$ and $%
n_{k}:=\max \left\{ z_{k},l_{k}\right\} \geq 1$ is a bounded sequence for
all $k.$

We will now prove (\ref{claim2}) when $p<N.$ The case $p\geq N$ shall
require only minor adjustments (in fact, in this case we can choose any
arbitrary, but fixed, $p_{s},q_{s}>p$ in the embedding $\mathbb{V}%
^{1,p}\left( \Omega \right) \subset \mathbb{X}^{p_{s},q_{s}}$). For each $%
k\geq 0$, we define%
\begin{equation*}
r_{k}:=\frac{N\left( p+m_{k}-1\right) -\left( N-p\right) \left(
1+m_{k}\right) }{N\left( p+m_{k}-1\right) -\left( N-p\right) \left(
1+m_{k-1}\right) },\text{ }s_{k}:=1-r_{k}.
\end{equation*}%
We aim to estimate the term on the right-hand side of (\ref{eq4}) in terms
of the $\mathbb{X}^{1+m_{k-1}}$-norm of $U.$ First, H\"{o}lder and Sobolev
inequalities (with the equivalent norm of Sobolev spaces in $W^{1,p}\left(
\Omega \right) \subset L^{p_{s}}\left( \Omega \right) $, $p_{s}=pN/\left(
N-p\right) $) yield%
\begin{eqnarray}
\int_{\Omega }\left\vert u\right\vert ^{1+m_{k}}dx &\leq &\left(
\int_{\Omega }\left\vert u\right\vert ^{\frac{\left( p+m_{k}-1\right) N}{N-p}%
}dx\right) ^{s_{k}}\left( \int_{\Omega }\left\vert u\right\vert
^{1+m_{k-1}}dx\right) ^{r_{k}}  \label{ee5} \\
&\leq &c\left( \int_{\Omega }\left\vert \nabla \left\vert u\right\vert ^{%
\frac{\left( p+m_{k}-1\right) }{p}}\right\vert ^{p}dx+\left( \int_{\Omega
}\left\vert u\right\vert ^{1+m_{k-1}}dx\right) ^{\alpha _{k}}\right) ^{%
\overline{s}_{k}}\left( \int_{\Omega }\left\vert u\right\vert
^{1+m_{k-1}}dx\right) ^{r_{k}},  \notag
\end{eqnarray}%
with%
\begin{align*}
\overline{s}_{k}& :=s_{k}\frac{N}{N-p}=\frac{\left( p-1\right) Nm_{k-1}}{%
\left( p-1\right) Nm_{k-1}+m_{k}+p_{1}}\in \left( 0,1\right) , \\
p_{1}& :=\left( p-1\right) N-\left( N-p\right) >0,\text{ }\alpha _{k}:=\frac{%
p+m_{k}-1}{1+m_{k-1}}\in \left[ 1,p\right] .
\end{align*}%
Applying Young's inequality on the right-hand side of (\ref{ee5}), we get%
\begin{equation}
Q_{\tau _{1}}\left( m_{k}\right) \int_{\Omega }\left\vert u\right\vert
^{1+m_{k}}dx\leq \frac{\gamma _{k}}{4}\int_{\Omega }\left\vert \nabla
\left\vert u\right\vert ^{\frac{p+m_{k}-1}{p}}\right\vert ^{p}dx+Q_{\tau
_{2}}\left( m_{k}\right) \left( \int_{\Omega }\left\vert u\right\vert
^{1+m_{k-1}}dx\right) ^{\max \left\{ z_{k},\alpha _{k}\right\} },
\label{ee5bis}
\end{equation}%
for some positive constants $\tau _{i}$ independent of $m$, and where%
\begin{equation*}
z_{k}:=r_{k}/\left( 1-\overline{s}_{k}\right) =\left( p_{1}+m_{k+1}\right)
/\left( p_{1}+m_{k}\right) \geq 1
\end{equation*}%
is bounded for all $k$. Note that we can choose $\tau _{2}$ to be some fixed
positive number since $Q_{\tau _{2}}$ also depends on $\gamma _{k}^{-1}\sim
m_{k}^{p-2}$. To treat the boundary terms in (\ref{eq4}), we define for $%
k\geq 1,$%
\begin{equation*}
y_{k}:=\frac{\left( N-1\right) \left( p+m_{k}-1\right) -\left( N-p\right)
\left( 1+m_{k}\right) }{\left( N-1\right) \left( p+m_{k}-1\right) -\left(
N-p\right) \left( 1+m_{k-1}\right) },\text{ }x_{k}:=1-y_{k}.
\end{equation*}%
On account of H\"{o}lder and Sobolev inequalities (see Section 2), we obtain%
\begin{align}
\int_{\Gamma }\left\vert v\right\vert ^{1+m_{k}}\frac{dS}{b}& \leq c\left(
\int_{\Gamma }\left\vert v\right\vert ^{\frac{\left( N-1\right) \left(
p+m_{k}-1\right) }{N-p}}dS\right) ^{x_{k}}\left( \int_{\Gamma }\left\vert
v\right\vert ^{1+m_{k-1}}\frac{dS}{b}\right) ^{y_{k}}  \label{ee6} \\
& \leq c\left( \int_{\Omega }\left\vert \nabla \left\vert u\right\vert ^{%
\frac{\left( p+m_{k}-1\right) }{p}}\right\vert ^{p}dx+\left( \int_{\Omega
}\left\vert u\right\vert ^{1+m_{k-1}}dx\right) ^{\alpha _{k}}\right) ^{%
\overline{x}_{k}}\left( \int_{\Gamma }\left\vert v\right\vert ^{1+m_{k-1}}%
\frac{dS}{b}\right) ^{y_{k}},  \notag
\end{align}%
with%
\begin{align*}
\overline{x}_{k}& :=\frac{N-1}{N-p}x_{k}=\frac{\left( N-1\right) \left(
p-1\right) m_{k-1}}{\left( N-1\right) \left( p-1\right) m_{k-1}+\left(
p-1\right) m_{k-1}+p_{2}}, \\
p_{2}& :=\left( N-1\right) \left( p-1\right) -\left( N-p\right) >0.
\end{align*}%
Since $\overline{x}_{k}\in \left( 0,1\right) $, we can apply Young's
inequality on the right-hand side of (\ref{ee6}), use the estimate for the $%
L^{1+m_{k}}\left( \Omega \right) $-norm of $u$ from (\ref{ee5bis}) in order
to deduce the following estimate:%
\begin{align}
& Q_{\tau _{3}}\left( m_{k}\right) \int_{\Gamma }\left\vert v\right\vert
^{1+m_{k}}\frac{dS}{b}  \label{ee6bis} \\
& \leq \frac{\gamma _{k}}{4}\int_{\Omega }\left\vert \nabla \left\vert
u\right\vert ^{\frac{p+m_{k}-1}{p}}\right\vert ^{p}dx+Q_{\tau _{4}}\left(
m_{k}\right) \left( \int_{\Omega }\left\vert u\right\vert
^{1+m_{k-1}}dx\right) ^{\max \left\{ l_{k},\alpha _{k}\right\} },  \notag
\end{align}%
for some positive constants $\tau _{3},\tau _{4}$ depending on $\tau
_{1},\tau _{2},$ but which are independent of $m$. The sequence%
\begin{equation*}
l_{k}:=\frac{y_{k}}{\left( 1-\overline{x}_{k}\right) }=\frac{\left(
p-1\right) m_{k+1}+pp_{2}}{\left( p-1\right) m_{k}+pp_{2}}\geq 1
\end{equation*}%
is bounded for all $k\geq 1$. Inserting estimates (\ref{ee5bis})-(\ref%
{ee6bis}) on the right-hand side of (\ref{eq4}), we obtain the following
inequality:%
\begin{equation}
\partial _{t}\mathcal{Y}_{k}+\gamma _{k}\int_{\Omega }\left\vert \nabla
\left\vert u\right\vert ^{\frac{p+m_{k}-1}{p}}\right\vert ^{p}dx\leq c\left(
m_{k}\right) ^{\sigma _{1}}\left( \mathcal{Y}_{k-1}+1\right) ^{n_{k}},
\label{e11}
\end{equation}%
for some positive constant $\sigma _{1}$ that depends on $\tau _{i}$; we
recall that $n_{k}=\max \left\{ z_{k},l_{k},\alpha _{k}\right\} \geq 1$, and
$\gamma _{k}\sim m_{k}^{2-p}$.

Let now $\zeta \left( s\right) $ be a positive function $\zeta :\mathbb{R}%
_{+}\rightarrow \left[ 0,1\right] $ such that $\zeta \left( s\right) =0$ for
$s\in \left[ 0,t-\mu /r_{k}\right] ,$ $\zeta \left( s\right) =1$ if $s\in %
\left[ t,+\infty \right) $ and $\left\vert d\zeta /ds\right\vert \leq
m_{k}/\mu $, if $s\in \left( t-\mu /r_{k},t\right) $. We define $Z_{k}\left(
s\right) =\zeta \left( s\right) \mathcal{Y}_{k}\left( s\right) $ and notice
that%
\begin{align}
\frac{d}{ds}Z_{k}\left( s\right) & \leq \zeta \left( s\right) \frac{d}{ds}%
\mathcal{Y}_{k}\left( s\right) +\frac{m_{k}}{\mu }\mathcal{Y}_{k}\left(
s\right)  \label{e11bis} \\
& =\zeta \left( s\right) \frac{d}{ds}\mathcal{Y}_{k}\left( s\right)
+Q_{1}\left( m_{k}\right) \left( \int_{\Omega }\left\vert u\right\vert
^{1+m_{k}}dx+\int_{\Gamma }\left\vert v\right\vert ^{1+m_{k}}\frac{dS}{b}%
\right) .  \notag
\end{align}%
The last two integrals in (\ref{e11bis}) can be estimated as in (\ref{ee5bis}%
) and (\ref{ee6bis}). Combining the above estimates and the fact that $%
Z_{k}\leq \mathcal{Y}_{k}$, we deduce the following inequality:%
\begin{equation}
\frac{d}{ds}Z_{k}\left( s\right) +cm_{k}Z_{k}\left( s\right) \leq
M_{k}\left( t,\mu \right) ,\text{ for all }s\in \left[ t-\mu /r_{k},+\infty
\right) .  \label{e12}
\end{equation}%
Note that $c=c\left( \mu \right) \sim \mu ^{-1}$ as $\mu \rightarrow 0$, and
$c\left( \mu \right) $ is bounded if $\mu $ is bounded away from zero.
Integrating (\ref{e12}) with respect to $s$ from $t-\mu /r_{k}$ to $t,$ and
taking into account the fact that $Z_{k}\left( t-\mu /r_{k}\right) =0,$ we
obtain that $\mathcal{Y}_{k}\left( t\right) =Z_{k}\left( t\right) \leq
M_{k}\left( t,\mu \right) \left( 1-e^{-c\mu }\right) $, which proves the
claim (\ref{claim2}).

\noindent \textbf{Step 3} (The iterative argument). Let now $\tau ^{^{\prime
}}>\tau >0$ be given with $\tau =t^{\#}$ as in (\ref{3.2}), and define $\mu
=p(\tau ^{^{\prime }}-\tau )\geq 1,$ $t_{0}=\tau ^{^{\prime }}=t^{\#}+1$ and
$t_{k}=t_{k-1}-\mu /m_{k},$ $k\geq 1$. Using (\ref{claim2}), we have%
\begin{equation}
\sup_{t\geq t_{k-1}}\mathcal{Y}_{k}\left( t\right) \leq c\left( \mu \right)
\left( m_{k}\right) ^{\sigma }(\sup_{s\geq t_{k}}\mathcal{Y}_{k-1}\left(
s\right) +1)^{n_{k}},\text{ }k\geq 1.  \label{e13}
\end{equation}%
Note that from (\ref{3.2}), we have%
\begin{equation}
\sup_{s\geq t_{1}=\tau }\left( \mathcal{Y}_{0}\left( s\right) +1\right) \leq
C_{0}+1=:\overline{C},  \label{e13bis}
\end{equation}%
and $c=c\left( \mu \right) $ is bounded away from zero. Thus, we can iterate
in (\ref{e13}) with respect to $k\geq 1$ and obtain that%
\begin{align}
\sup_{t\geq t_{k-1}}\mathcal{Y}_{k}\left( t\right) & \leq \left(
cm_{k}^{\sigma }\right) \left( cm_{k-1}^{\sigma }\right) ^{n_{k}}\left(
cm_{k-2}^{\sigma }\right) ^{n_{k}n_{k-1}}\cdot ...\cdot \left(
cm_{0}^{\sigma }\right) ^{n_{k}n_{k-1}...n_{0}}(\overline{C})^{\xi _{k}}
\label{e13tris} \\
& \leq c^{A_{k}}p^{\sigma B_{k}}\left( \overline{C}\right) ^{\xi _{k}},
\notag
\end{align}%
where $\xi _{k}:=n_{k}n_{k-1}...n_{0},$ and%
\begin{equation}
A_{k}:=1+n_{k}+n_{k}n_{k-1}+...+n_{k}n_{k-1}...n_{0},  \label{ak2}
\end{equation}%
\begin{equation}
B_{k}:=k+n_{k}\left( k-1\right) +n_{k}n_{k-1}\left( k-2\right)
+...+n_{k}n_{k-1}...n_{0}.  \label{bk2}
\end{equation}%
Without loss of generality, let us assume that $z_{k}\geq l_{k}\geq 1,$ for
each $k$. Then, $n_{k}=z_{k},$ and $\xi _{k}=\left( p_{1}+m_{k+1}\right)
/\left( p_{1}+m_{0}\right) $. The argument below also applies to the case
when $n_{k}=l_{k}$. Thus, we have%
\begin{equation}
A_{k}\leq \left( p_{1}+m_{k}\right) \sum_{i=1}^{\infty }\frac{1}{p_{1}+m_{i}}%
\text{ and }B_{k}\leq \left( p_{1}+m_{k}\right) \sum_{i=1}^{\infty }\frac{i}{%
p_{1}+m_{i}}.  \label{abk}
\end{equation}%
Therefore, since%
\begin{equation}
\sup_{t\geq t_{0}}\mathcal{Y}_{k}\left( t\right) \leq \sup_{t\geq t_{k-1}}%
\mathcal{Y}_{k}\left( t\right) \leq c^{A_{k}}p^{\sigma B_{k}}\left(
\overline{C}\right) ^{\xi _{k}}  \label{e14}
\end{equation}%
and the series in (\ref{abk})\ are convergent, we can take the $1+m_{k}$%
-root on both sides of (\ref{e14}) and let $k\rightarrow +\infty $. We deduce%
\begin{equation}
\sup_{t\geq t_{0}=\tau ^{^{\prime }}}\left\Vert U\left( t\right) \right\Vert
_{\mathbb{X}^{\infty }}\leq \lim_{k\rightarrow +\infty }\sup_{t\geq
t_{0}}\left( \mathcal{Y}_{k}\left( t\right) \right) ^{1/\left(
1+m_{k}\right) }\leq C_{1},  \label{linf}
\end{equation}%
for some positive constant $C_{1}$ independent of $t,$ $k$, $U,$ $\epsilon $
and initial data. The proof of Theorem \ref{linff} is now complete.
\end{proof}

\begin{remark}
\label{rem}(i) We can easily modify our argument in the proof of (\ref{3.3})
in order to show that the $\mathbb{X}^{\infty }$-norm of the solution $%
U\left( t\right) $ stays bounded for all time $t\geq 0,$ if $U_{0}$ is
bounded in the $\mathbb{X}^{\infty }$-norm. It suffices to note that in
place of the inequality (\ref{claim2}), we may use instead the inequality%
\begin{equation*}
\mathcal{Y}_{k}\left( t\right) \leq Q(\left\Vert U_{0}\right\Vert _{\mathbb{X%
}^{\infty }},\sup_{t>0}M_{k}\left( t,\mu \right) ),
\end{equation*}%
which is an immediate consequence of (\ref{e11}). Then arguing as in the
proof of Theorem \ref{linff}, we also have the estimate:%
\begin{equation}
\sup_{t\geq 0}\left\Vert U\left( t\right) \right\Vert _{\mathbb{X}^{\infty
}}\leq Q(\left\Vert U_{0}\right\Vert _{\mathbb{X}^{\infty }},\sup_{t\geq
0}\left\Vert U\left( t\right) \right\Vert _{\mathbb{X}^{2}}),
\label{linfbis}
\end{equation}%
for some positive monotone nondecreasing (in each of its variables)\
function $Q:\mathbb{R}^{2}\rightarrow \mathbb{R}_{+}$ independent of $%
\epsilon .$

(ii) By slightly refining the arguments in the proof of Theorem \ref{linff}
(in Step 3), it is also easy to show that, for each $\tau >0,$%
\begin{equation*}
\sup_{t\geq 2\tau }\left\Vert U\left( t\right) \right\Vert _{\mathbb{X}%
^{\infty }}\leq Q(\tau ^{-1},\sup_{t\geq \tau }\left\Vert U\left( t\right)
\right\Vert _{\mathbb{X}^{2}}).
\end{equation*}
\end{remark}

Arguing exactly as in the proof of Theorem \ref{T1}, we can also obtain the
following \emph{general} balance condition between the functions $f$,$g$,
implying boundedness of the solution.

\begin{proposition}
\label{linff2}Let the assumptions of Theorem \ref{T1} be satisfied, and let $%
h_{1}\in L^{\infty }\left( \Omega \right) $, $h_{2}\in L^{\infty }\left(
\Gamma \right) $. Suppose that there exist $\tau \geq 0$ and $y_{0}>0$, such
that for any $m\geq 1$ and $\left\vert y\right\vert \geq y_{0},$ it holds%
\begin{align}
& f\left( y\right) \left\vert y\right\vert ^{m-1}y+\left( \left\vert \Omega
\right\vert \lambda \right) ^{-1}g\left( y\right) \left\vert y\right\vert
^{m-1}y-\frac{\widetilde{C}_{\Omega }^{q}m^{-q/p}}{\left( \varepsilon
p\right) ^{q/p}q}\left\vert y\right\vert ^{m-1}\left\vert g^{^{\prime
}}\left( y\right) y+mg\left( y\right) \right\vert ^{q}  \label{weaker} \\
& \geq -Q_{\tau }\left( m\right) (\left\vert y\right\vert ^{m+1}+1),  \notag
\end{align}%
for some $\varepsilon \in (0,\frac{\nu }{q}),$ and some positive function $%
Q_{\tau }:\mathbb{R}_{+}\rightarrow \mathbb{R}_{+},$ $Q_{\tau }\left(
m\right) \sim m^{\tau },$ as $m\rightarrow \infty $. Then, the same
conclusion of Theorem \ref{linff} applies to any weak solution of problem (%
\ref{1.11})-(\ref{1.12}), (\ref{ini}).
\end{proposition}

\begin{proof}
Let us return to equations (\ref{eqn2})-(\ref{eqn3b}). We can write%
\begin{align}
& \left\langle f\left( u\right) ,\left\vert u\right\vert
^{m-1}u\right\rangle _{2}+\left\langle g\left( v\right) ,\frac{\left\vert
v\right\vert ^{m-1}v}{b}\right\rangle _{2,\Gamma }  \label{li1} \\
& =\left\langle f\left( u\right) +c_{B}g\left( u\right) ,\left\vert
u\right\vert ^{m-1}u\right\rangle _{2}  \notag \\
& -c_{B}\left\langle g\left( u\right) u\left\vert u\right\vert
^{m-1}-\lambda \left\langle g\left( u\right) u\left\vert u\right\vert
^{m-1},1/b\right\rangle _{1,\Gamma }\right\rangle _{2},  \notag
\end{align}%
where $c_{B}$ and $\lambda $ are as in Theorem \ref{T1}. Applying the
Poincare's inequality (\ref{PI}), we have%
\begin{align*}
& c_{B}\left\vert \left\langle g\left( u\right) u\left\vert u\right\vert
^{m-1}-\lambda \left\langle g\left( u\right) u\left\vert u\right\vert
^{m-1},1/b\right\rangle _{1,\Gamma }\right\rangle _{2}\right\vert \\
& \leq \widetilde{C}\left\Vert \nabla \left( g\left( u\right) u\left\vert
u\right\vert ^{m-1}\right) \right\Vert _{L^{1}\left( \Omega \right) } \\
& =\widetilde{C}\left\Vert \left\vert u\right\vert ^{m-1}\nabla u\left(
g^{^{\prime }}\left( u\right) u+mg\left( u\right) \right) \right\Vert
_{L^{1}\left( \Omega \right) } \\
& =\widetilde{C}\int_{\Omega }\left\vert \left( \left\vert u\right\vert ^{%
\frac{m-1}{p}}\nabla u\right) \left\vert u\right\vert ^{\frac{m-1}{q}}\left(
g^{^{\prime }}\left( u\right) u+mg\left( u\right) \right) \right\vert dx.
\end{align*}%
On account of standard H\"{o}lder and Young inequalities, we can estimate
the term on the right-hand side in terms of%
\begin{align}
& \widetilde{C}\left( \int_{\Omega }\left\vert u\right\vert ^{m-1}\left\vert
\nabla u\right\vert ^{p}dx\right) ^{1/p}\left( \int_{\Omega }\left\vert
u\right\vert ^{m-1}\left\vert g^{^{\prime }}\left( u\right) u+mg\left(
u\right) \right\vert ^{q}dx\right) ^{1/q}  \label{li3} \\
& =\widetilde{C}\left( \frac{p}{p+m-1}\right) \left( m\int_{\Omega
}\left\vert \nabla \left\vert u\right\vert ^{\frac{p+m-1}{p}}\right\vert
^{p}dx\right) ^{1/p}  \notag \\
& \times \left( \int_{\Omega }\left\vert u\right\vert ^{m-1}\left\vert
g^{^{\prime }}\left( u\right) u+mg\left( u\right) \right\vert ^{q}dx\right)
^{1/q}m^{-1/p}  \notag \\
& \leq \varepsilon m\left( \frac{p}{p+m-1}\right) ^{p}\int_{\Omega
}\left\vert \nabla \left\vert u\right\vert ^{\frac{p+m-1}{p}}\right\vert
^{p}dx  \notag \\
& +\frac{\widetilde{C}^{q}m^{-q/p}}{\left( \varepsilon p\right) ^{q/p}q}%
\int_{\Omega }\left\vert u\right\vert ^{m-1}\left\vert g^{^{\prime }}\left(
u\right) u+mg\left( u\right) \right\vert ^{q}dx.  \notag
\end{align}%
Recalling (\ref{li1}), on the basis of (\ref{li3}), we can estimate%
\begin{align}
& \left\langle f\left( u\right) ,\left\vert u\right\vert
^{m-1}u\right\rangle _{2}+\left\langle g\left( v\right) ,\frac{\left\vert
v\right\vert ^{m-1}v}{b}\right\rangle _{2,\Gamma }  \label{li4} \\
& \geq \left\langle f\left( u\right) +c_{B}g\left( u\right) ,\left\vert
u\right\vert ^{m-1}u\right\rangle _{2}-\frac{\widetilde{C}^{q}m^{-q/p}}{%
\left( \varepsilon p\right) ^{q/p}q}\left\langle \left\vert g^{^{\prime
}}\left( u\right) u+mg\left( u\right) \right\vert ^{q},\left\vert
u\right\vert ^{m-1}\right\rangle _{2}  \notag \\
& -\varepsilon m\left( \frac{p}{p+m-1}\right) ^{p}\int_{\Omega }\left\vert
\nabla \left\vert u\right\vert ^{\frac{p+m-1}{p}}\right\vert ^{p}dx.  \notag
\end{align}%
Therefore, combining (\ref{eqn2}) with (\ref{eqn3}), then using (\ref{eqn3b}%
) and (\ref{weaker}), we arrive at the following inequality%
\begin{equation}
\frac{d}{dt}\left\Vert U\right\Vert _{\mathbb{X}^{m+1}}^{m+1}+\gamma \left(
\nu -\varepsilon \right) \int_{\Omega }\left\vert \nabla \left\vert
u\right\vert ^{\frac{p+m-1}{p}}\right\vert ^{p}dx\leq Q_{\tau }\left(
m\right) \left( \left\Vert U\right\Vert _{\mathbb{X}^{m+1}}^{m+1}+1\right) ,
\end{equation}%
for some $\varepsilon \in \left( 0,\nu /q\right) .$ From this point on, the
proof goes on exactly as in the proof of Theorem \ref{linff} (cf. \textbf{%
Steps} \textbf{2} and \textbf{3)}. We omit the details.
\end{proof}

We will now verify the hypothesis in Proposition \ref{linff2} for functions
that satisfy the assumptions of Corollaries \ref{C1} and \ref{C1bis}.

\begin{corollary}
\label{C3}Assume that (H1) holds, and the functions $f,g\in C^{1}\left(
\mathbb{R},\mathbb{R}\right) $ satisfy all the assumptions of Corollary \ref%
{C1}-(i). Then, for any initial data $U_{0}$\ in $\mathbb{X}^{2}$, the
corresponding solution $U\left( t\right) $ of (\ref{1.11})-(\ref{1.12}), (%
\ref{ini}) belongs to $\mathbb{X}^{\infty }$, for each $t>0,$ and estimate (%
\ref{3.3}) holds.
\end{corollary}

\begin{proof}
As in the proof of Corollary \ref{C1}, the leading terms on the left-hand
side of (\ref{weaker}) are, for sufficiently large $\left\vert y\right\vert
\gg 1$ and any $m\geq 1,$%
\begin{equation}
c_{f}\left\vert y\right\vert ^{r_{1}+m-1}+\left( \left\vert \Omega
\right\vert \lambda \right) ^{-1}c_{g}\left\vert y\right\vert ^{r_{2}+m-1}-%
\frac{\widetilde{C}_{\Omega }^{q}m^{-q/p}}{\left( \varepsilon p\right)
^{q/p}q}c_{g}^{q}\left( r_{2}+m-1\right) ^{q}\left\vert y\right\vert
^{q\left( r_{2}-1\right) +m-1},  \label{list2}
\end{equation}%
for some $\varepsilon \in \left( 0,\nu /q\right) $. From Corollary \ref{C1}%
-(i), it holds $\gamma =\max \left( r_{2},q\left( r_{2}-1\right) \right)
<r_{1},$ so the coefficient of the highest order term in (\ref{list2}), for
any $m\geq 1,$ is $c_{f}>0$. Therefore, the desired claim follows
immediately from (\ref{weaker}).
\end{proof}

\begin{corollary}
\label{C2bis}Let $h_{2}=0$ and assume $f,$ $g\in C^{1}\left( \mathbb{R},%
\mathbb{R}\right) $ satisfy%
\begin{equation*}
\lim_{\left\vert y\right\vert \rightarrow \infty }f^{^{\prime }}\left(
s\right) =\left( r_{1}-1\right) c_{f}<0\text{ and }\lim_{\left\vert
y\right\vert \rightarrow \infty }g^{^{\prime }}\left( s\right) =\left(
r_{2}-1\right) c_{g}>0.
\end{equation*}%
Then, the conclusion of Theorem \ref{linff} applies.
\end{corollary}

\begin{proof}
The proof follows, for instance, from Theorem \ref{linff} since the
functions $f,$ $g$ satisfy (\ref{weaker2}).
\end{proof}

Having established that the weak solution is bounded for any positive times,
we also have the following.

\begin{proposition}
\label{higher}Let the assumptions of either Theorem \ref{linff} or
Proposition \ref{linff2} be satisfied. Then, any solution $U\left( t\right) $
of (\ref{1.11})-(\ref{1.12}), (\ref{ini}) belongs to $\mathbb{V}^{p}$, for
each $t>0$, and the following estimate holds:%
\begin{equation}
\sup_{t\geq t_{1}}\left( \left\Vert U\left( t\right) \right\Vert _{\mathbb{V}%
^{p}}^{p}+\int_{t}^{t+1}\left\Vert \partial _{t}U\left( s\right) \right\Vert
_{\mathbb{X}^{2}}^{2}ds\right) \leq C_{2},  \label{wip}
\end{equation}%
for some positive constant $C_{2},$ independent of $t,$ $\epsilon $ and
initial data.
\end{proposition}

\begin{proof}
It suffices to show (\ref{wip}). We first recall that, using assumption (\ref%
{1.1.2}) and the fact that $Tr_{D}:W^{1,p}\left( \Omega \right) \rightarrow
W^{1-1/p,p}\left( \Gamma \right) $ is a bounded map, from (\ref{3.2}) we get%
\begin{equation}
\sup_{t\geq t^{\#}}\int_{t}^{t+1}\left( \left\Vert U\left( s\right)
\right\Vert _{\mathbb{V}^{p}}^{p}+\left\Vert u\left( s\right) \right\Vert
_{L^{r_{1}}\left( \Omega \right) }^{r_{1}}\right) ds\leq C_{0}^{^{\prime }},
\label{wipbis}
\end{equation}%
for some positive constant $C_{0}^{^{\prime }}$ independent of time and
initial data.

Let us now multiply equation (\ref{1.11}) by $\partial _{t}u\left( t\right)
, $ then integrate over $\Omega ,$ and multiply equation (\ref{1.12}) by $%
\partial _{t}v\left( t\right) /b\left( x\right) $ and integrate over $\Gamma
$. Adding the relations that we obtain, we deduce after standard
transformations,%
\begin{align}
& \frac{1}{2}\frac{d}{dt}\left[ \left\langle A(\left\vert \nabla u\left(
t\right) \right\vert ^{2}),1\right\rangle _{2}+2\left\langle F\left( u\left(
t\right) ,1\right) \right\rangle _{2}+2\left\langle G\left( v\left( t\right)
\right) ,\frac{1}{b}\right\rangle _{2,\Gamma }\right.  \label{inter} \\
& \left. -2\left\langle h_{1},u\left( t\right) \right\rangle
_{2}-2\left\langle \frac{h_{2}}{b},v\left( t\right) \right\rangle _{2,\Gamma
}\right]  \notag \\
& =-\left\Vert \partial _{t}u\right\Vert _{L^{2}\left( \Omega \right)
}^{2}-\left\Vert \partial _{t}v\right\Vert _{L^{2}\left( \Gamma ,dS/b\right)
}^{2},  \notag
\end{align}%
for all $t\geq t_{+}$ (with $t_{+}$ as in (\ref{3.3})), where we have set%
\begin{equation*}
A\left( \left\vert y\right\vert ^{2}\right) =\int_{0}^{\left\vert
y\right\vert ^{2}}a\left( s\right) ds,\text{ }F\left( y\right)
=\int_{0}^{y}f\left( s\right) ds\text{, }G\left( y\right)
=\int_{0}^{y}g\left( s\right) ds.
\end{equation*}%
Next, let us define%
\begin{align}
\mathcal{E}\left( t\right) & :=\left\langle A(\left\vert \nabla u\left(
t\right) \right\vert ^{2}),1\right\rangle _{2}+2\left\langle F\left( u\left(
t\right) ,1\right) \right\rangle _{2}+2\left\langle G\left( v\left( t\right)
\right) ,1/b\right\rangle _{2,\Gamma }  \label{ener} \\
& -2\left\langle h_{1},u\left( t\right) \right\rangle _{2}-2\left\langle
h_{2}/b,v\left( t\right) \right\rangle _{2,\Gamma }+C_{F,G}.  \notag
\end{align}%
Here the constant $C_{F,G}>0$ is taken large enough in order to ensure that $%
\mathcal{E}\left( t\right) $ is nonnegative (recall that $F\left( u\right) $
and $G\left( v\right) $ are both bounded by (\ref{3.3})). On the other hand,
on account of (\ref{3.3}), one can easily check, using the fact%
\begin{equation*}
\left\langle A(\left\vert \nabla u\right\vert ^{2}),1\right\rangle _{2}\geq
c_{p}\left\Vert \nabla u\right\Vert _{L^{p}\left( \Omega \right) }^{p},\text{
}c_{p}>0,
\end{equation*}%
that there exists a positive constant $c,$ independent of $t$ and the
initial data, such that
\begin{equation}
\left\Vert \nabla u\left( t\right) \right\Vert _{L^{p}}^{p}-c\leq \mathcal{E}%
\left( t\right) ,  \label{e15}
\end{equation}%
for $\,t\geq \max \left\{ t_{+},t^{\#}\right\} .$ From (\ref{inter}), we
have
\begin{equation}
\frac{d\mathcal{E}\left( t\right) }{dt}+2\left\Vert \partial
_{t}u\right\Vert _{L^{2}\left( \Omega \right) }^{2}+2\left\Vert \partial
_{t}v\right\Vert _{L^{2}\left( \Gamma ,dS/b\right) }^{2}=0,\text{ \ }\forall
\,t\geq \max \left\{ t_{+},t^{\#}\right\} .  \label{4.28}
\end{equation}%
Then, exploiting estimates (\ref{3.3}), (\ref{wipbis}) and (\ref{e15}), we
can apply to (\ref{4.28}) the uniform Gronwall's lemma (see, e.g., \cite{T})
and find a time $t_{1}\geq 1$, depending on $t_{+},$ $t^{\#},$ such that%
\begin{equation}
\left\Vert \nabla u\left( t\right) \right\Vert _{L^{p}\left( \Omega \right)
}^{p}\leq c,\qquad \forall \,t\geq t_{1},  \label{4.29}
\end{equation}%
for some positive constant $c$. Summing up, we conclude by observing that (%
\ref{wip}) follows from (\ref{4.29}) and the boundedness of the trace map $%
Tr_{D}:W^{1,p}\left( \Omega \right) \rightarrow W^{1-1/p,p}\left( \Gamma
\right) .$ The proof is finished.
\end{proof}

Finally, the above dissipative estimates and standard compactness results
(see, e.g., \cite[Section 2]{GW}), allow us to conclude the following.

\begin{theorem}
\label{C2}Let the assumptions of either Theorem \ref{linff} or Proposition %
\ref{linff2}\ be satisfied and, in addition, let (H2) hold. Then, the
dynamical system $\left( \mathcal{S}_{2}\left( t\right) ,\mathbb{X}%
^{2}\right) $ generated by the initial value problem (\ref{1.11})-(\ref{1.12}%
), (\ref{ini})\ possesses the\ global attractor $\mathcal{A}_{gl}\subset
\mathbb{X}^{2}$, which is a bounded subset of $\mathbb{V}^{p}\cap \mathbb{X}%
^{\infty }$. Moreover,%
\begin{equation}
\lim_{t\rightarrow +\infty }dist_{\mathbb{X}^{s_{1},s_{2}}}\left( \mathcal{S}%
_{2}\left( t\right) B,\mathcal{A}_{gl}\right) =0,  \label{4.30}
\end{equation}%
for any finite $s_{1},s_{2}\geq 2$, for all bounded subsets $B$\ of $\mathbb{%
X}^{2}$.
\end{theorem}

\subsection{Trajectory dynamical systems}

In the final part of this section, we shall devote our attention to
constructing the \textquotedblleft usual\textquotedblright\ \emph{weak}
trajectory attractor and verify (using the maximum principle established in
Theorem \ref{linff}) that any solution, belonging to the attractor, is
bounded so that uniqueness holds on the attractor. We will employ a slightly
different construction (compared to e.g., \cite{CV} and references therein)
of the trajectory attractor, which also looks more natural from the physical
point of view. Namely any weak solution $U\left( t\right) $ of (\ref{1.11})-(%
\ref{1.12}), (\ref{ini}) is included in the trajectory phase-space of the
problem if and only if it can be obtained in the limit, as $\epsilon
\rightarrow 0$, of the corresponding solutions $U_{\epsilon }\left( t\right)
$ of the approximate system (\ref{ap})-(\ref{ini2}).

In order to define the trajectory dynamical system for weak solutions
without uniqueness we need to introduce first the appropriate functional
framework. First, let us recall estimates (\ref{dissi}) and (\ref{dissi2})
which hold for any smooth solution $U=U_{\epsilon }$ of the approximate
problem (\ref{ap})-(\ref{ini2}). By a standard application of Gronwall's
inequality (see, e.g., \cite{T}), we get the well-known estimate%
\begin{align}
& \left\Vert U\left( t\right) \right\Vert _{\mathbb{X}^{2}}^{2}+\int_{s}^{t}%
\left( \left\Vert U\left( r\right) \right\Vert _{\mathbb{V}%
^{p}}^{p}+\left\Vert u\left( r\right) \right\Vert _{L^{r_{1}}\left( \Omega
\right) }^{r_{1}}+\widetilde{c}\left\Vert v\left( r\right) \right\Vert
_{L^{r_{2}}\left( \Gamma \right) }^{r_{2}}\right) dr  \label{tdis} \\
& \leq \left\Vert U\left( s\right) \right\Vert _{\mathbb{X}%
^{2}}^{2}e^{-c\left( t-s\right) }+c\left( 1+\left\Vert h_{1}\right\Vert
_{L^{r_{1}^{^{\prime }}}}^{r_{1}^{^{\prime }}}+\left\Vert h_{2}\right\Vert
_{L^{r_{2}^{^{\prime }}}\left( \Gamma \right) }^{r_{2}^{^{\prime }}}\right)
\left( 1-e^{-c\left( t-s\right) }\right) ,  \notag
\end{align}%
for all $t\geq s\geq 0,$ and some appropriate positive constant $c.$ Here
and below $\widetilde{c}=1,$ if (H3b) is assumed and $\widetilde{c}=0$, when
(H3a) holds. Let $\Theta _{+}^{w,loc}$ denote the local weak topology in the
space%
\begin{equation*}
L^{\infty }\left( \mathbb{R}_{+};\mathbb{X}^{2}\right) \cap L^{p}\left(
\mathbb{R}_{+};\mathbb{V}^{p}\right) \cap \left( L^{r_{1}}\left( \mathbb{R}%
_{+};L^{r_{1}}\left( \Omega \right) \right) \times L^{r_{2}}\left( \mathbb{R}%
_{+};L^{r_{2}}\left( \Gamma \right) \right) \right) .
\end{equation*}%
By definition, a sequence $U_{n}\left( t\right) \rightarrow U\left( t\right)
,$ as $n\rightarrow \infty $, in the topology of $\Theta _{+}^{w,loc}$ if,
for every $T>0,$%
\begin{align*}
U_{n}\left( t\right) & \rightarrow U\left( t\right) \text{ }\ast \text{%
-weakly in }L^{\infty }\left( \left[ 0,T\right] ;\mathbb{X}^{2}\right) , \\
U_{n}\left( t\right) & \rightarrow U\left( t\right) \text{ weakly in }%
L^{p}\left( \left[ 0,T\right] ;\mathbb{V}^{p}\right) , \\
U_{n}\left( t\right) & \rightarrow U\left( t\right) \text{ weakly in }%
L^{r_{1}}\left( \left[ 0,T\right] ;L^{r_{1}}\left( \Omega \right) \right)
\times L^{r_{2}}\left( \left[ 0,T\right] ;L^{r_{2}}\left( \Gamma \right)
\right) .
\end{align*}%
We recall that $\Theta _{+}^{w,loc}$ is a Hausdorff and Frechet-Urysohn
space with a countable topology space (see, e.g., \cite{RR}). Next, let $%
\Theta _{+}^{b}$ be the Banach space defined as%
\begin{equation*}
\Theta _{+}^{b}:=L^{\infty }\left( \mathbb{R}_{+};\mathbb{X}^{2}\right) \cap
L^{p}\left( \mathbb{R}_{+};\mathbb{V}^{p}\right) \cap \left( L^{r_{1}}\left(
\mathbb{R}_{+}\times \Omega \right) \times L^{r_{2}}\left( \mathbb{R}%
_{+}\times \Gamma \right) \right) .
\end{equation*}%
Note that the unit ball of $\Theta _{+}^{b}$ is compact in the local weak
topology of $\Theta _{+}^{w,loc}$ (see \cite{RR}).

\begin{definition}
A function $U\in \Theta _{+}^{b}$ is a solution of (\ref{1.11})-(\ref{1.12}%
), (\ref{ini}) with $U_{0}\in \mathbb{X}^{2}$ if it solves (\ref{1.11})-(\ref%
{1.12}) in the sense of Definition \ref{weak}, and there exists a sequence $%
\epsilon _{n}\rightarrow 0,$ a sequence $u_{0\epsilon _{n}}\in C^{\infty
}\left( \overline{\Omega }\right) ,$ $U_{0\epsilon _{n}}=\binom{u_{0\epsilon
_{n}}}{u_{0\epsilon _{n}\mid \Gamma }}$ and a sequence $U_{\epsilon _{n}}$
of classical solutions of the approximate problem (\ref{ap})-(\ref{ini2})
with $\epsilon =\epsilon _{n}$ such that%
\begin{equation}
U_{0}=\mathbb{X}_{w}^{2}-\lim_{n\rightarrow \infty }U_{0\epsilon _{n}}\text{
and }U=\Theta _{+}^{w,loc}-\lim_{n\rightarrow \infty }U_{\epsilon _{n}}.
\label{cucu}
\end{equation}
\end{definition}

Note that (\ref{cucu}) implies in a standard way the weak convergence of%
\begin{equation}
\mathcal{B}_{p,\epsilon _{n}}U_{\epsilon _{n}}\rightarrow \mathcal{B}_{p}U%
\text{ in }L_{loc}^{q}\left( \mathbb{R}_{+};\left( \mathbb{V}^{p}\right)
^{\ast }\right) ,  \label{ann}
\end{equation}%
and consequently, the weak-star convergence of $\partial _{t}U_{\epsilon
_{n}}\rightarrow \partial _{t}U$ in%
\begin{equation}
L_{loc}^{q}\left( \mathbb{R}_{+};\left( \mathbb{V}^{p}\right) ^{\ast
}\right) +(L_{loc}^{r_{1}^{^{\prime }}}\left( \mathbb{R}_{+};L^{r_{^{^{%
\prime }}1}}\left( \Omega \right) \right) \times L_{loc}^{r_{2}^{^{\prime
}}}(\mathbb{R}_{+};L^{r_{2}^{^{\prime }}}\left( \Gamma \right) )).
\label{ann2}
\end{equation}%
This gives the strong convergence $U_{\epsilon _{n}}\left( t\right)
\rightarrow U\left( t\right) $\ in $C_{loc}\left( \mathbb{R}_{+};\left(
\mathbb{V}^{k,p}\right) ^{\ast }\right) $ (see Section 2). Thus, any
solution $U$ of (\ref{1.11})-(\ref{1.12}), (\ref{ini}) is weakly continuous
with values in $\mathbb{X}^{2}$ (see Remark \ref{imp}), and for any $t\geq 0$%
, we have the weak-convergence%
\begin{equation}
U_{\epsilon _{n}}\left( t\right) \rightharpoondown U\left( t\right)
\label{wc}
\end{equation}%
in the space $\mathbb{X}^{2}$. It is important that we do not require the
strong convergence in (\ref{wc}) even for $t=0.$

We can now summarize the results in Section 2 by stating the following.

\begin{proposition}
Let the assumptions of Theorem \ref{T1} hold. Then, for every $U_{0}\in
\mathbb{X}^{2}$, there exists at least one globally defined weak solution $%
U\left( t\right) ,$ $t\in \mathbb{R}_{+},$ of the degenerate problem (\ref%
{1.11})-(\ref{1.12}), with $U\left( 0\right) =U_{0},$ which can be obtained
as a weak limit (\ref{cucu}) of the corresponding solutions $U_{\epsilon
_{n}}\left( t\right) $ of the approximate non-degenerate parabolic system (%
\ref{ap})-(\ref{ini2}).
\end{proposition}

In order to construct the global attractor for the dynamical system
associated with the degenerate parabolic system (\ref{1.11})-(\ref{1.12}),
we need the following definition.

\begin{definition}
Let
\begin{equation*}
M_{U}\left( t\right) :=\inf \left\{ \underset{n\rightarrow \infty }{\lim
\inf \left\Vert U_{\epsilon _{n}}\left( t\right) \right\Vert _{\mathbb{X}%
^{2}}:}U=\Theta _{+}^{w,loc}-\lim_{n\rightarrow \infty }U_{\epsilon _{n}},%
\text{ }U_{0}=\mathbb{X}_{w}^{2}-\lim_{n\rightarrow \infty }U_{\epsilon
_{n}}\left( 0\right) \right\} ,
\end{equation*}%
where the external infimum is taken over all possible sequences of solutions
of the approximate problem (\ref{ap})-(\ref{ini2}), which converges as $%
\epsilon _{n}\rightarrow 0$ to the given solution $U$ of the limit problem (%
\ref{1.11})-(\ref{1.12}), (\ref{ini}).
\end{definition}

Some simple properties of this $M$-functional are stated below.

\begin{proposition}
\label{mu}Let $U$ be a solution of problem (\ref{1.11})-(\ref{1.12}), (\ref%
{ini}) and let $M_{U}\left( t\right) $ be the associated functional. Then,

(a) $\left\Vert U\left( t\right) \right\Vert _{\mathbb{X}^{2}}\leq
M_{U}\left( t\right) ,$ for all $t\in \mathbb{R}_{+}.$

(b) The following estimate holds:
\begin{align}
& \left( M_{U}\left( t\right) \right) ^{2}+\int_{s}^{t}\left( \left\Vert
U\left( r\right) \right\Vert _{\mathbb{V}^{p}}^{p}+\left\Vert u\left(
r\right) \right\Vert _{L^{r_{1}}\left( \Omega \right) }^{r_{1}}+\widetilde{c}%
\left\Vert v\left( r\right) \right\Vert _{L^{r_{2}}\left( \Gamma \right)
}^{r_{2}}\right) dr  \label{cucubis} \\
& \leq \left( M_{U}\left( s\right) \right) ^{2}e^{-c\left( t-s\right)
}+c\left( 1+\left\Vert h_{1}\right\Vert _{L^{r_{1}^{^{\prime
}}}}^{r_{1}^{^{\prime }}}+\left\Vert h_{2}\right\Vert _{L^{r_{2}^{^{\prime
}}}\left( \Gamma \right) }^{r_{2}^{^{\prime }}}\right) \left( 1-e^{-c\left(
t-s\right) }\right) ,  \notag
\end{align}%
where $t\geq s\geq 0$.

(c) $M_{T\left( h\right) U}\left( t\right) \leq M_{U}\left( t+h\right) ,$
for all $h\geq 0$, where $\left( T\left( h\right) U\right) \left( t\right)
:=U\left( t+h\right) .$
\end{proposition}

\begin{proof}
(a) is immediate since the norm $\left\Vert \cdot \right\Vert _{\mathbb{X}%
^{2}}$ is weakly lower semicontinuous and the convergence of $U_{\epsilon
_{n}}$ to $U$ in $\Theta _{+}^{w,loc}$ implies the weak convergence $%
U_{\epsilon _{n}}\left( t\right) \rightarrow U\left( t\right) $ for every $%
t. $ To prove (b), we note that due to the energy estimates for the
approximate parabolic system (\ref{ap})-(\ref{ini2}) (cf. Section 2), we have%
\begin{align}
& \left\Vert U_{\epsilon _{n}}\left( t\right) \right\Vert _{\mathbb{X}%
^{2}}^{2}+\int_{s}^{t}\left( \left\Vert U_{\epsilon _{n}}\left( r\right)
\right\Vert _{\mathbb{V}^{p}}^{p}+\left\Vert u_{\epsilon _{n}}\left(
r\right) \right\Vert _{L^{r_{1}}\left( \Omega \right) }^{r_{1}}+\widetilde{c}%
\left\Vert v_{\epsilon _{n}}\left( r\right) \right\Vert _{L^{r_{2}}\left(
\Gamma \right) }^{r_{2}}\right) dr  \label{cucu2} \\
& \leq \left\Vert U_{\epsilon _{n}}\left( s\right) \right\Vert _{\mathbb{X}%
^{2}}^{2}e^{-c\left( t-s\right) }+c\left( 1+\left\Vert h_{1}\right\Vert
_{L^{r_{1}^{^{\prime }}}}^{r_{1}^{^{\prime }}}+\left\Vert h_{2}\right\Vert
_{L^{r_{2}^{^{\prime }}}\left( \Gamma \right) }^{r_{2}^{^{\prime }}}\right)
\left( 1-e^{-c\left( t-s\right) }\right) ,  \notag
\end{align}%
for every $U_{\epsilon _{n}}$. By definition of $M_{U}$, for every $\delta
>0 $, we can find an approximating sequence $U_{\epsilon _{n}}$ such that%
\begin{equation*}
\underset{n\rightarrow \infty }{\lim \inf }\left\Vert U_{\epsilon
_{n}}\left( t\right) \right\Vert _{\mathbb{X}^{2}}\leq M_{U}\left( t\right)
+\delta .
\end{equation*}%
Passing to the limit, as $n\rightarrow \infty $, in (\ref{cucu2}), we have%
\begin{align*}
& \left( M_{U}\left( t\right) \right) ^{2}+\int_{s}^{t}\left( \left\Vert
U\left( r\right) \right\Vert _{\mathbb{V}^{p}}^{p}+\left\Vert u\left(
r\right) \right\Vert _{L^{r_{1}}\left( \Omega \right) }^{r_{1}}+\widetilde{c}%
\left\Vert v\left( r\right) \right\Vert _{L^{r_{2}}\left( \Gamma \right)
}^{r_{2}}\right) dr \\
& \leq \underset{n\rightarrow \infty }{\lim \inf }\left( \left\Vert
U_{\epsilon _{n}}\left( t\right) \right\Vert _{\mathbb{X}^{2}}^{2}+%
\int_{s}^{t}\left( \left\Vert U_{\epsilon _{n}}\left( r\right) \right\Vert _{%
\mathbb{V}^{p}}^{p}+\left\Vert u_{\epsilon _{n}}\left( r\right) \right\Vert
_{L^{r_{1}}\left( \Omega \right) }^{r_{1}}+\widetilde{c}\left\Vert
v_{\epsilon _{n}}\left( r\right) \right\Vert _{L^{r_{2}}\left( \Gamma
\right) }^{r_{2}}\right) dr\right) \\
& \leq \left( M_{U}\left( s\right) +\delta \right) ^{2}e^{-c\left(
t-s\right) }+c\left( 1+\left\Vert h_{1}\right\Vert _{L^{r_{1}^{^{\prime
}}}}^{r_{1}^{^{\prime }}}+\left\Vert h_{2}\right\Vert _{L^{r_{2}^{^{\prime
}}}\left( \Gamma \right) }^{r_{2}^{^{\prime }}}\right) \left( 1-e^{-c\left(
t-s\right) }\right) ,
\end{align*}%
and since $\delta \rightarrow 0$ is arbitrary, we get the desired
inequality. The third assertion is also immediate since the infimum in the
definition of $M_{T\left( h\right) U}\left( t\right) $ is taken over the
larger set of admissible approximating sequences than the infimum in the
definition of $M_{U\left( t+h\right) }$.
\end{proof}

We are now ready to construct the trajectory-phase space, the trajectory
semigroup and the kernel associated with the degenerate problem (\ref{1.11}%
)-(\ref{1.12}). To this end, let $\mathcal{K}_{+}\subset \Theta _{+}^{b}$ be
the set of all solutions of (\ref{1.11})-(\ref{1.12}), in the sense of
Definition \ref{weak}, which corresponds to all $U_{0}\in \mathbb{X}^{2}$,
and let
\begin{equation}
T\left( h\right) :\mathcal{K}_{+}\rightarrow \mathcal{K}_{+},\text{ }h\geq 0,
\label{cucu3}
\end{equation}%
be%
\begin{equation}
\left( T\left( h\right) U\right) \left( t\right) :=U\left( t+h\right) .
\label{cucu4}
\end{equation}%
We shall refer to $\mathcal{K}_{+}$ and $T\left( h\right) :\mathcal{K}%
_{+}\rightarrow \mathcal{K}_{+}$ as the trajectory phase space and the
trajectory dynamical system, respectively, associated with the degenerate
parabolic system (\ref{1.11})-(\ref{1.12}). In addition, we endow the set $%
\mathcal{K}_{+}$ with the topology induced by the embedding $\mathcal{K}%
_{+}\subset \Theta _{+}^{w,loc}$ and we will say that a set $B\subset
\mathcal{K}_{+}$ is $M$-bounded if
\begin{equation*}
M_{B}\left( 0\right) :=\sup_{U\in B}M_{U}\left( 0\right) <\infty .
\end{equation*}%
Note that any $M$-bounded set $B\subset \mathcal{K}_{+}$ is bounded in the
norm of $\Theta _{+}^{b}$. Finally, a kernel $\mathcal{K\subset }\Theta
^{b}, $%
\begin{equation*}
\Theta ^{b}:=L^{\infty }\left( \mathbb{R};\mathbb{X}^{2}\right) \cap
L^{p}\left( \mathbb{R};\mathbb{V}^{p}\right) \cap \left( L^{r_{1}}\left(
\mathbb{R}\times \Omega \right) \times L^{r_{2}}\left( \mathbb{R}\times
\Gamma \right) \right)
\end{equation*}%
consists of all complete (defined for all $t\in \mathbb{R}$) bounded
solutions of (\ref{1.11})-(\ref{1.12}), (\ref{ini}) which can be obtained as
the weak limit, as $\epsilon _{n}\rightarrow 0$, of the appropriate
solutions of the approximate non-degenerate parabolic system (\ref{ap})-(\ref%
{ini2}). Namely, $U\in \mathcal{K}$ if and only if there exists a sequence $%
\epsilon _{n}\rightarrow 0,$ a sequence of times $t_{n}\rightarrow -\infty ,$
and a bounded sequence of initial data $u_{0\epsilon _{n}}\in C^{\infty
}\left( \overline{\Omega }\right) ,$ $\left\Vert U_{0\epsilon
_{n}}\right\Vert _{\mathbb{X}^{2}}\leq C$, such that the corresponding
solutions $U_{\epsilon _{n}}$ of (\ref{ap}) on the interval $[t_{n},+\infty
) $ with initial data $U_{\epsilon _{n}}\left( t_{n}\right) =U_{0\epsilon
_{n}} $ converges weakly in $\Theta ^{b}$ to the complete solution $U$
considered.

We now recall the definition of the global attractor for the trajectory
dynamical system $\left( T\left( h\right) ,\mathcal{K}_{+}\right) $ (see
\cite{CV} for more details; cf. also \cite{MZ, Z})$.$

\begin{definition}
\label{deftraj}A set $\mathcal{A}_{tr}\subset \mathcal{K}_{+}$ is a (weak)
trajectory attractor associated with the degenerate parabolic system (\ref%
{1.11})-(\ref{1.12}) (= global attractor for the trajectory dynamical system
$\left( T\left( h\right) ,\mathcal{K}_{+}\right) $) if the following
conditions are satisfied:

(i) $\mathcal{A}_{tr}$ is compact in $\mathcal{K}_{+}$ and is $M$-bounded;

(ii) It is strictly invariant, i.e., $T\left( h\right) \mathcal{A}_{tr}=%
\mathcal{A}_{tr}$, $h>0$;

(iii) It attracts the images of bounded ($M$-bounded) sets as $h\rightarrow
\infty $, i.e., for every $B$ bounded in $\mathcal{K}_{+}$ and every
neighborhood $O(\mathcal{A}_{tr})$ of $\mathcal{A}_{tr}$ (in the topology of
$\Theta _{+}^{w,loc}$), there exists $h_{0}=h_{0}\left( B,O\right) $ such
that $T(h)B\subset O\left( \mathcal{A}_{tr}\right) ,$ $\forall h\geq h_{0}$.
\end{definition}

The next theorem can be considered as the second main result of this section.

\begin{theorem}
\label{traj}Let $p\in (\frac{2N}{N+2},\infty )\cap \left( 1,\infty \right) $%
. Let all the assumptions of Theorem \ref{T1} be satisfied. Then, the
degenerate parabolic problem (\ref{1.11})-(\ref{1.12}), (\ref{ini})
possesses a trajectory attractor $\mathcal{A}_{tr}\subset \Theta _{+}^{b}$
and the following description holds:%
\begin{equation}
\mathcal{A}_{tr}=\Pi _{t\geq 0}\left( \mathcal{K}\right) .  \label{com}
\end{equation}%
Here and below, $\Pi _{t\in I}\left( f\right) $ denotes the restriction on $%
I $ of a function $f$ defined on $\mathbb{R}$.
\end{theorem}

\begin{proof}
According to general theory \cite{CV, MZ, Z}, we are only required to check
that the trajectory dynamical system is continuous and that it possesses a
compact and $M$-bounded absorbing set. The continuity is immediate since $%
T\left( h\right) $ are continuous on $\Theta _{+}^{w,loc}$. The estimate (b)
in Proposition \ref{mu} guarantees that the set%
\begin{equation}
\mathcal{B}:=\left\{ U\in \mathcal{K}_{+}:\left( M_{U}\left( 0\right)
\right) ^{2}\leq 2c\left( 1+\left\Vert h_{1}\right\Vert _{L^{r_{1}^{^{\prime
}}}}^{r_{1}^{^{\prime }}}+\left\Vert h_{2}\right\Vert _{L^{r_{2}^{^{\prime
}}}\left( \Gamma \right) }^{r_{2}^{^{\prime }}}\right) \right\}  \label{abs}
\end{equation}%
will be absorbing for the semigroup $T\left( h\right) :\mathcal{K}%
_{+}\rightarrow \mathcal{K}_{+}$. Moreover, this set is semi-invariant. This
follows from Proposition \ref{mu}, (c) since%
\begin{align*}
\left( M_{T\left( h\right) U}\left( 0\right) \right) ^{2}& \leq \left(
M_{U}\left( h\right) \right) ^{2} \\
& \leq \left( M_{U}\left( 0\right) \right) ^{2}e^{-c\left( t-s\right)
}+c\left( 1+\left\Vert h_{1}\right\Vert _{L^{r_{1}^{^{\prime
}}}}^{r_{1}^{^{\prime }}}+\left\Vert h_{2}\right\Vert _{L^{r_{2}^{^{\prime
}}}\left( \Gamma \right) }^{r_{2}^{^{\prime }}}\right) \left( 1-e^{-c\left(
t-s\right) }\right) \\
& \leq 2c\left( 1+\left\Vert h_{1}\right\Vert _{L^{r_{1}^{^{\prime
}}}}^{r_{1}^{^{\prime }}}+\left\Vert h_{2}\right\Vert _{L^{r_{2}^{^{\prime
}}}\left( \Gamma \right) }^{r_{2}^{^{\prime }}}\right) ,
\end{align*}%
for all $U\in \mathcal{B}$. Therefore, $T\left( h\right) \mathcal{B}\subset
\mathcal{B}$. It remains to show that $\mathcal{B}$ is compact. Due to the
inequality (b) in Proposition \ref{mu}, the set $\mathcal{B}$ is bounded in $%
\Theta _{+}^{b}$, and therefore precompact in $\Theta _{+}^{w,loc}$. Thus,
we only need to show that $\mathcal{B}$ is sequentially closed, i.e., if $%
U_{n}\in \mathcal{B}$ and $U=\Theta _{+}^{w,loc}-\lim_{n\rightarrow \infty
}U_{n},$ then $U\in \mathcal{B}$ as well.

For $U_{n}\in \mathcal{B}$, $M_{U_{n}}\left( 0\right) $ is bounded. By
estimate (b) in Proposition \ref{mu} the sequence $M_{U_{n}}\left( t\right) $
is bounded for all $t\geq 0$. Moreover, since every $U_{n}$ is a solution of
the degenerate problem (\ref{1.11})-(\ref{1.12}), there exists a sequence $%
U_{n,\epsilon _{l}}$ of solutions to the approximate non-degenerate problem (%
\ref{ap})-(\ref{ini2}) such that $\epsilon _{l}=\epsilon _{n,l}\rightarrow 0$
as $k\rightarrow \infty $, and%
\begin{equation*}
U_{n}=\Theta _{+}^{w,loc}-\lim_{k\rightarrow \infty }U_{n,\epsilon _{l}}.
\end{equation*}%
Without loss of generality, we may assume that%
\begin{equation}
\left\vert M_{U_{n}}\left( 0\right) -\left\Vert U_{n,\epsilon _{l}}\left(
0\right) \right\Vert _{\mathbb{X}^{2}}\right\vert \leq \frac{1}{n},\text{ }%
e_{l}\leq \frac{1}{n},  \label{star}
\end{equation}%
for all $l\in \mathbb{N}$, and we may also suppose that%
\begin{equation*}
\underset{n\rightarrow \infty }{\lim \inf }M_{U_{n}}\left( 0\right)
=\lim_{n\rightarrow \infty }M_{U_{n}}\left( 0\right)
\end{equation*}%
(we may pass to a subsequence in $n$ if necessary). It remains to show that
we can extract from $\left\{ U_{n,\epsilon _{l}}\right\} _{n,l\in \mathbb{N}%
} $ a one parametric sequence which will converge to the limit function $U$.
To this end, recall that the topology of $\Theta _{+}^{w,loc}$ is metrizable
on every bounded set of $\Theta _{+}^{b}$ (see \cite{RR}). Let $\delta >0$
be such that all $U_{n}$ belong to the closed ball $B_{\delta }$ of $\Theta
_{+}^{b}$. Evidently, $U\in B_{\delta }\subset B_{2\delta }$, and we may
also assume that $U_{n,\epsilon _{l}}\in B_{2\delta },$ for all $n,l\in
\mathbb{N}$. Indeed, the sequence $U_{n,\epsilon _{l}}\left( 0\right) $ is
uniformly bounded in $n,l$, due to (\ref{star}), and since $M_{U_{n}}\left(
0\right) $ is bounded, then recalling estimate (\ref{cucu2}), we also get
that $U_{n,\epsilon _{l}}\left( t\right) $ is uniformly bounded with respect
to $n,l$ and $t$. Let $d\left( \cdot ,\cdot \right) $ be a metric on $%
B_{2\delta }$. Thus,%
\begin{equation*}
\lim_{n\rightarrow \infty }d\left( U_{n},U\right) =0,\text{ }%
\lim_{k\rightarrow \infty }d\left( U_{n},U_{n,\epsilon _{l}}\right) =0,
\end{equation*}%
for every $n$. Therefore, for any $n$, there exists $l_{0}=l_{0}\left(
n\right) $ such that $d\left( U,U_{n,\epsilon _{l}}\right) \leq 1/n$, for
all $l\geq l_{0}$. Thus, we have $d(U,U_{n,\epsilon _{l_{0}}})\rightarrow 0$
as $n\rightarrow \infty ,$ and therefore%
\begin{equation*}
U=\Theta _{+}^{w,loc}-\lim_{n\rightarrow \infty }U_{n,\epsilon _{l_{0}}}.
\end{equation*}%
Moreover, thanks to (\ref{star}), $\epsilon _{l_{0}}\rightarrow 0$ as $%
n\rightarrow \infty $ and so $U$ is a solution of the degenerate problem, and%
\begin{equation*}
M_{U}\left( 0\right) \leq \underset{n\rightarrow \infty }{\lim \inf }%
\left\Vert U_{n,\epsilon _{l}}\left( 0\right) \right\Vert _{\mathbb{X}^{2}}=%
\underset{n\rightarrow \infty }{\lim \inf }M_{U_{n}}\left( 0\right) .
\end{equation*}%
Thus, $\mathcal{B}$ is indeed a compact semi-invariant absorbing set for $%
\left( T\left( h\right) ,\mathcal{K}_{+}\right) ,$ and the desired attractor
can now be found in a standard way, as the $\omega $-limit set of $\mathcal{B%
}:$%
\begin{equation*}
\mathcal{A}_{tr}=\omega \left( \mathcal{B}\right) =\bigcap\nolimits_{h\geq
0}T\left( h\right) \mathcal{B}\text{.}
\end{equation*}%
The description (\ref{com}) is a standard corollary of this explicit formula
and the diagonalization procedure described above. This completes the proof.
\end{proof}

What is the connection between the dynamical system $\left( \mathcal{S}%
_{2}\left( t\right) ,\mathbb{X}^{2}\right) $ introduced at the end of
Section 2, and the trajectory dynamical system $\left( T\left( h\right) ,%
\mathcal{K}_{+}\right) $ constructed here? It turns out that, under the
assumptions of Proposition \ref{uniq}, the solution $U(t)$ of the system (%
\ref{1.11})-(\ref{1.12}) is unique and, consequently, this parabolic system
generates a semigroup in the classical phase space $\mathbb{X}^{2},$ in a
standard way by the formula (\ref{2.15}). If we define the map%
\begin{equation*}
\Pi _{0}:\mathcal{K}_{+}\rightarrow \mathbb{X}^{2}\text{, }\Pi _{0}\left(
U\left( t\right) \right) :=U\left( 0\right) ,
\end{equation*}%
we see that the map $\Pi _{0}$ is one-to-one and, in fact, $\Pi _{0}$
defines a Lipschitz homeomorphism between $\mathcal{K}_{+}$ and $\mathbb{X}%
_{w}^{2}$ (i.e., $\mathbb{X}^{2}$ endowed with the weak topology).
Therefore, when uniqueness holds (for instance, if we require that the
functions $f,g$ satisfy (H2)),%
\begin{equation*}
\mathcal{S}_{2}\left( t\right) =\Pi _{0}T\left( t\right) \Pi _{0}^{-1},
\end{equation*}%
the trajectory dynamical system $\left( T\left( t\right) ,\mathcal{K}%
_{+}\right) $\ is conjugated to the classical dynamical system $\mathcal{S}%
_{2}\left( t\right) $ defined on the phase space $\mathbb{X}^{2}$ endowed
with the weak topology. However, we note that assumption (H2) is quite
restrictive as it does not allow for a competing scenario between the
nonlinearities $f,g$, as proposed in e.g., Corollaries \ref{C1} and \ref%
{C1bis}. So without assumption (H2) the uniqueness problem for (\ref{1.11})-(%
\ref{1.12}) is not known, and the classical semigroup $\mathcal{S}_{2}\left(
t\right) $ can be defined as a semigroup of multi-valued maps only (see,
e.g., \cite{MZ} for further details). The trajectory dynamical approach
allows us to avoid the use of multivalued maps and to apply the usual theory
of global attractors to investigate the long term behavior of the degenerate
parabolic system (\ref{1.11})-(\ref{1.12}), (\ref{ini}).

As a consequence of Theorem \ref{traj}, we can also state the following.

\begin{corollary}
Let the assumptions of Theorem \ref{traj} hold and let $B\subset \mathcal{K}%
_{+}$ be an arbitrary $M$-bounded set. Then, for every $T\in \mathbb{R}_{+}$
and every $s\in (0,k)$, the following convergence holds:%
\begin{equation}
\lim_{h\rightarrow \infty }dist_{C\left( \left( h,T+h\right) ;\left( \mathbb{%
V}^{s,p}\right) ^{\ast }\right) }\left( B_{\mid \left( h,T+h\right) },%
\mathcal{A}_{tr\mid \left( h,T+h\right) }\right) =0.  \label{convtraj}
\end{equation}
\end{corollary}

\begin{proof}
Indeed, from the fact that every $U\in \mathcal{K}_{+}$ is a weak solution
in the sense of Definition \ref{weak} (cf. (\ref{VF})), we can express and
estimate (as in Section 2) the time derivative of $U$, i.e., $\partial
_{t}U\in L_{loc}^{s}\left( \mathbb{R}_{+};\left( \mathbb{V}^{k,p}\right)
^{\ast }\right) .$ Next, since the embedding%
\begin{equation}
\left\{ U:U\in L_{loc}^{\infty }\left( \mathbb{R}_{+};\mathbb{X}^{2}\right) ,%
\text{ }\partial _{t}U\in L_{loc}^{s}\left( \mathbb{R}_{+};\left( \mathbb{V}%
^{k,p}\right) ^{\ast }\right) \right\} \subset C_{loc}\left( \mathbb{R}%
_{+};\left( \mathbb{V}^{s,p}\right) ^{\ast }\right) ,  \label{emb}
\end{equation}%
is compact for every $0<s<k$, then $\mathcal{K}_{+}\subset C_{loc}\left(
\mathbb{R}_{+};\left( \mathbb{V}^{s,p}\right) ^{\ast }\right) $ is also
compact in the sense that every $M$-bounded subset of $\mathcal{K}_{+}$ is a
precompact set in the set $C_{loc}\left( \mathbb{R}_{+};\left( \mathbb{V}%
^{s,p}\right) ^{\ast }\right) .$ Thus, the above convergence (\ref{convtraj}%
) is an immediate corollary of Definition \ref{deftraj}, (iii).
\end{proof}

In the sequel, we shall also verify, under additional assumptions on the
nonlinearities (which still allow for a competing behavior between $f$ and $%
g $), that every solution $U$ of the degenerate parabolic system is
uniformly bounded on the (weak) trajectory attractor $\mathcal{A}_{tr}$.
Then, using this fact we can establish that the solution of (\ref{1.11})-(%
\ref{1.12}), (\ref{ini}) is unique on the attractor. From now, we will
always assume that $p\geq 2.$

\begin{theorem}
\label{lc}Let the assumptions of either Theorem \ref{linff} or Proposition %
\ref{linff2} be satisfied. Then, for every complete solution $U\in \mathcal{K%
}$, we have $U\in L^{\infty }\left( \mathbb{R};\mathbb{X}^{\infty }\cap
\mathbb{V}^{p}\right) $ and the following estimate holds:%
\begin{equation}
\left\Vert U\left( t\right) \right\Vert _{\mathbb{X}^{\infty }\cap \mathbb{V}%
^{p}}\leq Q(1+\left\Vert h_{1}\right\Vert _{L^{\infty }\left( \Omega \right)
}+\left\Vert h_{2}\right\Vert _{L^{\infty }\left( \Gamma \right) }),
\label{linfcomp}
\end{equation}%
for all $t\in \mathbb{R}$, for some monotone non-decreasing function $Q$\
independent of $U$, $t$ and initial data. Moreover, $U\in C_{loc}\left(
\mathbb{R};\mathbb{X}^{2}\right) ,$ for every $U\in \mathcal{K}$.
\end{theorem}

\begin{proof}
The proof is essentially based on the maximum principle deduced in Section
2, and the description of $\mathcal{K}$ from Theorem \ref{traj}. Let $U\in
\mathcal{K}$ be an arbitrary complete solution, i.e., let $\epsilon
_{n}\rightarrow 0$, $t_{n}\rightarrow -\infty $ and $U_{\epsilon _{n}}\left(
t\right) ,$ $t\geq t_{n}$ be the sequence of solutions of the approximate
parabolic system (\ref{ap}) with $U_{\epsilon _{n}}\left( t_{n}\right)
=U_{0\epsilon _{n}},$ where $\left\Vert U_{0\epsilon _{n}}\right\Vert _{%
\mathbb{X}^{2}}\leq C$, uniformly with respect to $n$. Let us now fix an
arbitrary $T\in \mathbb{R}$. Then from the convergence $U_{\epsilon
_{n}}\rightarrow U$, we also know that%
\begin{equation}
U_{\epsilon _{n}}\left( T\right) \rightarrow U\left( T\right)  \label{pass}
\end{equation}%
strongly in $\mathbb{X}^{2}$ (passing to a subsequence in $n$, if
necessary), since the embedding $\mathcal{K}_{+}\subset L_{loc}^{p}\left(
\mathbb{R}_{+};\mathbb{X}^{2}\right) $ is compact (i.e., any $M$-bounded
subset of $\mathcal{K}_{+}$ is precompact in $L_{loc}^{p}\left( \mathbb{R}%
_{+};\mathbb{X}^{2}\right) $). It follows from Theorem \ref{traj} and
estimate (\ref{cucu2}) that%
\begin{align}
& \left\Vert U_{\epsilon _{n}}\left( T\right) \right\Vert _{\mathbb{X}%
^{2}}^{2}+\int_{T}^{T+1}\left( \left\Vert U_{\epsilon _{n}}\left( r\right)
\right\Vert _{\mathbb{V}^{p}}^{p}+\left\Vert u_{\epsilon _{n}}\left(
r\right) \right\Vert _{L^{r_{1}}\left( \Omega \right) }^{r_{1}}+\widetilde{c}%
\left\Vert v_{\epsilon _{n}}\left( r\right) \right\Vert _{L^{r_{2}}\left(
\Gamma \right) }^{r_{2}}\right) dr  \label{cucu5} \\
& \leq 2c\left( 1+\left\Vert h_{1}\right\Vert _{L^{r_{1}^{^{\prime
}}}}^{r_{1}^{^{\prime }}}+\left\Vert h_{2}\right\Vert _{L^{r_{2}^{^{\prime
}}}\left( \Gamma \right) }^{r_{2}^{^{\prime }}}\right) ,  \notag
\end{align}%
where the constant $c>0$ is independent of $\epsilon _{n}$ and $T.$ By
modifying the proof of Theorem \ref{linff} in a suitable way, we arrive at
the following inequality for the approximate solutions $U_{\epsilon _{n}}$,%
\begin{align}
\sup_{t\geq t_{n}}\left\Vert U_{\epsilon _{n}}\left( t\right) \right\Vert _{%
\mathbb{X}^{\infty }}& \leq Q(c+\left\Vert U\right\Vert _{L^{\infty }\left(
\mathbb{R};\mathbb{X}^{2}\right) }),  \label{e12q} \\
& \leq Q(1+\left\Vert h_{1}\right\Vert _{L^{\infty }\left( \Omega \right)
}+\left\Vert h_{2}\right\Vert _{L^{\infty }\left( \Gamma \right) }),  \notag
\end{align}%
for some monotone function $Q$ independent of $\epsilon _{n},$ $T,$ $t.$
Passing to the limit, in a standard way in (\ref{e12q}), we may think that $%
U_{\epsilon _{n}}\left( T\right) \rightharpoondown U_{0}\left( T\right) $
weakly-star in $\mathbb{X}^{\infty },$ for some $U_{0}\in \mathbb{X}^{\infty
}$ such that $U_{0}$ satisfies%
\begin{equation}
\left\Vert U_{0}\left( T\right) \right\Vert _{\mathbb{X}^{\infty }}\leq
Q(1+\left\Vert h_{1}\right\Vert _{L^{\infty }\left( \Omega \right)
}+\left\Vert h_{2}\right\Vert _{L^{\infty }\left( \Gamma \right) }).
\label{e12tt}
\end{equation}%
This together with (\ref{pass}) gives $U\left( T\right) \in \mathbb{X}%
^{\infty }$ such that $U\left( T\right) $ satisfies the analogue of (\ref%
{e12tt}).

Finally, it remains to prove that $U\left( T\right) $ is also bounded in $%
\mathbb{V}^{p}$ by essentially arguing as in the proof of Proposition \ref%
{higher}. Note once again that in contrast to the limit case $\epsilon
_{n}=0,$ the parabolic system (\ref{ap}) is non-degenerate if $\epsilon
_{n}>0$, and we have enough regularity of $U_{\epsilon _{n}}$ to justify the
multiplication by the test functions $\partial _{t}U_{\epsilon _{n}}$ in the
weak formulation (\ref{VF}) (cf. Definition \ref{weak}). Analogous to (\ref%
{inter})-(\ref{4.28}), we get%
\begin{equation}
\partial _{t}\mathcal{E}_{\epsilon _{n}}\left( t\right) +2\left\Vert
\partial _{t}u_{\epsilon _{n}}\left( t\right) \right\Vert _{L^{2}\left(
\Omega \right) }^{2}+2\left\Vert \partial _{t}v_{\epsilon _{n}}\left(
t\right) \right\Vert _{L^{2}\left( \Gamma ,dS/b\right) }^{2}=0,  \label{e13b}
\end{equation}%
for the energy $\mathcal{E}_{\epsilon _{n}}$ defined in (\ref{ener}).
Therefore, Gronwall's inequality applied to (\ref{e13b}) yields as in (\ref%
{4.28})-(\ref{4.29}),%
\begin{align}
& \mathcal{E}_{\epsilon _{n}}\left( T\right) +\int_{T}^{T+1}\left(
\left\Vert \partial _{t}u_{\epsilon _{n}}\left( r\right) \right\Vert
_{L^{2}\left( \Omega \right) }^{2}+\left\Vert \partial _{t}v_{\epsilon
_{n}}\left( r\right) \right\Vert _{L^{2}\left( \Gamma ,dS/b\right)
}^{2}\right) dr  \label{e14b} \\
& \leq Q(1+\left\Vert h_{1}\right\Vert _{L^{\infty }\left( \Omega \right)
}+\left\Vert h_{2}\right\Vert _{L^{\infty }\left( \Gamma \right) }),  \notag
\end{align}%
which in light of (\ref{e15}) gives that $U_{\epsilon _{n}}\left( T\right) $
is uniformly (in $T,\epsilon _{n}$) bounded in $\mathbb{V}^{p}$. Passing now
to the limit as $\epsilon _{n}\rightarrow 0$, we obtain in a standard way
that $U\left( T\right) $ satisfies (\ref{e14b}) as well. The desired
inequality (\ref{linfcomp}) follows immediately from (\ref{e12tt}) and (\ref%
{e14b}). The last assertion in the theorem is a standard corollary of the
energy identity (\ref{coo}). So, the proof is complete.
\end{proof}

The (forward) uniqueness theorem holds for bounded solutions of (\ref{1.11}%
)-(\ref{1.12}), (\ref{ini}).

\begin{theorem}
\label{aa}Let the assumptions of Theorem \ref{lc} be satisfied, and let $%
f,g\in C^{1}\left( \mathbb{R},\mathbb{R}\right) $. Consider two functions $%
U_{1}\left( t\right) ,\,U_{2}\left( t\right) \in C\left( \left[ 0,T\right] ;%
\mathbb{X}^{2}\right) $ which solve (\ref{1.11})-(\ref{1.12}), (\ref{ini})
in the sense of Definition \ref{weak}. In addition, let $U_{1},U_{2}\in
L^{\infty }\left( \left[ 0,T\right] ;\mathbb{X}^{\infty }\right) $. Then, $%
U_{1}\left( 0\right) =U_{2}\left( 0\right) $ implies that $U_{1}\left(
t\right) =U_{2}\left( t\right) ,$ for all $t\in \left[ 0,T\right] .$

\begin{proof}
The proof is an immediate consequence of inequality (\ref{uniqqq}) (cf.
Proposition \ref{uniq}). Indeed, for $E\left( t\right) :=\left\Vert
U_{1}\left( t\right) -U_{2}\left( t\right) \right\Vert _{\mathbb{X}%
^{2}}^{2}, $ we have $E\in C^{1}\left( 0,T\right) $, and%
\begin{equation}
\frac{1}{2}\partial _{t}E\left( t\right) \leq Q\left( \left\Vert U_{1}\left(
t\right) \right\Vert _{\mathbb{X}^{\infty }},\left\Vert U_{2}\left( t\right)
\right\Vert _{\mathbb{X}^{\infty }}\right) E\left( t\right) \leq cE\left(
t\right) ,  \label{uniqq}
\end{equation}%
which yields the desired claim on account of the application of Gronwall's
inequality.
\end{proof}
\end{theorem}

\begin{remark}
\label{buniq}We note that in contrast to the non-degenerate case of $a\left(
s\right) \equiv \nu $, we do \emph{not} know whether the backwards
uniqueness theorem holds for bounded\ solutions of the parabolic system (\ref%
{1.11})-(\ref{1.12}) for $p\neq 2$. Namely, if the equality $U_{1}\left(
T\right) =U_{2}\left( T\right) $ holds for some $T\geq 0$, then we have $%
U_{1}\left( t\right) =U_{2}\left( t\right) ,$ for $t\leq T$, as well.
Indeed, in the former case the parabolic system (\ref{1.11})-(\ref{1.12}) is
just a reaction-diffusion equation with dynamic boundary conditions, so we
can establish additional regularity of the weak solutions in $L^{\infty
}\left( \mathbb{R}_{+};\mathbb{V}^{2,2}\right) \cap W^{1,\infty }\left(
\mathbb{R}_{+};\mathbb{V}^{1,2}\right) $, following \cite[Theorem 2.3]{GalNS}
(see also below). Thus, exploiting a well-known theorem (see, e.g., \cite[%
Theorem 11.10]{R}; cf. also \cite[Chapter III]{T}), we can easily establish
the backwards uniqueness result in this case.
\end{remark}

Finally, we observe that since a bounded weak solution $U\left( t\right) \in
L^{\infty }\left( \mathbb{R}_{+};\mathbb{X}^{\infty }\cap \mathbb{V}%
^{p}\right) $ is unique, we may define a global attractor $\mathcal{A}_{gl}$
for the parabolic system (\ref{1.11})-(\ref{1.12}), (\ref{ini}) by the
standard expression%
\begin{equation}
\mathcal{A}_{gl}:=\Pi _{0}\mathcal{A}_{tr},  \label{proj}
\end{equation}%
and define a classical semigroup on this attractor via%
\begin{equation}
\mathcal{S}_{t}:\mathcal{A}_{gl}\rightarrow \mathcal{A}_{gl},\text{ }%
\mathcal{S}_{t}U\left( 0\right) =U\left( t\right) ,  \label{proj2}
\end{equation}%
Here, $U\left( t\right) $ is the unique (bounded) weak solution of (\ref%
{1.11})-(\ref{1.12}), (\ref{ini}), such that $U\left( t\right) $ satisfies
the energy identity (\ref{coo}). We also note that estimate (\ref{linfcomp})
gives a uniform estimate of the $L^{\infty }$-norm of the trajectories
belonging to the attractor $\mathcal{A}_{gl}$. Therefore, the growth rate of
the nonlinearities $f$, $g$ with respect to $U$ becomes nonessential for
further investigations of global attractors and we can study them exactly as
in, e.g., \cite{R,T}.

\subsection{Strong trajectory attractors: the semilinear case}

Let $\left( T_{2}\left( t\right) ,\mathcal{K}_{+}\right) $ be the trajectory
dynamical system associated with the reaction-diffusion equation (\ref%
{1.11bb}), subject to the dynamic boundary condition (\ref{1.12bb}) (see
Section 1). In this section, we shall establish additional regularity
estimates for the weak solutions of (\ref{1.11bb})-(\ref{1.12bb}), and
obtain as a by-product, that the weak trajectory attractor $\mathcal{A}%
_{tr}, $ constructed in Theorem \ref{traj}, is in fact a \emph{strong}
trajectory attractor. In order to do so we will verify, for every $U\in
\mathcal{A}_{tr} $, that the attraction property and the compactness holds
not only in the weak topology of $\Theta _{+}^{w,loc}$, but also in the
strong topology of $\Theta _{+}^{s,loc}$. The definition of a strong
trajectory attractor is obtained by replacing the weak attraction condition
(iii)\ in Definition \ref{deftraj} by the condition of strong attraction in
the topology of $\Theta _{+}^{s,loc}$.

We have the following proposition, whose proof goes essentially as in \cite[%
Theorem 2.3]{GalNS}.

\begin{proposition}
\label{propsem}Let the assumptions of either Theorem \ref{linff} or
Proposition \ref{linff2} be satisfied, and $f,g\in C^{1}\left( \mathbb{R},%
\mathbb{R}\right) .$ Every bounded complete weak solution $U\in \mathcal{K}$
of problem (\ref{1.11bb})-(\ref{1.12bb}) belongs to $L^{\infty }\left(
\mathbb{R};\mathbb{V}^{2,2}\right) \cap W^{1,\infty }\left( \mathbb{R};%
\mathbb{V}^{1,2}\right) ,$ and the following inequality holds for $t\in
\mathbb{R}$,%
\begin{align}
& \left\Vert U\left( t\right) \right\Vert _{\mathbb{V}^{2,2}}^{2}+\left\Vert
\partial _{t}U\left( t\right) \right\Vert _{\mathbb{V}^{1,2}}^{2}+%
\int_{t}^{t+1}\left\Vert \partial _{t}^{2}U\left( r\right) \right\Vert _{%
\mathbb{X}^{2}}^{2}dr  \label{rege} \\
& \leq Q(1+\left\Vert h_{1}\right\Vert _{L^{\infty }\left( \Omega \right)
}+\left\Vert h_{2}\right\Vert _{L^{\infty }\left( \Gamma \right) }),  \notag
\end{align}%
where $Q$ is a monotone nondecreasing function, independent of $t,U$ and the
data.
\end{proposition}

\begin{proof}
Indeed, having established the $L^{\infty }$-estimate (\ref{linfcomp}), (\ref%
{rege}) can be easily derived using a standard technique for parabolic
equations with dynamic boundary conditions (see \cite[Theorem 2.3]{GalNS}
for further details; cf. also \cite{CGGM, GG1, GG2}).
\end{proof}

Consequently, we have shown the following.

\begin{theorem}
\label{lctraj}Under the assumptions of Proposition \ref{propsem}, the weak
attractor $\mathcal{A}_{tr}$, constructed in Theorem \ref{traj}, is a strong
trajectory attractor for the trajectory dynamical system $\left( T_{2}\left(
t\right) ,\mathcal{K}_{+}\right) $. Moreover, $\mathcal{A}_{tr}$ is compact
in $C_{loc}\left( \mathbb{R}_{+};\mathbb{V}^{2-s,2}\right) ,$ for any $s\in
(0,1]$, and the set%
\begin{equation*}
\partial _{t}\mathcal{A}_{tr}:=\left\{ \partial _{t}U:U\in \mathcal{A}%
_{tr}\right\}
\end{equation*}%
is compact in $C_{loc}\left( \mathbb{R}_{+};\mathbb{V}^{1-l,2}\right) ,$ for
any $l\in (0,1/2).$
\end{theorem}

\begin{proof}
It is sufficient to show that the set $T\left( t_{\#}\right) \mathcal{B}$,
where $t_{\#}\geq 1$ is sufficiently large and $\mathcal{B}$ is the
absorbing set (\ref{abs}) for the semigroup $\left\{ T\left( t\right)
\right\} $ is compact in the strong topology of the space $\Theta
_{+}^{w,loc}$. We observe that for every $t\geq t_{\#}$ and any $U\in
\mathcal{K}_{+},$ we have%
\begin{equation}
T\left( t\right) U\in L^{\infty }\left( \mathbb{R}_{+};\mathbb{V}^{2,2}\cap
\mathbb{X}^{\infty }\right) \cap W^{1,\infty }\left( \mathbb{R}_{+};\mathbb{V%
}^{1,2}\right) \cap W^{2,2}\left( \mathbb{R}_{+};\mathbb{X}^{2}\right) .
\label{reg_t}
\end{equation}%
Thus, it is enough to prove that the kernel $\mathcal{K}_{+}$, for the
attractor $\mathcal{A}_{tr}$ defined by (\ref{com}), is compact in the
strong local topology of the space $\Theta _{+}^{b}$. We first note that,
due to estimates (\ref{linfcomp}), (\ref{rege}) and (\ref{reg_t}), we have%
\begin{equation*}
\mathcal{K}_{+}\text{ is bounded in }C_{b}\left( \mathbb{R}_{+};\mathbb{V}%
^{2-s,2}\right) ,\text{ for any }s\in (0,1],
\end{equation*}%
and moreover, it is compact in the local topology%
\begin{equation}
\mathcal{K}_{+}\subset C_{loc}\left( \mathbb{R}_{+};\mathbb{V}%
^{2-s,2}\right) .  \label{comp1}
\end{equation}%
This follows from the following embedding%
\begin{equation*}
\left\{ U:U\in L_{loc}^{\infty }\left( \mathbb{R}_{+};\mathbb{V}%
^{2,2}\right) ,\text{ }\partial _{t}U\in L_{loc}^{\infty }\left( \mathbb{R}%
_{+};\mathbb{V}^{1,2}\right) \right\} \subset C_{loc}\left( \mathbb{R}_{+};%
\mathbb{V}^{2-s,2}\right) .
\end{equation*}%
which is compact. Thus, in view of (\ref{comp1}) and the boundedness of $%
\mathcal{K}_{+}$ in $L^{\infty }\left( \mathbb{R}_{+};\mathbb{X}^{\infty
}\right) $, we immediately see that $\mathcal{K}_{+}$ is also compact in the
(strong)\ local topology of $\Theta _{+}^{b},$ i.e., we have%
\begin{equation*}
\mathcal{K}_{+}\overset{c}{\subset }\Theta _{+,loc}^{b}:=L_{loc}^{\infty
}\left( \mathbb{R}_{+};\mathbb{X}^{2}\right) \cap L_{loc}^{2}\left( \mathbb{R%
}_{+};\mathbb{V}^{1,2}\right) \cap \left( L_{loc}^{r_{1}}\left( \mathbb{R}%
_{+};L^{r_{1}}\left( \Omega \right) \right) \times L_{loc}^{r_{2}}\left(
\mathbb{R}_{+};L^{r_{2}}\left( \Gamma \right) \right) \right) .
\end{equation*}%
The second statement follows analogously using the compactness of the
following embedding%
\begin{equation*}
\left\{ \partial _{t}U:U\in L_{loc}^{\infty }\left( \mathbb{R}_{+};\mathbb{V}%
^{1,2}\right) ,\text{ }\partial _{t}^{2}U\in L_{loc}^{2}\left( \mathbb{R}%
_{+};\mathbb{X}^{2}\right) \right\} \subset C_{loc}\left( \mathbb{R}_{+};%
\mathbb{V}^{1-l,2}\right) .
\end{equation*}%
The proof is finished.
\end{proof}

\begin{corollary}
\label{c1}Under the validity of assumptions of Theorem \ref{lctraj}, the
reaction-diffusion equation (\ref{1.11bb}), with the dynamic boundary
condition (\ref{1.12bb}) possesses a global attractor $\mathcal{A}_{gl}(=%
\mathcal{K}\left( 0\right) ),$ defined by (\ref{proj})-(\ref{proj2}), with $%
\mathcal{A}_{gl}$ bounded in $\mathbb{V}^{2,2}\cap \mathbb{X}^{\infty }$.
\end{corollary}

\begin{remark}
\label{reg}It is worth mentioning that one can also establish more
regularity of the weak solution $u\in \mathcal{C}^{s}\left( \overline{\Omega
}\right) $ as much as it is allowed by the regularity of $\Omega $, $f,$ $g$
and the external sources $h_{1},$ $h_{2}$. Taking advantage of the
regularity result in Corollary \ref{c1}, we can prove that the global
attractor $\mathcal{A}_{gl}$ is finite dimensional, by establishing the
existence of a more refined object called exponential attractor $\mathcal{E}%
_{gl}$. However, since the associated solution semigroup $\mathcal{S}_{t}$
happens to be (uniformly quasi-) differentiable with respect to the initial
data, on the attractor $\mathcal{A}_{gl}$, we can instead employ a volume
contraction argument (see, e.g., \cite{GalNS}).
\end{remark}

We can now extend the results in \cite{GalNS} for the case of nonlinear
boundary conditions, without requiring that the restrictive condition (H2)
holds. More precisely, we can establish the following upper-bound on the
dimension of the global attractor $\mathcal{A}_{gl}$.

\begin{theorem}
\label{main1}Provided that $f,g\in C^{2}\left( \mathbb{R},\mathbb{R}\right) $
satisfy the assumptions of Corollary \ref{c1} for as long as (\ref{weaker2})
holds, the fractal dimension of $\mathcal{A}_{gl}=\mathcal{K}\left( 0\right)
$ admits the estimate%
\begin{equation}
\dim _{F}\mathcal{A}_{gl}\leq \mathcal{C}_{1}\left( 1+\nu ^{-\left(
N-1\right) }\right) ,\text{ for }N\geq 2  \label{updim}
\end{equation}%
and%
\begin{equation}
\dim _{F}\mathcal{A}_{gl}\leq \mathcal{C}_{1}\left( 1+\nu ^{-1/2}\right) ,%
\text{ for }N=1.  \label{updim2}
\end{equation}%
where $\mathcal{C}_{1}$ depends only on $\Omega $, $\Gamma $ and the sources
$h_{1},h_{2}$.
\end{theorem}

\begin{proof}
First, it is easy to establish that the flow $\mathcal{S}_{t}:\mathcal{A}%
_{gl}\rightarrow \mathcal{A}_{gl}$ generated by the reaction-diffusion
equation (\ref{1.11bb}) and dynamic boundary condition (\ref{1.12bb}) is
uniformly differentiable on $\mathcal{A}_{gl},$ with differential%
\begin{equation}
\mathbf{L}\left( t,U\left( t\right) \right) :\Phi =\binom{\xi _{1}}{\xi _{2}}%
\in \mathbb{X}^{2}\mapsto V=\binom{v}{\varphi }\in \mathbb{X}^{2},
\label{4.9}
\end{equation}%
where $V$ is the unique solution to%
\begin{align}
\partial _{t}v& =\nu \Delta v-f^{^{\prime }}\left( u\left( t\right) \right)
v,\text{ }\left( \partial _{t}\varphi +\nu b\partial _{\mathbf{n}%
}v+g^{^{\prime }}(u\left( t\right) _{\mid \Gamma })v\right) _{\mid \Gamma
}=0,  \label{var} \\
V\left( 0\right) & =\Phi .  \notag
\end{align}%
Indeed, the uniform differentiability result follows from the assumptions on
$f,g$ and is a consequence of the boundedness of $\mathcal{A}_{gl}$ into $%
\mathbb{V}^{2,2}\cap \mathbb{X}^{\infty }$ (see \cite{GalNS}). In order to
deduce (\ref{updim})-(\ref{updim2}), it is sufficient (see, e.g., \cite[%
Chapter III, Definition 4.1]{CV})\ to estimate the $j$-trace of the operator%
\begin{equation*}
\mathbf{L}\left( t,U\left( t\right) \right) =\left(
\begin{array}{cc}
\nu \Delta -f^{^{\prime }}\left( u\left( t\right) \right) & 0 \\
-b\nu \partial _{\mathbf{n}} & -g^{^{\prime }}\left( v\left( t\right) \right)%
\end{array}%
\right) .
\end{equation*}%
We have%
\begin{align*}
Trace\left( \mathbf{L}\left( t,U\left( t\right) \right) Q_{m}\right) &
=\sum_{j=1}^{m}\left\langle \mathbf{L}\left( t,U\left( t\right) \right)
\varphi _{j},\varphi _{j}\right\rangle _{\mathbb{X}^{2}} \\
& =\sum_{j=1}^{m}\left\langle \nu \Delta \varphi _{j},\varphi
_{j}\right\rangle _{2}-\sum_{j=1}^{m}\left\langle \nu \partial _{\mathbf{n}%
}\varphi _{j},\varphi _{j}\right\rangle _{2,\Gamma } \\
& -\sum_{j=1}^{m}\left\langle f^{^{\prime }}\left( u\left( t\right) \right)
\varphi _{j},\varphi _{j}\right\rangle _{2}-\sum_{j=1}^{m}\left\langle
g^{^{\prime }}\left( v\left( t\right) \right) \varphi _{j},\varphi
_{j}\right\rangle _{2,\Gamma },
\end{align*}%
where the set of real-valued functions $\varphi _{j}\in \mathbb{X}^{2}\cap
\mathbb{V}^{1,2}$ is an orthonormal basis in $Q_{m}\mathbb{X}^{2}$. By
Theorem \ref{lc}, it follows that every bounded complete trajectory $U\left(
t\right) $, $t\in \mathbb{R}$, for the dynamical system $(\mathcal{S}_{t},%
\mathcal{K}\left( 0\right) )$ is uniformly bounded in $\mathbb{V}^{2,2}\cap
\mathbb{X}^{\infty }$, namely, it holds:%
\begin{equation}
\sup_{t\in \mathbb{R}}\left\Vert U\left( t\right) \right\Vert _{\mathbb{X}%
^{\infty }}\leq \mathcal{C},  \label{compllinf}
\end{equation}%
where the positive constant $\mathcal{C}$ is independent of $U\left(
t\right) ,$ $\nu $, $\Omega ,$ $\Gamma $, but depends on the $L^{\infty }$%
-norms of the external forces $h_{1}$, $h_{2}$ and the constants in (\ref%
{weakly}) (this follows from the usual description of the global attractor,
see, e.g., \cite{BV}, and the main estimates). Thus, exploiting (\ref{arbp}%
)\ once more there holds%
\begin{equation}
\max \sup_{t\in \mathbb{R}}\left\{ \left\Vert f^{^{\prime }}\left( u\left(
t\right) \right) \right\Vert _{L^{\infty }\left( \Omega \right) },\left\Vert
g^{^{\prime }}\left( v\left( t\right) \right) \right\Vert _{L^{\infty
}\left( \Gamma \right) }\right\} \leq C_{\star },  \label{compllinf2}
\end{equation}%
for some $C_{\star }>0$ which depends on $\mathcal{C}$ and the growth rates $%
r_{1},$ $r_{2}\geq 1$; we find%
\begin{equation*}
Trace\left( \mathbf{L}\left( t,U\right) Q_{m}\right) \leq -\nu
\sum_{j=1}^{m}\left\Vert \nabla \varphi _{j}\right\Vert _{2}^{2}+C_{\star }m.
\end{equation*}%
From \cite[Proposition 5.5]{GalNS}, we obtain%
\begin{align*}
Trace\left( \mathbf{L}\left( t,U\right) Q_{m}\right) & \leq -\nu
c_{1}C_{W}\left( \Omega ,\Gamma \right) m^{\frac{1}{N-1}+1}+(c_{1}\nu
C_{W}\left( \Omega ,\Gamma \right) +C_{\star })m \\
& =:\rho \left( m\right) ,
\end{align*}%
for some absolute positive constant $c_{1}$ independent of the parameters of
the problem ($C_{W}$ is given explicitly in \cite[Theorem 5.4]{GalNS}). The
function $\rho \left( y\right) $ is concave. The root of the equation $\rho
\left( y\right) =0$ is%
\begin{equation}
y^{\ast }=\left( 1+\frac{C_{\star }}{\nu c_{1}C_{W}\left( \Omega ,\Gamma
\right) }\right) ^{N-1}.  \label{updim3}
\end{equation}%
Thus, we can apply a well-known result \cite[Chapter VIII, Theorem 3.1]{CV}
to deduce that $\dim _{F}\mathcal{A}_{gl}\leq y^{\ast },$ from which (\ref%
{updim}) follows. The case $N=1$ is similar.
\end{proof}

\begin{remark}
In the competing scenario (\ref{NB}), we can also explicitly estimate from
above the dimension of the global attractor $\mathcal{A}_{gl}=\mathcal{K}%
\left( 0\right) $ for the dynamical system $\left( \mathcal{S}_{t},\mathcal{A%
}_{gl}\right) $ defined by (\ref{proj})-(\ref{proj2}). However, in this case
the upper bound does not seem to be as sharp as in (\ref{updim})-(\ref%
{updim2}). Indeed, by Theorem \ref{lc} (see also Theorem \ref{linff}), it
follows that every bounded complete trajectory $U\left( t\right) $, $t\in
\mathbb{R}$, for the dynamical system $(\mathcal{S}_{t},\mathcal{K}\left(
0\right) )$ satisfies (\ref{compllinf}) with constant $\mathcal{C}=\mathcal{C%
}_{\nu }\sim \nu ^{\sigma _{1}}$ as $\nu \rightarrow 0^{+},$ for some $%
\sigma _{1}<0$; by assumption (H3), this implies that the constant $C_{\star
}$ in (\ref{compllinf2}) behaves as $C_{\star }=C_{\star }\left( \nu \right)
\sim \nu ^{\sigma _{2}}$, for some $\sigma _{2}<-1$ depending on $\sigma
_{1} $ and $r_{1},r_{2}$. Hence, in this case there seems to be a
discrepancy between the upper bound (\ref{updim3}) and the lower bound in (%
\ref{dim}).
\end{remark}

\subsection{A blowup result}

As pointed out at the beginning of this article, nonlinear dissipative
boundary conditions cannot prevent blowup of some solutions of (\ref{1.11bb}%
)-(\ref{1.12bb}) when the non-dissipative interior term $f$ is superlinear,
i.e., when $f$ satisfies (\ref{asss}) for some $r_{1}>2$. We will follow
some arguments similar to ones presented in \cite{RB} for nonlinear Robin
boundary conditions, by constructing some subsolutions to (\ref{1.11bb})-(%
\ref{1.12bb}) which become unbounded in finite time at some points of the
boundary $\Gamma $. We begin with the following notion.

\begin{definition}
\label{sub}A function $v:\overline{\Omega }\times \left( 0,T\right)
\rightarrow \mathbb{R}$ is a subsolution of (\ref{1.11bb})-(\ref{1.12bb}) if
it satisfies%
\begin{equation}
\left\{
\begin{array}{ll}
\partial _{t}v-\nu \Delta v+f\left( v\right) -h_{1}\left( x\right) \leq 0, &
\text{in }\Omega \times (0,T), \\
\partial _{t}v+\nu b\partial _{\mathbf{n}}v+g\left( v\right) -h_{2}\left(
x\right) \leq 0, & \text{on }\Gamma \times \left( 0,T\right)%
\end{array}%
\right.  \label{sub1}
\end{equation}%
and%
\begin{equation}
v\left( 0\right) \leq u_{0}\text{ in }\overline{\Omega }.  \label{sub2}
\end{equation}%
Analogously, the function $v$ is called a supersolution if the inequalities
in (\ref{sub1})-(\ref{sub2}) are reversed.
\end{definition}

From \cite[Section 7]{Be2}, we have the following

\begin{proposition}
\label{P32}Let $u$ be a solution of (\ref{1.11bb})-(\ref{1.12bb}), and let $%
v $ and $\widetilde{v}$ be a subsolution and supersolution, respectively, of
(\ref{1.11bb})-(\ref{1.12bb}), in the sense of Definition \ref{sub}. Then,
\begin{equation*}
v\left( x,t\right) \leq u\left( x,t\right) \leq \widetilde{v}\left(
x,t\right) ,
\end{equation*}%
for all $x\in \overline{\Omega }$ and for as long as they exist. In
particular, if $f\left( 0\right) \leq 0$ and $g\left( 0\right) \leq 0,$ and
if $u_{0}\geq 0,$ then the solution of (\ref{1.11bb})-(\ref{1.12bb})
satisfies $u\left( x,t\right) \geq 0,$ for all $x\in \overline{\Omega }$,
for as long as it exists.
\end{proposition}

We aim to construct subsolutions by comparing solutions of (\ref{1.11bb})-(%
\ref{1.12bb}) with classical solutions fulfilling certain Dirichlet
conditions on the time lateral boundary $\Gamma \times \left( 0,T\right) $.
For that purpose, the following result is very useful (see, \cite{FM, We};
cf. also \cite{Be2, RB}). In what follows, it suffices to consider the case $%
h_{1}\equiv 0,$ $h_{2}\equiv 0$.

\begin{proposition}
\label{P33}Assume that there are a $C^{1}$-concave decreasing function $%
h\left( s\right) $ and a number $s_{0}\geq 0$ such that
\begin{equation*}
\underset{s\rightarrow \infty }{\lim \sup }h^{^{\prime }}\left( s\right)
<-\lambda _{1,D},
\end{equation*}%
where $\lambda _{1,D}>0$ is the first eigenvalue of $-\nu \Delta $ with
Dirichlet boundary conditions. Moreover, suppose that $h\left( s\right) <0$
for all $s>s_{0}$ such that
\begin{equation}
\int_{s_{0}}^{\infty }\frac{ds}{\left\vert h\left( s\right) \right\vert }%
<\infty .  \label{int}
\end{equation}%
Then, there are positive (smooth and locally well-defined) solutions $v$\ of
the reaction-diffusion equation%
\begin{equation}
\begin{array}{ll}
\partial _{t}v-\nu \Delta v+f\left( v\right) =0, & \text{in }\Omega \times
(0,T),%
\end{array}
\label{RD}
\end{equation}%
subject to the boundary and initial conditions%
\begin{equation}
\left\{
\begin{array}{l}
v=0\text{, on }\Gamma \times \left( 0,T\right) , \\
v\left( 0\right) :=v_{0},\text{ in }\Omega ,%
\end{array}%
\right.  \label{D}
\end{equation}%
that blowup in finite time.
\end{proposition}

\begin{theorem}
\label{blowup}Let the assumptions of Proposition \ref{P33} be satisfied, and
assume in addition that $f\left( s\right) \leq h\left( s\right) <0$ for all $%
s>s_{0}$ such that (\ref{int}) holds. Then, for any nonlinear function $g,$
there exist solutions of (\ref{1.11bb})-(\ref{1.12bb}) that blowup in finite
time. Moreover, there exists a positive function $w_{0}\left( x\right) $
such that all solutions with initial data $u_{0}$ greater or equal than $%
w_{0}+v_{0}$, blowup in finite time.
\end{theorem}

\begin{proof}
Let $\varphi _{1}$ be the principal eigenfunction associated with $\lambda
_{1,D}>0$ such that $\left\Vert \varphi _{1}\right\Vert _{L^{1}\left( \Omega
\right) }=1$. It is well-known that $\varphi _{1}>0$ in $\Omega $ by the
maximum principle. Thus, by choosing $A:=\max \left\{ s_{0},s_{0}^{^{\prime
}},0\right\} $ such that $h^{^{\prime }}\left( s\right) <-\lambda _{1,D}$,
for all $s>s_{0}^{^{\prime }},$ we can define $w\left( x\right) :=\delta
\varphi _{1}\left( x\right) >0$, for some $\delta >0,$ such that%
\begin{equation}
\nu b\partial _{\mathbf{n}}w+g\left( A\right) \leq 0\text{ on }\Gamma .
\label{ch1}
\end{equation}%
Note that it is always possible to fix $\delta >0$ since, by the maximum
principle once again, we have $\partial _{\mathbf{n}}\varphi _{1}<0$ on $%
\Gamma $. Let us now define%
\begin{equation*}
\underline{u}\left( x,t\right) =w_{0}\left( x\right) +v\left( x,t\right) ,
\end{equation*}%
where $w_{0}\left( x\right) :=w\left( x\right) +A$ and $v$ is a solution of (%
\ref{RD})-(\ref{D}). Arguing as in the proof of \cite[Theorem 4.3]{RB}, we
can easily establish that%
\begin{equation*}
\partial _{t}\underline{u}-\nu \Delta \underline{u}+f\left( \underline{u}%
\right) \leq 0\text{, in }\Omega \times \left( 0,T\right) .
\end{equation*}%
On the other hand, on $\Gamma \times \left( 0,T\right) $, we have%
\begin{equation*}
\partial _{t}\underline{u}+\nu b\partial _{\mathbf{n}}\underline{u}+g\left(
\underline{u}\right) =\nu b\partial _{\mathbf{n}}w+\nu b\partial _{\mathbf{n}%
}v+g\left( A\right) ,
\end{equation*}%
but since $\partial _{\mathbf{n}}v\leq 0$ on $\Gamma $, by the choice of $A$
and $w$ (cf. (\ref{ch1})), it follows%
\begin{equation*}
\partial _{t}\underline{u}+\nu b\partial _{\mathbf{n}}\underline{u}+g\left(
\underline{u}\right) \leq 0\text{, on }\Gamma \times \left( 0,T\right) .
\end{equation*}%
Consequently, we deduce that $\underline{u}$ is a subsolution of (\ref%
{1.11bb})-(\ref{1.12bb}). Therefore, on account of Propositions \ref{P32}
and \ref{P33}, all solutions $u=u\left( x,t\right) $ of (\ref{1.11bb})-(\ref%
{1.12bb}) that satisfy $u_{0}\geq w_{0}+v_{0}$ in $\overline{\Omega }$,
blowup in finite time. The proof is complete.
\end{proof}

\begin{corollary}
Assume that $f$ satisfies%
\begin{equation*}
\underset{s\rightarrow \infty }{\lim \sup }\frac{f\left( s\right) }{s\left(
\ln s\right) ^{l}}<0,
\end{equation*}%
for some $l>1.$ Then the conclusion of Theorem \ref{blowup} applies.
\end{corollary}

\section{Appendix}

In this section, we state two basic results which are essential to the
analysis of problem (\ref{1.11})-(\ref{1.12}). The first lemma is just a
variation of a result in \cite[Section 5]{GW1} using well-known facts about
nonlinear forms and maximal monotone operators in Sobolev spaces. For the
convenience of the reader, we give below a simple proof of that result.

\begin{lemma}
\label{lem-hemi-2}Let $a\left( s\right) $ satisfy all the assumptions in
Theorem \ref{T1}, and let $U\in \mathbb{V}^{p}$ be fixed. Then the
functional $V\mapsto \mathcal{B}_{p}(U,V),$ defined by (\ref{calc}), belongs
to $\left( \mathbb{V}^{p}\right) ^{\ast }$. Moreover, $\mathcal{B}_{p}$ is
strictly monotone, hemicontinuous and coercive.
\end{lemma}

\begin{proof}
We will make use of \cite[Definitions 2.1 and 2.2]{GW1}. In our case, for $%
p>p_{0},$%
\begin{equation*}
V:=\mathbb{V}^{p}\subset \mathbb{X}^{2}=\left( \mathbb{X}^{2}\right) ^{\ast
}\subset V^{\ast }=\left( \mathbb{V}^{p}\right) ^{\ast }.
\end{equation*}%
Let $U=\binom{u}{u_{\mid \Gamma }}\in \mathbb{V}^{p}$ be fixed. It is clear
that $\mathcal{B}_{p}(U,\cdot )$ is linear. Let $V=\binom{v}{v_{\mid \Gamma }%
}\in \mathbb{V}^{p}$. Exploiting (\ref{1.1.2}), we obtain%
\begin{equation}
\left\vert \mathcal{B}_{p}(U,V)\right\vert \leq \Vert u\Vert
_{W^{1,p}(\Omega )}^{p-1}\Vert V\Vert _{\mathbb{V}^{p}}.  \label{bounded2}
\end{equation}%
This implies that $\mathcal{B}_{p}(U,\cdot )\in \left( \mathbb{V}^{p}\right)
^{\ast },$ for every $U\in \mathbb{V}^{p}$.

Next, let $U,V\in \mathbb{V}^{p}$. Then, recalling (\ref{1.1.7}), we have%
\begin{align}
& \mathcal{B}_{p}(U,U-V)-\mathcal{B}_{p}(V,U-V)  \label{strictt} \\
& =\int_{\Omega }\left( a\left( \left\vert \nabla u\right\vert ^{2}\right)
\nabla u-a\left( \left\vert \nabla v\right\vert ^{2}\right) \nabla v\right)
\cdot \nabla (u-v)dx+\int_{\Omega }\left( |u|^{p-2}u-|v|^{p-2}v\right)
(u-v)dx  \notag \\
& \geq \int_{\Omega }\left( |u|+|v|\right) ^{p-2}|u-v|^{2}dx\geq 0,  \notag
\end{align}%
which shows that $\mathcal{B}_{p}$ is monotone. This estimate also shows that%
\begin{equation*}
\mathcal{B}_{p}(U,U-V)-\mathcal{B}_{p}(V,U-V)>0,
\end{equation*}%
for all $U,V\in \mathbb{V}^{p}$ with $U\neq V$. Thus, $\mathcal{B}_{p}$ is
strictly monotone. The continuity of the norm function implies that $%
\mathcal{B}_{p}$ is hemicontinuous. Finally, it is easy to deduce that%
\begin{equation}
\lim_{\Vert U\Vert _{\mathbb{V}^{p}}\rightarrow +\infty }\frac{\mathcal{B}%
_{p}(U,U)}{\Vert U\Vert _{\mathbb{V}^{p}}}=+\infty ,  \label{coerc}
\end{equation}%
which shows that $\mathcal{B}_{p}$ is coercive. The proof is finished.
\end{proof}

The following lemma is useful in controlling surface integrals by means of
volume integrals. The proof follows from an application of the divergence
theorem to the function $w=\left( \left\vert u\right\vert ^{s}\mathbf{\xi }%
\right) \cdot \mathbf{n,}$ for a smooth vector field $\mathbf{\xi }\in
C^{1}\left( \overline{\Omega },\mathbb{R}^{N}\right) $ such that $\mathbf{%
\xi }\cdot \mathbf{n}=1.$

\begin{lemma}
\label{boundary}Let $p>1$, $s>1$ and $u\in W^{1,p}\left( \Omega \right) .$
Then for every $\varepsilon >0,$ there exists a positive constant $%
C_{\varepsilon }=C\left( \varepsilon ,s,p\right) ,$ independent of $u$, such
that%
\begin{equation*}
\left\Vert u\right\Vert _{L^{s}\left( \Gamma \right) }^{s}\leq \varepsilon
\left\Vert \nabla u\right\Vert _{L^{p}\left( \Omega \right)
}^{p}+C_{\varepsilon }\left( \left\Vert u\right\Vert _{L^{\gamma }\left(
\Omega \right) }^{\gamma }+1\right) ,
\end{equation*}%
where $\gamma =\max \left( s,p\left( s-1\right) /\left( p-1\right) \right) .$
\end{lemma}

\begin{acknowledgement}
The author wishes to express his thanks to Martin Meyries for bringing his
attention to paper \cite{RB}, and for all the interesting and fruitful
discussions on the issues presented in this article.
\end{acknowledgement}

\end{document}